\documentclass[10pt]{article}

\usepackage{amsfonts}
\usepackage{amsmath}
\usepackage{amssymb}
\usepackage{amsbsy}
\usepackage{amsopn}
\usepackage{amstext}
\usepackage{bm}
\usepackage{latexsym}
\usepackage{algorithm,algorithmic}
\usepackage[colorlinks=true]{hyperref}
\usepackage{graphicx}
\usepackage{color}
\usepackage{epstopdf}
\graphicspath{{./images/}}
\usepackage{subfig}

\usepackage{threeparttable}

\usepackage{verbatim}
\usepackage{indentfirst}

\usepackage{cite}
\usepackage{hyperref}

\usepackage[notref,notcite]{showkeys}

\usepackage{cancel}
\usepackage{float}
\usepackage{subfig}
\usepackage[section]{placeins}

\newtheorem{remark}{Remark}[section]

\numberwithin{equation}{section}
\numberwithin{lemma}{section}	
\numberwithin{Defi}{section}
\numberwithin{theorem}{section}
\numberwithin{remark}{section}
\numberwithin{corollary}{section}
\numberwithin{figure}{section}
\numberwithin{table}{section}

\def\NN{\hbox{\rlap{I}\kern.16em N}}
\def\NC{\hbox{\rlap{\kern.24em\raise.1ex\hbox{\vrule height1.3ex width.9pt}}C}}

\let\div\relax
\DeclareMathOperator{\div}{div}

\let\curl\relax
\DeclareMathOperator{\curl}{curl}

\title{FieldTNN-based machine learning method for Maxwell eigenvalue problems}
\author{Jiantao Jiang\thanks{Beijing Computational Science Research Center, Beijing 100193, China({\tt jiangjiantao@csrc.ac.cn}).}
\and Yanli Wang\thanks{Beijing Computational Science Research Center, Beijing 100193, China({\tt ylwang@csrc.ac.cn}).}
\and Yifan Wang\thanks{School of Mathematical Sciences, Peking University, Beijing 100871, China({\tt wangyifan1994@pku.edu.cn}).}
\and Hehu Xie\thanks{LSEC, NCMIS, Institute of Computational Mathematics, Academy of Mathematics and Systems Science, Chinese Academy of Sciences, Beijing 100190, China, and School of Mathematical Sciences, University of Chinese Academy of Sciences, Beijing 100049, China({\tt hhxie@lsec.cc.ac.cn}).}}
\date{\today}

\begin{document}
\maketitle
\tableofcontents
\newpage

\begin{abstract}
The aim of this paper is to introduce a FieldTNN-based machine learning method for solving the Maxwell eigenvalue problem in both 2D and 3D domains, including both tensor and non-tensor computational regions. First, we extend the existing TNN-based approach to address the Maxwell eigenvalue problem, a fundamental challenge in electromagnetic field theory. Second, we tackle non-tensor computational domains, which represents a novel and significant contribution of this work. Third, we incorporate the divergence-free condition into the optimization process, allowing for the automatic filtering of spurious eigenpairs. Numerical examples are presented to demonstrate the efficiency and accuracy of our algorithm, underscoring its potential for broader applications in computational electromagnetics.

	

\end{abstract}
\section{Introduction}
In this paper, we develop a high-precision machine learning method for solving the Maxwell eigenvalue problem, a fundamental model problem in computational electromagnetism, with applications in areas such as electromagnetic waveguides and cavity resonances \cite{Costabel2003,FemEig2016}. 
Denote by $\nabla\times$ and $\nabla\cdot$ the $\curl$ and $\div$ operators, respectively, we consider the following Maxwell eigenvalue problems: Find eigenpairs $(\lambda,\mathbf{E})$ such that
\begin{equation}\label{eq:Maxwelleig}
\begin{aligned}
\nabla\times\left(\mu^{-1}\nabla\times\mathbf{E}\right)&=\lambda\varepsilon\mathbf{E},&\text{in}~\Omega,\\
\nabla\cdot\left(\varepsilon\mathbf{E}\right)&=0,&\text{in}~\Omega,\\
\mathbf{E}\times\mathbf{n}&=0,&\text{on}~\partial\Omega,
\end{aligned}    
\end{equation}
where $\Omega\subset\mathbb{R}^{d},~d = 2,3$, $\mathbf{n}$ is the unit outward normal to the boundary $\Omega$, $\mu$ and $\varepsilon$ denote the magnetic permeability and dielectric permittivity, respectively. 
The electric field $\mathbf{E}$ are called the Maxwell eigenfunctions, and the values of $\lambda$ are referred to resonant frequencies of $\Omega$ or Maxwell eigenvalues.

The computation of Maxwell eigenproblem has long presented a compelling challenge for the numerical analysis community. 
In addition to its broad and critical practical applications, ranging from signal processing and cardiovascular and neural biology to nuclear fusion—the Maxwell eigenvalue problem is distinguished by mathematical features that markedly differentiate it from traditional elliptic eigenvalue problems \cite{Costabel2003}. 
The classical numerical methods, including finite element method (see e.g., \cite{Hiptmair_2002,MonkMaxwell} and references there in) and spectral method (see e.g., \cite{IEEE_SEM,CICP_2015,MMAMaxwell} and the references there in), for solving Maxwell's equations have garnered significant attention. 
However, despite considerable progress, the Maxwell eigenvalue problem remains particularly challenging due to its inherent complexity, such as the vector nature of the fields, the presence of spurious modes, and the requirement to enforce divergence-free conditions. 
These challenges become more pronounced in high-frequency regimes and in complex geometries, where standard methods may struggle to achieve the desired accuracy or computational efficiency.

In recent decades, various machine learning-based algorithms have been proposed to solve eigenvalue problems of partial differential equations (PDEs). Lagaris et al.  pioneered the application of neural networks to solve PDEs by constructing networks that satisfy boundary conditions, successfully addressing simple PDEs eigenvalue problems \cite{LAGARIS19971}. 
The Deep Ritz method, proposed by E and Yu, also demonstrated effectiveness in solving eigenvalue problems by minimizing the variational form of the energy functional to obtain approximation solutions \cite{DRM2018}. 
Han et al. introduced a deep neural network approach to solve high-dimensional eigenvalue problems by reformulating them as fixed point problems using stochastic differential equations \cite{JCPNNeig2020}. 
Li et al. proposed a neural-network-based solver for the static Schr\"{o}dinger equation, effectively addressing the challenge of obtaining multiple excited-state energy eigenvalues and wave functions in multivariable quantum systems \cite{PRANNeig2021}. 
Due to the mesh-free nature of machine learning, it becomes feasible to compute partial differential equations in complex domains, particularly for the Maxwell eigenvalue problem in such regions, which bears substantial engineering significance and practical applications. 
Consequently, employing machine learning to solve the Maxwell eigenvalue problem should be a focal point of research. 
Unfortunately, to the best of our knowledge, there have been few studies addressing this problem using machine learning to date.

Recently, a type of Tensor Neural Network (TNN) and corresponding machine learning methods for solving high-dimensional PDEs \cite{JCM_TNN,JCP_TNN,Posteriori_TNN,Interpolation_TNN} have been developed. 
In these works, the trial and test functions in Galerkin scheme are constructed by the TNN, which can be viewed as basis functions generated by a neural network. 
The TNN has demonstrated its efficiency and accuracy in solving various Schr\"{o}dinger equations and high dimensional boundary value problems and eigenvalue  problems. 
The structure of TNN, its approximation properties, and the computational complexity of the related integration have been thoroughly discussed in \cite{JCM_TNN}.
Additionally, the construction of the TNN method for solving multiple eigenpairs in high-dimensional eigenvalue problems has been proposed in \cite{JCP_TNN}. 
A novel machine learning method that combines tensor neural networks with a posteriori error estimators to solve high-dimensional boundary value problems and eigenvalue problems with high accuracy and efficiency is introduced in \cite{Posteriori_TNN}. 
Li et al. introduced a tensor neural network-based interpolation method for approximating high-dimensional functions without tensor product structure, enabling the design of efficient and accurate numerical methods for high-dimensional integration and solving partial differential equations \cite{Interpolation_TNN}.
However, research on applying the TNN method to vector field equations, such as Maxwell's equations or Stokes equations, has been limited so far since the constraint of divergence free. 
Therefore, in this paper, we propose a Field Tensor Neural Network (FieldTNN) and corresponding machine-learning methods for solving the Maxwell eigenvalue problems in $2$D and $3$D domains. 
Furthermore, previous works on TNN have primarily focused on tensor-type domains, including those that can be transformed into tensor-type domains through coordinate transformations. 
In this work, we extend the computational domain to include non-tensor domains, particularly for the singular domains.

Building on the aforementioned considerations, we extend the TNN-based machine learning method to the Maxwell eigenvalue problem in various domains, including square and cube cavities, 2D and 3D L-shaped domains and inhomogeneous media within square cavity. 
Thus, we are proceed to construct a new type of TNN: field tensor neural network (FieldTNN), which can solve the PDEs in vector field. 
For the Maxwell eigenvalue problem in complex domains, such as L-shaped cavities, we adopt the idea of domain decomposition and spectral element method \cite{SEM_JCP}, we utilize FieldTNN-based machine learning method to train compact support basis functions in each subdomain. 
A key challenge in the Maxwell eigenvalue problem is handling the divergence-free condition. If this condition is not integrated into the numerical scheme or basis functions, spurious eigenpairs may arise. 
To address this issue, a penalty term in the loss function is considered, ensuring that spurious eigenpairs yield a large loss value, while the loss for ``real'' eigenpairs remains small. 
Thus, the divergence-free constraint is naturally embedded into the minimization problem.

The organization of this paper is as follows. 
In Section \ref{sec:Pre}, we present the preliminaries, focusing on the mathematical framework relevant to Maxwell eigenvalue problems. 
The architecture of the FieldTNN is introduced in Section \ref{sec:FieldTNN}, with a focus on its structural design and the procedure for computing multiple eigenpairs. 
Section \ref{sec:nontensor} extends the discussion to TNN-based machine learning methods applied to non-tensor domains and provides a detailed analysis of both TNN and FieldTNN. 
The application of FieldTNN to solve Maxwell eigenvalue problems is detailed in Section \ref{sec:Maxwell}, including the quadrature scheme for the variational form, the implementation of divergence-free and boundary conditions, and the algorithmic framework. 
Section \ref{sec:Numerical} presents numerical examples, demonstrating the efficiency of the proposed method across various cavity configurations, including square, L-shaped, inhomogeneous, and cubic cavities in both 2D and 3D. 
Finally, Section \ref{sec:conclusion} concludes the paper with a summary of our findings and potential future research.
\section{Preliminaries}\label{sec:Pre}
In this section, we will introduce some notations, definitions and fundamental concepts related to Sobolev spaces, which is essential for the analysis of the Maxwell eigenvalue problems. 
Let $\Omega \subset \mathbb{R}^{d}$ be a bounded, simply connected Lipschitz  domain, where $d$ denotes the dimension and $\mathbf{x} = (x_{1}, x_{2})$ in 2D or $\mathbf{x} = (x_{1}, x_{2}, x_{3})$ in 3D. The mathematical treatment of electromagnetic fields requires a rigorous framework provided by Sobolev spaces, particularly for ensuring the well-posedness of variational formulations.

The space $L^2(\Omega)$ consists of all measurable functions $v(x): \Omega \rightarrow \mathbb{R}$ that are square-integrable:
\begin{align*}
&L^2(\Omega):= \left\{v(x): \int_{\Omega}v^{2}\mathrm{d}x < \infty \right\},
\end{align*}
where the $L^{2}$-inner product is defined as $(u, v) = \int_{\Omega} u v\mathrm{d}x$, 
and the associated norm is denoted by $\|u\| = \left(\int_{\Omega} u^{2} \mathrm{d}x\right)^{1/2}$.
The Sobolev space $\mathrm{H}^s(\Omega)$, defined on $\Omega$, consists of functions whose weak derivatives up to order $s$ are square-integrable. The inner product in $\mathrm{H}^s(\Omega)$ is defined as $(u, v)_{\mathrm{H}^s(\Omega)} = \sum\limits_{|\alpha| = 0}^{s} \int_{\Omega} D^\alpha u(x) D^\alpha v(x)\mathrm{d}x$,
where $D^\alpha$ represents the weak derivative of the function and $\alpha$ is a multi-index. The corresponding norm is given by $\|u\|_{\mathrm{H}^s(\Omega)} = \bigg( \sum\limits_{|\alpha| = 0}^{s} \|D^\alpha u\|^2 \bigg)^{1/2}$. 
For the case $s = 0$, the Sobolev space $\mathrm{H}^0(\Omega)$ coincides with the space $L^2(\Omega)$. 

Based on these, we define the space $\mathbf{H}(\div; \Omega)$ in the form:
\begin{align*}
\mathbf{H}(\div; \Omega) : = \left\{\mathbf{v} \in L^2(\Omega)^{d} \mid \nabla \cdot \mathbf{v} \in L^2(\Omega)\right\}
\end{align*}
equipped with the inner product $(\mathbf{u}, \mathbf{v})_{\mathbf{H}(\div; \Omega)} = (\mathbf{u}, \mathbf{v}) + (\nabla \cdot \mathbf{u}, \nabla \cdot \mathbf{v})$
and the associated norm $\|\mathbf{u}\|_{\mathbf{H}(\div; \Omega)} = \Big((\mathbf{u}, \mathbf{u})_{\mathbf{H}(\div; \Omega)}\Big)^{1/2}$.
The space $\mathbf{H}(\curl; \Omega)$ is defined as
\begin{align*}
&\mathbf{H}(\curl;\Omega):=\Big\{\mathbf{v} \in L^2(\Omega)^{d} \mid \nabla \times \mathbf{v} \in L^2(\Omega)^{d}\Big\}
\end{align*}
equipped with the inner product $(\mathbf{u}, \mathbf{v})_{\mathbf{H}(\curl; \Omega)}=(\mathbf{u}, \mathbf{v})+(\nabla \times \mathbf{u}, \nabla \times \mathbf{v})$.
Correspondingly, this inner product leads to the norm $\|\mathbf{u}\|_{\mathbf{H}(\curl;\Omega)} = \Big((\mathbf{u}, \mathbf{u})_{\mathbf{H}(\operatorname{curl} ; \Omega)}\Big)^{1/2}$.
To address perfect electric conductor (PEC) boundary condition in \eqref{eq:Maxwelleig}, we introduce a subspace $\mathbf{H}_0(\curl; \Omega)$ of $\mathbf{H}(\curl; \Omega)$:
\begin{align*}
&\mathbf{H}_0(\curl;\Omega):=\Big\{\mathbf{v} \in \mathbf{H}(\curl;\Omega) \mid \mathbf{v} \times \mathbf{n} = 0 \ \text {on}\  \partial \Omega\Big\}.
\end{align*}
This subspace imposes a homogeneous tangential boundary condition, which is significant in the study of conducting boundaries in electromagnetism. 
In the context of Maxwell eigenproblems, the divergence-free condition for the electric field $\mathbf{E}$, given by $\nabla \cdot (\varepsilon \mathbf{E}) = 0$, must be satisfied. To this end, we define another subspace of $\mathbf{H}(\curl;\Omega)$: 
\begin{align*} 
&\mathbf{Y}(\Omega):= \left\{\mathbf{v} \in \mathbf{H}_{0}(\curl; \Omega) \mid \nabla \cdot (\varepsilon \mathbf{v}) = 0 \right\}, 
\end{align*} 
which incorporates the divergence-free condition into the function space itself. 
Then the variational formulation of the Maxwell eigenvalue problem \eqref{eq:Maxwelleig} is: Find non-trivial $\mathbf{E} \in \mathbf{Y}(\Omega)$ and corresponding eigenvalue $0\neq \lambda  \in \mathbb{R}$ such that 
\begin{align}
&\mathcal{A}(\mathbf{E}, \mathbf{V}) = \lambda \mathcal{B}(\mathbf{E}, \mathbf{V}),~\forall~\mathbf{V} \in \mathbf{Y}(\Omega),\label{eq:vform1}
\end{align} 
where the bilinear forms $\mathcal{A}(\cdot, \cdot)$ and $\mathcal{B}(\cdot, \cdot)$ are defined as 
\begin{align*}
&\mathcal{A}(\mathbf{E}, \mathbf{V}) = \left(\mu^{-1} \nabla \times \mathbf{E},~\nabla \times \mathbf{V} \right),~   \mathcal{B}(\mathbf{E}, \mathbf{V}) = \left(\varepsilon\mathbf{E}, \mathbf{V}\right).
\end{align*}
Here, $\mathbf{E} = (\mathrm{E}_1, \mathrm{E}_2, \mathrm{E}_3)$ and $\mathbf{V} = (\mathrm{V}_1, \mathrm{V}_2, \mathrm{V}_3)$ represent the trial and test function, respectively. Due to the difficulties in imposing the divergence-free condition directly, we alternatively consider a relaxed variational formulation: Find $ 0\neq\mathbf{E}  \in \mathbf{H}_0(\operatorname{curl}; \Omega)$ and $0\neq \lambda \in \mathbb{R}$ such that 
\begin{align}
&\mathcal{A}(\mathbf{E}, \mathbf{V}) = \lambda \mathcal{B}(\mathbf{E}, \mathbf{V}),~\forall~\mathbf{V} \in \mathbf{H}_{0}(\operatorname{curl} ; \Omega).\label{eq:vform2}
\end{align}

Let $\left(\mathbf{E}_{h}, \lambda_{h}\right)$ be a discrete eigenpair obtained by discretizing the eigenvalue problem \eqref{eq:vform2} using a certain numerical method. 
Upon solving this discrete problem, both ``real'' and spurious eigenpairs are obtained. 
As mentioned in \cite{CMAME2008}, the occurrence of spurious eigenpairs arises from solutions to \eqref{eq:vform2} that are not divergence-free. 
Interestingly, when a spurious eigenvalue is simple, the corresponding eigenfunction is inherently curl-free, as discussed in \cite{Costabel2003,CMAME2008}. This property suggests a potential strategy for developing a training method to filter out spurious solutions in machine learning, which will be introduced in Section \ref{subsec:div_free}.

In the following sections, we build on these foundational concepts to develop and analyze machine learning methods tailored for solving Maxwell eigenvalue problems efficiently, particularly in the context of complex geometries and material properties.
\section{FieldTNN-based machine learning method}\label{sec:FieldTNN}
The TNN architecture is distinguished by its low-rank structure, which is constructed by the tensor product of $d$ subnetworks, each handling one-dimensional inputs and multiple outputs. 
This structure enables the design of an efficient and precise quadrature scheme for high-dimensional integration associated with TNN, including the computation of inner products between two TNNs. 
In the previous work \cite{JCM_TNN}, a detailed introduction for the TNN framework and a numerical integration method was provided, demonstrating that its computational complexity grows polynomially with respect to the dimensionality when solving partial differential equations. 
Based on that, we propose FieldTNN in this paper, which can solve vector field equations, including Maxwell's equations as well as the associated eigenvalue problems and so on. 

\subsection{FieldTNN architecture}\label{subsec:FieldTNN}
\begin{figure}[htb]
\centering
\includegraphics[width=.7\textwidth]{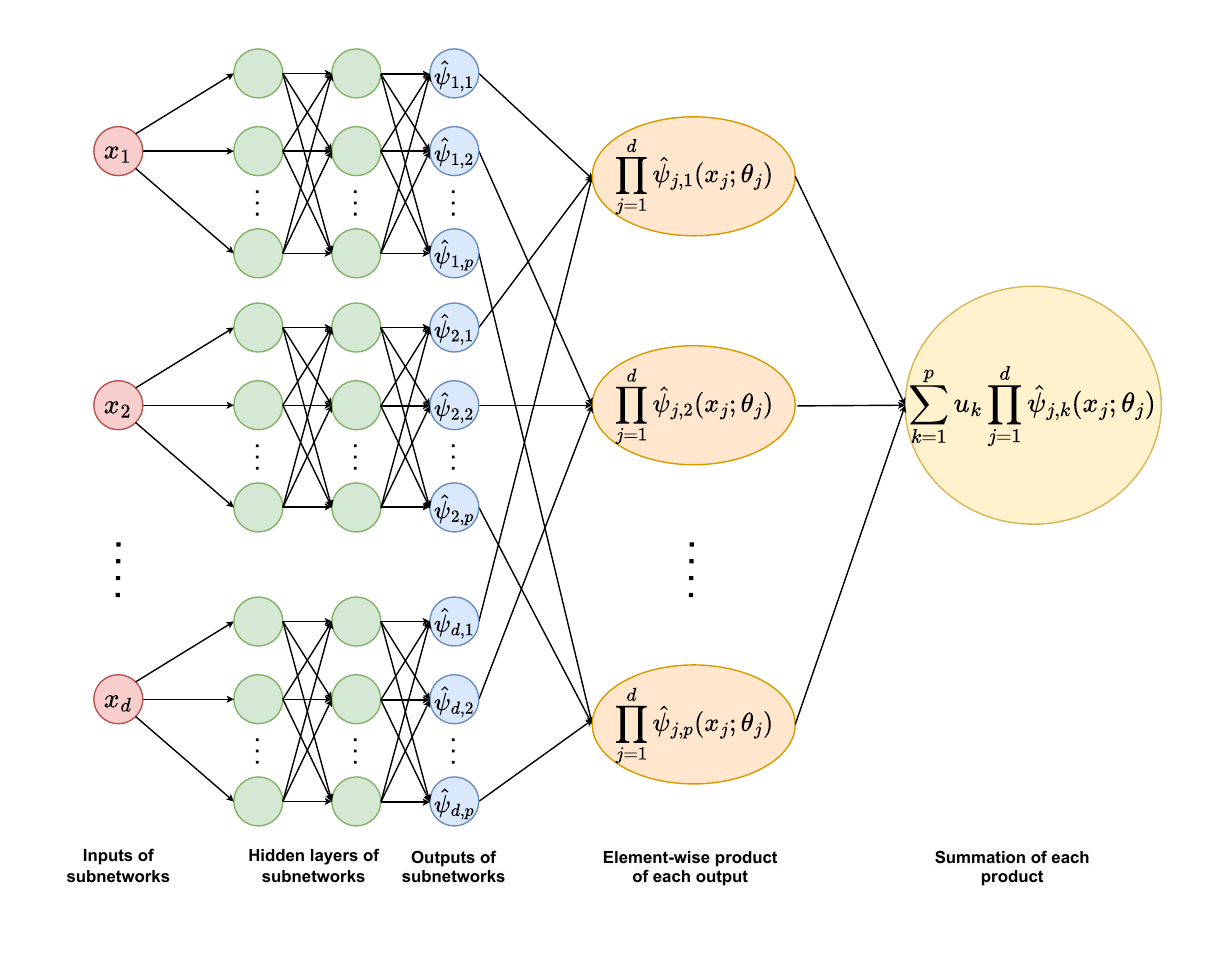}~
\caption{Architecture of TNN. The red circles $x_{j}~(j=1,\cdots,d)$ are the input data, then the blue circles $\hat{\psi}_{j,k}(x_{j};\theta_{j})~(j=1,\cdots,d,~k=1,\cdots,p)$ are obtain by the hidden layers of each subnetworks (the green cirlces), i.e., the FNN as shown in Eq. \eqref{eq:FNN}. The orange ovals $\prod\limits_{j=1}^{d} \hat{\psi}_{j, k}\left(x_j ; \theta_{j}\right)~(k=1,\cdots,p)$ are obtained by the scalar multiplication of the blue circles. The summation of $p$ oranges ovals form the final output of FieldTNN.}
\label{fig:TNN}
\end{figure}

In this paper, denote by $W_{j}^{(\ell)} \in \mathbb{R}^{N_{\ell} \times N_{\ell-1}}$ and $b_{j}^{(\ell)} \in \mathbb{R}^{N_{\ell}}$ the weights and biases in the $\ell$-th layer of $j$-th subnetwork, 
where $N_{\ell}$ indicate the number of neurons in the $\ell$-th layer, $N_{0} = 1$, $N_{L} = p$ denotes the rank, and $\ell =1,\cdots,L$, $j = 1,\cdots,d$. 
Then, we write $W_{j} = \big\{W_{j}^{(1)},\cdots, W_{j}^{(L)}\big\}$ and $b_{j} = \big\{b_{j}^{(1)},\cdots,b_{j}^{(L)}\big\}$.
Denote by 
\begin{align*}
&\bigodot\limits_{i = 1}^{L} g^{i} = g^{L}\circ g^{L-1}\circ\cdots g^{2}\circ g^{1},
\end{align*}
the consecutive composition of $g^{i}$, $i =1,2,\cdots,L$.
Subsequently, we introduce a $L$ layers fully connected neural network (FNN) architectures of $j$-th subnetwork as follows:
\begin{align}
\psi_{j}(x_{j}; \theta_{j}) &= \left(W_{j}^{(L)} \bigodot\limits_{i = 1}^{L-1}\sigma\left(W_{j}^{(i)}\cdot + b_{j}^{(i)}\right)  + b_{j}^{(L)}\right)\left(x_{j}\right),\label{eq:FNN}
\end{align}
where $j=1,\cdots,d$, $\sigma(\cdot)$ denote activation function and $\theta_{j} :=\{W_{j}, b_{j}\}$.
Then the TNN is defined as
\begin{align}
&\Psi(\mathbf{x}; \Theta) = \sum_{k=1}^p u_{k} \hat{\psi}_{1, k}\left(x_1 ; \theta_{1}\right) \hat{\psi}_{2, k}\left(x_2 ; \theta_{2}\right) \cdots \hat{\psi}_{d, k}\left(x_d ; \theta_{d}\right) = \sum_{k=1}^p u_{k} \prod_{j=1}^d \hat{\psi}_{j, k}\left(x_j ; \theta_{j}\right), \label{eq:TNN}
\end{align}
where $\mathbf{x} = \left(x_{1}, \cdots, x_{d}\right)$, $\mathbf{u} = \left\{u_{k}\right\}_{k=1}^p$, and $\Theta=\left\{\mathbf{u}, \theta_{1}, \cdots, \theta_{d}\right\}$ denotes the parameters of whole neural network architecture. The normalized functions $\hat{\psi}_{j, k}\left(x_{j}; \theta_{j}\right)$ are defined as follows:
\begin{align}
&\hat{\psi}_{j, k}\left(x_j;\theta_{j}\right)=\frac{\psi_{j, k}\left( x_{j}; \theta_{j} \right)}{\left\|\psi_{j, k}\left(x_{j};\theta_{j}\right)\right\|},\label{eq:normal}
\end{align}
where $j = 1,\cdots,d$, $k=1, \cdots, p$.
As mentioned in \cite{JCP_TNN} that the TNN architecture \eqref{eq:TNN} has better numerical stability during the training process than the one defined in \cite{JCM_TNN}, though they are mathematically equivalent.
For providing readers with an intuitive understanding of TNN, we present the corresponding network architecture of TNN in Figure \ref{fig:TNN}.
As shown in Figure \ref{fig:TNN}, we observe that the parameters of each rank of TNN exhibit a correlation with the FNN, thereby ensuring the stability of TNN-based machine learning method.
Moreover, we also use $\hat{\psi}_{j}(x_{j}; \theta_{j})$ to denote the output of the $j$-th subnetwork, i.e.,
\begin{align}
\hat{\psi}_{j}(x_{j}; \theta_{j}):= \left(\hat{\psi}_{j,1}(x_{j}; \theta_{j}),\cdots,\hat{\psi}_{j,p}(x_{j}; \theta_{j})\right)^{\top},\ j=1,2,\cdots,d.\label{eq:TNNjth}
\end{align}

\begin{figure}[!h]
\centering
\includegraphics[width=.9\textwidth]{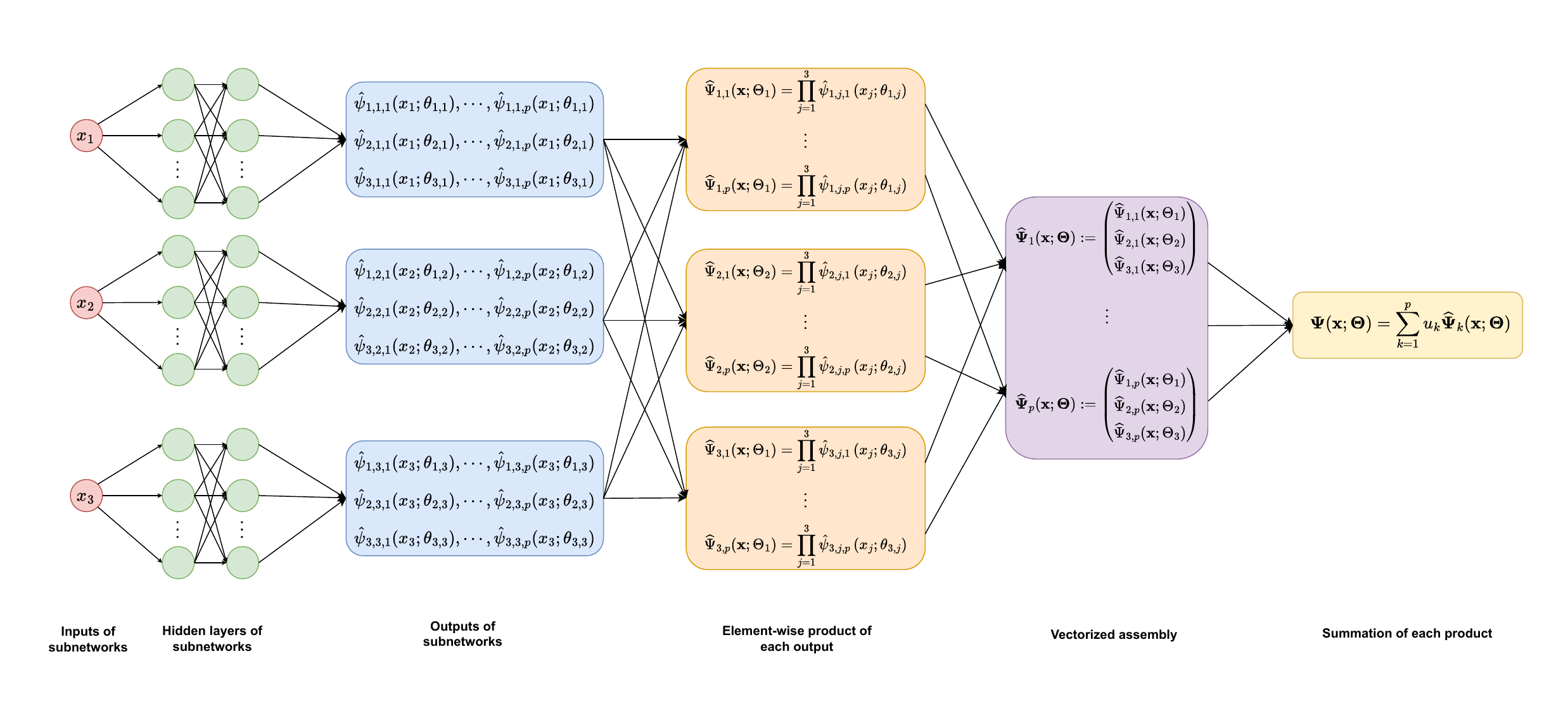}~
\caption{Architecture of FieldTNN. The red circles $x_{j}~(j=1,\cdots,d)$ are the input data, then the blue squircles $\hat{\psi}_{i,j,k}(x_{j};\theta_{i,j})~(i,j=1,\cdots,d,~k=1,\cdots,p)$ are obtain by the hidden layers each subnetworks (the green cirlces), i.e., the FNN as shown in Eq. \eqref{eq:FNN}, where $N_{L} = dp$. The orange squircles $\widehat{\Psi}_{i,k}(\mathbf{x}; \Theta_{i})~(i=1,\cdots,d,~k=1,\cdots,p)$ are obtained by the scalar multiplication of the blue ones. Then the purple squircles are derived by vectorized assembly from the orange ones. The summation of $p$ elements in pruple squircles form the final output of TNN.}
\label{fig:FieldTNN}
\end{figure}

For solving the Maxwell eigenvalue problems, we need a new type of tensor neural network, field tensor neural network (FieldTNN), to represent the appropriate vector functions in electromagnetic fields.
For example, a three-dimensional FieldTNN is shown in Figure \ref{fig:FieldTNN}.
To be specific,
we first define three-dimensional rank-one scalar functions $\widehat{\Psi}_{i,k}(\mathbf{x}; \Theta_{i})$ as follows
\begin{align}
&\widehat{\Psi}_{i,k}(\mathbf{x}; \Theta_{i}) = \hat{\psi}_{i,1, k}\left(x_{1} ; \theta_{i,1}\right) \hat{\psi}_{i,2, k}\left(x_{2} ; \theta_{i,2}\right) \hat{\psi}_{i,3, k}\left(x_{3} ; \theta_{i,3}\right) = \prod_{j = 1}^{3}\hat{\psi}_{i,j, k}\left(x_{j} ; \theta_{i,j}\right),\label{eq:FieldTNNik}
\end{align}
where $i = 1,2,3$, $k = 1,\cdots,p$, and $\Theta_{i} = \left\{\alpha, \theta_{i,1}, \theta_{i,2}, \theta_{i,3}\right\}$.
Then, we use $\mathbf{\widehat{\Psi}}_{k}(\mathbf{x} ; \mathbf{\Theta})$ to denote three-dimensional vector functions which are collated by the scalar functions \eqref{eq:FieldTNNik} as follows
\begin{equation}\label{eq:FieldTNNk}
\mathbf{\widehat{\Psi}}_{k}(\mathbf{x} ; \mathbf{\Theta}) :=
\left(
\begin{array}{l}
\widehat{\Psi}_{1,k}(\mathbf{x}; \Theta_{1})\\
\widehat{\Psi}_{2,k}(\mathbf{x}; \Theta_{2})\\
\widehat{\Psi}_{3,k}(\mathbf{x}; \Theta_{3})
\end{array}
\right)
=
\left(\begin{array}{l}
\hat{\psi}_{1,1,k}\left(x_{1} ; \theta_{1,1}\right) \hat{\psi}_{1,2,k}\left(x_{2} ; \theta_{1,2}\right) \hat{\psi}_{1,3,k}\left(x_{3} ; \theta_{1,3}\right) \\
\hat{\psi}_{2,1,k}\left(x_{1} ; \theta_{2,1}\right) \hat{\psi}_{2,2,k}\left(x_{2} ; \theta_{2,2}\right) \hat{\psi}_{2,3,k}\left(x_{3} ; \theta_{2,3}\right) \\
\hat{\psi}_{3,1,k}\left(x_{1} ; \theta_{3,1}\right) \hat{\psi}_{3,2,k}\left(x_{2} ; \theta_{3,2}\right) \hat{\psi}_{3,3,k}\left(x_{3} ; \theta_{3,3}\right)
\end{array}\right),
\end{equation}
where $\mathbf{\Theta} := \left\{\Theta_{1},\Theta_{2},\Theta_{3}\right\}$. 
Finally, we can construct FieldTNN in the following form, 	
\begin{align}\label{eq:FieldTNN}
& \mathbf{\Psi}(\mathbf{x} ; \mathbf{\Theta}) = \sum_{k = 1}^{p} u_{k} \mathbf{\widehat{\Psi}}_{k}(\mathbf{x} ; \mathbf{\Theta}).
\end{align}
Under the above definition, it is easy to represent differential operators related to FieldTNN, especially the divergence operator and the curl operator that is needed for solving the Maxwell eigenvalue problem.
More specifically,
we can write $\nabla \cdot \mathbf{\Psi}(\mathbf{x} ; \mathbf{\Theta})$  in the following form,
\begin{align*}
\nabla \cdot \mathbf{\Psi}(\mathbf{x} ; \mathbf{\Theta}) = & \sum_{k = 1}^{p}u_{k}\left(\frac{\partial \widehat{\Psi}_{1,k}}{\partial x_{1}} + \frac{\partial \widehat{\Psi}_{2,k}}{\partial x_{2}} + \frac{\partial \widehat{\Psi}_{3,k}}{\partial x_{3}}\right)\\
= & \sum_{k = 1}^{p} u_{k}\left(\frac{\partial \hat{\psi}_{1,1,k}\left(x_{1} ; \theta_{1,1}\right)}{\partial x_{1}} \hat{\psi}_{1,2,k}\left(x_{2} ; \theta_{1,2}\right) \hat{\psi}_{1,3,k}\left(x_{3} ; \theta_{1,3}\right)\right. \\
&\qquad \qquad + \hat{\psi}_{2,1, k}\left(x_{1} ; \theta_{2,1}\right) \frac{\partial \hat{\psi}_{2,2, k}\left(x_{2} ; \theta_{2,2}\right)}{\partial x_{2}}  \hat{\psi}_{2,3, k}\left(x_{3} ; \theta_{2,3} \right)\\
&\qquad \qquad + \left.\hat{\psi}_{3,1, k}\left(x_{1} ; \theta_{3,1}\right) \hat{\psi}_{3,2, k}\left(x_{2} ; \theta_{3,2}\right) \frac{\partial \hat{\psi}_{3,3, k}\left(x_{3} ; \theta_{3,3} \right)}{\partial x_{3}}\right),
\end{align*}
and $\nabla \times \mathbf{\Psi}(\mathbf{x}; \Theta_{i}) = \sum \limits_{k = 1}^{p}u_{k}\left(\frac{\partial \widehat{\Psi}_{3,k}}{\partial x_{2}}-\frac{\partial \widehat{\Psi}_{2,k}}{\partial x_{3}},~\frac{\partial \widehat{\Psi}_{1,k}}{\partial x_{3}}-\frac{\partial \widehat{\Psi}_{3,k}}{\partial x_{1}},~\frac{\partial \widehat{\Psi}_{2,k}}{\partial x_{1}}-\frac{\partial \widehat{\Psi}_{1,k}}{\partial x_{2}}\right)$ can be derived in the same manner.

\subsection{Computing multi-eigenpairs by FieldTNN}\label{subsec:multi}
In this section, we provide a brief introduction to a FieldTNN-based discretization method for solving multi-eigenpairs of eigenvalue problems of vector field equations. 
In order to discuss our method within an abstract framework, we first introduce some necessary notations. Supposed that Sobolev spaces
$\mathbb{V} \subset \mathbb{W}$
are two Hilbert spaces, denote by $\mathcal{A}(\cdot, \cdot)$ and $\mathcal{B}(\cdot, \cdot)$ the two symmetric positive definite
bilinear forms defined on $\mathbb{V} \times \mathbb{V}$ and $\mathbb{W} \times \mathbb{W}$,
respectively. Subsequently, we define induced norms of the spaces $\mathbb{V}$ and $\mathbb{W}$ as follows,
\begin{align*}
\|\mathbf{v}\|_{\mathcal{A}}^{2} = \mathcal{A}(\mathbf{v}, \mathbf{v}),~\forall~\mathbf{v} \in \mathbb{V},~\text{and}~
\|\mathbf{w}\|_{\mathcal{B}}^{2} = \mathcal{B}(\mathbf{w}, \mathbf{w}),~\forall~\mathbf{w} \in \mathbb{W}.
\end{align*}
And suppose that the norm $\|\cdot\|_{\mathcal{A}}$ is relatively compact with respect to the norm $\|\cdot\|_{\mathcal{B}}$ \cite{Conway}. Then we introduce the following variational eigenvalue problem: 
Find $(\lambda, \mathbf{u}) \in \mathbb{R} \times \mathbb{V}$ such that $\mathcal{B}(\mathbf{u}, \mathbf{u})=1$ and
\begin{align}
\mathcal{A}(\mathbf{u}, \mathbf{v})=\lambda \mathcal{B}(\mathbf{u}, \mathbf{v}), \ \forall~\mathbf{v} \in \mathbb{V}. \label{eq:abstract}
\end{align}
In order to solve the above variational eigenvalue problem, taking advantage of the characteristics of  FieldTNN structure \eqref{eq:FieldTNN}, we can construct a subspace $\mathbb{V}_{p}$ of the Hilbert space $\mathbb{V}$ by FieldTNN basis functions as follows:
\begin{equation}\label{eq:FTNN_space}
\mathbb{V}_{p}:=\operatorname{span}\left\{\mathbf{\Psi}_{k}(\mathbf{x};\mathbf{\Theta}) = \left(\begin{array}{l} 
\prod\limits_{j=1}^{3}\psi_{1,j,k}\left(x_{j};\theta_{1,j}\right) \\
\prod\limits_{j=1}^{3}\psi_{2,j,k}\left(x_{j};\theta_{2,j}\right) \\
\prod\limits_{j=1}^{3}\psi_{3,j,k}\left(x_{j};\theta_{3,j}\right) 
\end{array}\right),~k = 1, \cdots, p\right\},
\end{equation}
where the FieldTNN basis functions $\mathbf{\Psi}_{k}(\mathbf{x},\mathbf{\Theta})$ are defined in \eqref{eq:FieldTNNk}.
Then, FieldTNN-Galerkin approximation scheme is: Find $(\lambda, \mathbf{u}_{\scriptscriptstyle \text{NN}}) \in \mathbb{R} \times \mathbb{V}_{p}$ such that $\mathcal{B}(\mathbf{u}_{\scriptscriptstyle \text{NN}}, \mathbf{u}_{\scriptscriptstyle \text{NN}})=1$ and
\begin{align}
\mathcal{A}(\mathbf{u}_{\scriptscriptstyle \text{NN}}, \mathbf{v}_{\scriptscriptstyle \text{NN}}) = \lambda_{\scriptscriptstyle \text{NN}} \mathcal{B}(\mathbf{u}_{\scriptscriptstyle \text{NN}}, \mathbf{v}_{\scriptscriptstyle \text{NN}}), \quad \forall~\mathbf{v}_{\scriptscriptstyle \text{NN}} \in \mathbb{V}_{p}. \label{eq:Field_Galerkin}
\end{align}
By some simple calculations, we can derive the following generalized eigenvalue problems
\begin{align}
&\mathbf{S}\mathrm{U} = \lambda_{\scriptscriptstyle \text{NN}} \mathbf{M}\mathrm{U}, \label{eq:GEQ}
\end{align}
where $\mathrm{U} \in \mathbb{R}^{p \times 1}$ is column vector, and $\mathbf{S}, \mathbf{M} \in \mathbb{R}^{p \times p}$ denote ``stiffness'' and mass matrix with entries
\begin{align}
\mathrm{s}_{\scriptscriptstyle jk}=\mathcal{A}(\mathbf{\Psi}_{k}(\mathbf{x} ; \mathbf{\Theta}),\mathbf{\Psi}_{j}(\mathbf{x} ; \mathbf{\Theta})),\ \ \ \mathrm{m}_{\scriptscriptstyle jk}=\mathcal{B}(\mathbf{\Psi}_{k}(\mathbf{x} ; \mathbf{\Theta}),\mathbf{\Psi}_{j}(\mathbf{x} ; \mathbf{\Theta})),\ \ \ j,k=1,\cdots,p,
\label{eq:entriesAB}
\end{align}
respectively.
By solving the generalized eigenvalue problem \eqref{eq:GEQ}, we can obtain the leading $p$ eigenvalues:
\begin{align}
&\lambda_{1,\scriptscriptstyle \text{NN}} \leq \cdots \leq \lambda_{p,\scriptscriptstyle \text{NN}}
\label{eq:FTNNeigva}
\end{align}
and the corresponding eigenvectors
\begin{align}
&\mathrm{U}_{1}, \cdots, \mathrm{U}_{p}, \label{eq:FTNNeigve}
\end{align}
where $\mathrm{U}_{\ell} = \left(u_{1\ell},u_{2\ell},\cdots,u_{p\ell}\right)^{\top}$, $\ell=1,\cdots,p$.
Then $\lambda_{1,\scriptscriptstyle \text{NN}}, \cdots, \lambda_{p,\scriptscriptstyle \text{NN}}$ can be chosen as the approximations to the leading $p$ eigenvalues $\lambda_{1}, \cdots$, $\lambda_{p}$.
And the corresponding approximation eigenfunctions $\mathbf{u}_{1,\scriptscriptstyle \text{NN}}, \cdots, \mathbf{u}_{p,\scriptscriptstyle \text{NN}}$ to the leading $p$ eigenfunctions of problem \eqref{eq:abstract} can be obtained by the linear combination the FieldTNN basis functions as follows
\begin{align*}
&\mathbf{u}_{\ell,\scriptscriptstyle \text{NN}}(\mathbf{x},\mathbf{\Theta}) = \sum_{k=1}^{p} u_{k\ell}\mathbf{\Psi}_{k}(\mathbf{x} ; \mathbf{\Theta}),\ \  \ell=1,\cdots,p.
\end{align*}
During the training of  FieldTNN, the single loss function $\mathcal{L}\left(\lambda_{\ell,\scriptscriptstyle \text{NN}}, \mathbf{u}_{\ell,\scriptscriptstyle \text{NN}}\right)$ for the $\ell$-th eigenpairs is defined by
\begin{align}
&\mathcal{L}\left(\lambda_{\ell,\scriptscriptstyle \text{NN}}, \mathbf{u}_{\ell,\scriptscriptstyle \text{NN}}\left(\mathbf{x} ; \mathbf{\Theta}\right)\right) = \lambda_{\ell,\scriptscriptstyle \text{NN}} + \beta \cdot {\rm constraint},\ \ \ell = 1,\cdots,p,\label{eq:singleLoss}
\end{align}
where the constrained term on the right side
will be specified for different problem, i.e., the penalty for the divergence-free condition in the Maxwell equations, and $\beta > 0$ is a penalty parameter.
For solving leading eigenpairs, we sort the single loss functions \eqref{eq:singleLoss} and obtain
\begin{align*}
&\mathcal{L}\left(\lambda_{1,\scriptscriptstyle \text{NN}}, \mathbf{u}_{1,\scriptscriptstyle \text{NN}}(\mathbf{x} ; \mathbf{\Theta})\right) \leq \mathcal{L}\left(\lambda_{2,\scriptscriptstyle \text{NN}}, \mathbf{u}_{2,\scriptscriptstyle \text{NN}}(\mathbf{x} ; \mathbf{\Theta})\right) \leq \cdots \leq \mathcal{L}\left(\lambda_{p,\scriptscriptstyle \text{NN}}, \mathbf{u}_{p,\scriptscriptstyle \text{NN}}(\mathbf{x} ; \mathbf{\Theta})\right).
\end{align*}
Then the total loss function $\widehat{\mathcal{L}}\left(\mathbf{\Theta}\right)$ of the leading $\mathcal{M}$ eigenpairs is defined as 
\begin{align}
&\widehat{\mathcal{L}}\left(\mathbf{\Theta}\right) = \sum_{m = 1}^{\mathcal{M}}\mathcal{L}\left(\lambda_{m,\scriptscriptstyle \text{NN}}, \mathbf{u}_{m,\scriptscriptstyle \text{NN}}(\mathbf{x} ; \mathbf{\Theta})\right),~\mathcal{M} \leq p,\label{eq:totalLoss}
\end{align}
which is a direct generalization of minimum-maximum principle \cite{Book_Babuska} and Rayleigh-Ritz method \cite{Book_Saad}.
The loss function \eqref{eq:totalLoss} can be automatically differentiated by backward propagation
which is supported by open-source deep learning frameworks such as TensorFlow \cite{TensorFlow} and PyTorch \cite{PyTorch}. In this paper, the gradient descent (GD) method is adopted to update all trainable parameters
\begin{align}
&\mathbf{\Theta} - \eta \nabla_{\mathbf{\Theta}}\widehat{\mathcal{L}}\left(\mathbf{\Theta}\right) \rightarrow \mathbf{\Theta},\label{eq:update}
\end{align}
where $\eta$ is the learning rate and adjusted by the ADAM optimizer. After sufficient training steps, we obtain parameters $\mathbf{\Theta}^{*}$ such that the loss function \eqref{eq:totalLoss} arrives its minimum value under the required tolerance.
\section{TNN-based machine learning method for problems in non-tensor domain}\label{sec:nontensor}
In previous work on TNN \cite{JCP_TNN, JCM_TNN}, the discussion of computational domains was limited to tensor-type domains.
In this section, we define TNN on the non-tensor domains and propose a corresponding TNN-based machine learning method, which is one of the innovation points of this paper.
Our approach is inspired by the work of Kwan and Shen \cite{SEM_JCP}, which introduced a spectral element discretization for elliptic problems using modal basis functions.
In their work, the computational domain is decomposed into some subdomains—whether they are line segments, squares, or cubes—and then a combination of compactly supported basis functions and hat basis functions are employed for implementation.
This inspires us to construct locally compactly supported TNNs and use them as basis functions to discrete partial differential equations.
Thanks to the variable separation nature of TNN, the key idea is to construct one-dimensional, locally compactly supported subnetworks.
For readers to understand, we present a detailed construction here.

Denote $\sigma(x) = \max\{0,x\}$ as the Rectified Linear Unit (ReLU) function.
Then, we define an auxiliary function by a combination of ReLUs as follows,
\begin{align}\label{eq:def_g}
&g_{[a,b]}(x) = \sigma(x - a) - \sigma(x - b) + a,
\end{align}
one can verify that 
\begin{align}\label{eq:gxab}
& g_{[a,b]}(x) = 
\begin{cases} 
a, & x < a, \\
x, & a \leq x \leq b,\\
b, & x > b.
\end{cases} 
\end{align}
By functions composition, any function $f(\cdot)$ can be focused on the interval $[a,b]$, i.e. the following property holds
\begin{align}\label{eq:fx}
& f\left(g_{[a,b]}(x)\right) = 
\begin{cases} 
f(a), & x < a, \\
f(x), & a \leq x \leq b,\\
f(b), & x > b.
\end{cases} 
\end{align}
Specially, if $f$ satisfies $f(a)=f(b)=0$, then support of the composite function $f\left(g_{[a,b]}(x)\right)$ is $[a,b]$.
According to this property, constructing one-dimensional locally compactly supported subnetworks can be implemented by constructing subnetworks with a zero value at the end of the interval.
More specific, for $j = 1,~2,~\cdots,~d$, we define the $j$-th subnetwork $\tilde{\psi}_{j}(x_{j}; \theta_{j})$ as follows
\begin{align*}
\tilde{\psi}_{j}(x_{j}; \theta_{j}) :& = \hat{\psi}_{j}(x_{j}; \theta_{j})(x_{j}-a_{j})(b_{j}-x_{j})\\
& = \left(\hat{\psi}_{j,1}(x_{j}; \theta_{j})(x_{j}-a_{j})(b_{j}-x_{j}),\cdots,\hat{\psi}_{j,p}(x_{j}; \theta_{j})(x_{j}-a_{j})(b_{j}-x_{j})\right)^{\top},
\end{align*}
where $\hat{\psi}_{j}(x_{j}; \theta_{j})$  is an FNN from $\mathbb{R}$ to $\mathbb{R}^{p}$ with appropriate activation functions and is defined in \eqref{eq:TNNjth}.
Then we have
\begin{align*}
\tilde{\psi}_{j,k}\big(g_{[a_{j},b_{j}]}\left(x_{j}\right); \theta_{j}\big) =
\begin{cases} 
\tilde{\psi}_{j,k}(a_{j}; \theta_{j}) = \hat{\psi}_{j,k}(a_{j};\theta_{j})(a_{j}-a_{j})(b_{j}-a_{j}) = 0, & x_{j} \leq a_{j}, \\
\tilde{\psi}_{j,k}\big(x_{j};\theta_{j}\big), & a_{j} < x_{j} < b_{j},\\
\tilde{\psi}_{j,k}(b_{j}; \theta_{j}) = \hat{\psi}_{j,k}(b_{j};\theta_{j})(b_{j}-a_{j})(b_{j}-b_{j}) = 0, & x_{j} \geq b_{j},
\end{cases}
\end{align*}
where $j = 1,~\cdots,~d$, $k=1,\cdots,p$ and $g_{[a_{j},b_{j}]}(x)$ is defined in \eqref{eq:gxab}.
Finally, we obtain the following TNN $\Phi\left(\mathbf{x}; \Theta\right)$ which is compactly supported on $\prod\limits_{j=1}^d[a_j,b_j]$:
\begin{align}
\Phi\left(\mathbf{x}; \Theta\right) & = \sum_{k = 1}^{p}u_{k}\prod_{j = 1}^{d}\tilde{\psi}_{j,k}(g_{[a_{j},b_{j}]}(x_{j}); \theta_{j}).\label{eq:suppTNN}
\end{align}
\begin{remark}
Noting that the definition \eqref{eq:def_g} can be written in the following neural network format
\begin{eqnarray}
g_{[a,b]}(x)=
\begin{pmatrix}
1 & -1
\end{pmatrix}
\sigma\left(
\begin{pmatrix}
1\\ 1
\end{pmatrix}
x+
\begin{pmatrix}
-a\\ -b
\end{pmatrix}
\right)+a,
\end{eqnarray}
i.e. $g_{[a,b]}(\cdot)$ is an FNN mapping form $\mathbb R$ to $\mathbb R$ with one hidden layer and two hidden neurons.
Then the composite function $\hat{\psi}_{j}\big(g_{[a_{j},b_{j}]}\left(\cdot\right); \theta_{j}\big)$ can be obtained by simply concatenating two FNNs $g_{[a_{j},b_{j}]}(\cdot)$ and $\hat{\psi}_{j}(\cdot)$.

It is worth mentioning that in the implementation, we do not need actually to concatenate the two networks.
When calculating entries in the ``stiffness'' matrix and the mass matrix
\eqref{eq:entriesAB}, we use the assembling technique of ``stiffness'' matrix in finite element method \cite{NMPDES_2017}.
The ``stiffness'' matrix of the whole domain is obtained by summing the corresponding ``stiffness`` matrix on each subdomain according to the local compact support property of the TNN basis function.
\end{remark}

\begin{figure}[H]
\centering
\subfloat[{L-shaped domain $\Omega = [-1,1]^{2}/\left\{[0,1]\times[-1,0]\right\}$.}]{\includegraphics[width=0.23\linewidth]{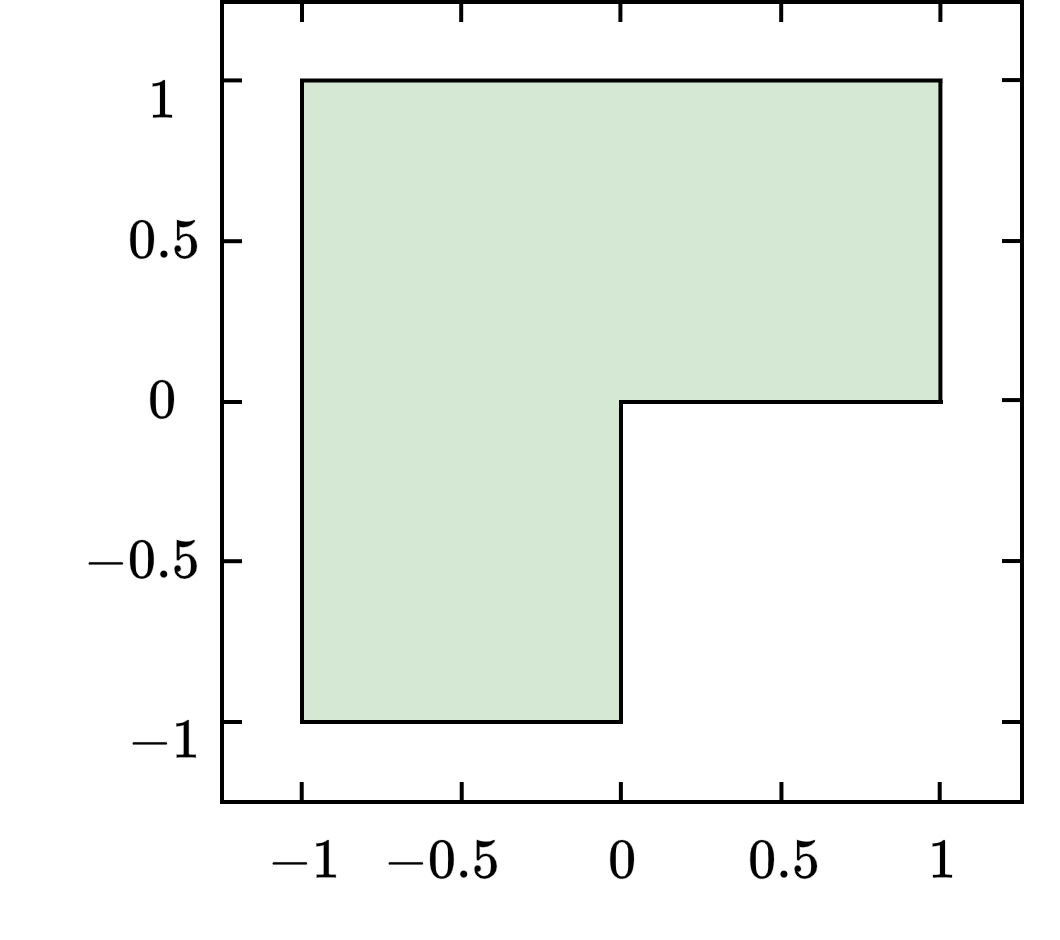}}\hfill
\subfloat[{Decomposition of the L-shaped domain.}]{\includegraphics[width=0.23\linewidth]{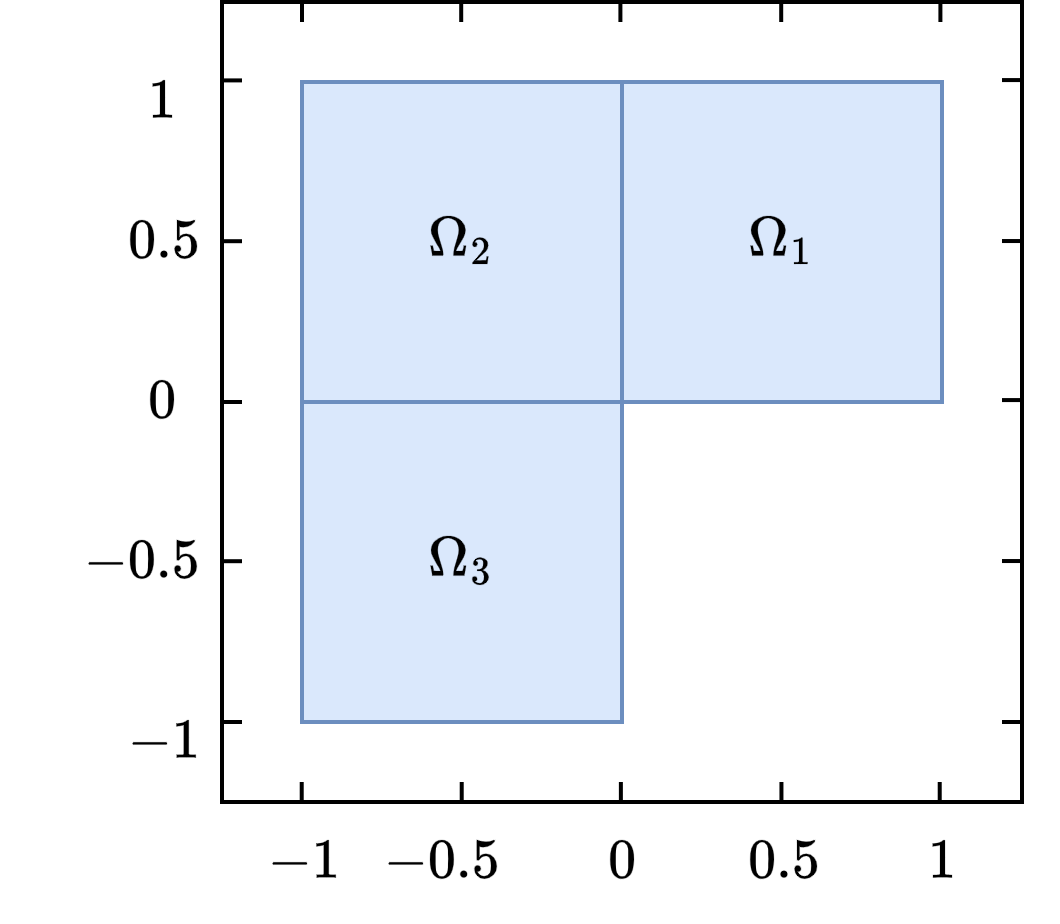}}\hfill
\subfloat[{Local compact supported TNNs in the L-shaped domain.}]{\includegraphics[width=0.3\linewidth]{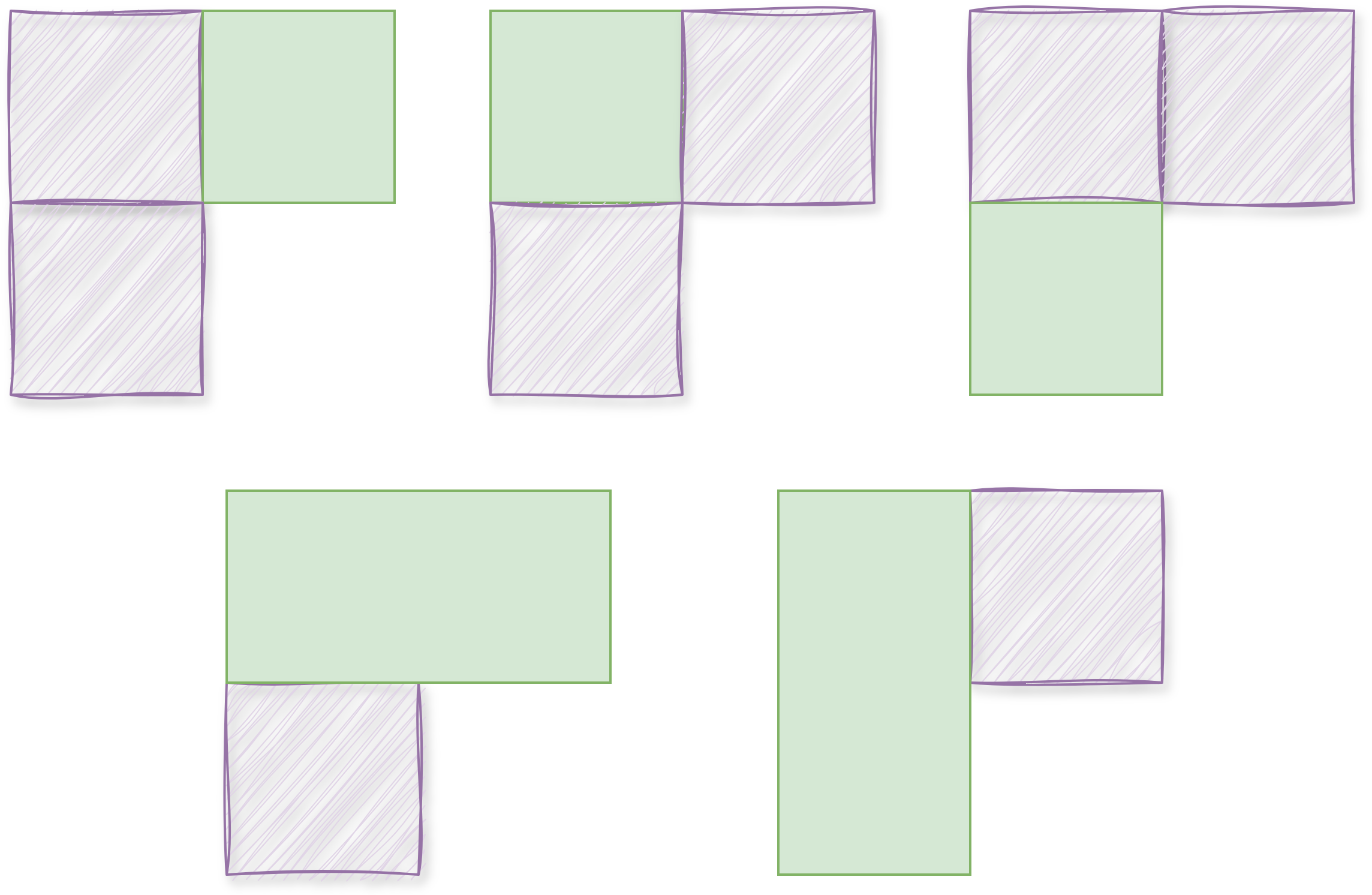}}\\
\caption{Illustration for the L-shaped domain.}
\label{fig:sec4}
\end{figure}

Here, taking the L-shaped domain as an example, we propose a novel TNN-based machine learning method for solving problems in non-tensor domains.
For L-shaped domain shown in Figure \ref{fig:sec4} $(1)$,
we first decompose the computational domain $\Omega$ into $3$ subdomains, as illustrated in Figure \ref{fig:sec4} $(2)$, i.e., $\Omega_{1} = [0,1]\times[0,1]$, $\Omega_{2} = [-1,0]\times[0,1]$, and $\Omega_{3} = [-1,0]\times[-1,0]$, and further define two union domain: $\Omega_{4} = \Omega_{1} \bigcup \Omega_{2}$, $\Omega_{5} = \Omega_{2} \bigcup \Omega_{3}$. 
In each domain $\Omega_{\ell}~(\ell=1,\cdots,5)$, we will train local compact supported TNNs $\Phi(\mathbf{x}; \Theta)$ defined in \eqref{eq:suppTNN}, where $\mathbf{x} \in \Omega_{\ell}$. 
An illustration of the local compact supported TNNs is provided in Figure \ref{fig:sec4} $(3)$, where functions are local compact supported in green subdomain and zero in purple ones.
Simple calculation yields the following block matrix form:
\begin{align}
\mathbf{A} = \begin{pmatrix}
\star        & \mathbf{0} & \mathbf{0} & \star      & \mathbf{0}  \\
\mathbf{0}   & \star      & \mathbf{0} & \star      & \star \\
\mathbf{0}   & \mathbf{0} & \star      & \mathbf{0} & \star \\
\star        & \star      & \mathbf{0} & \star      & \star \\
\mathbf{0}   & \star      & \star      & \star      & \star
\end{pmatrix},\label{eq:SPMatrix}
\end{align}
where $\star$ denotes non-zero matrix, and $\mathbf{0}$ is zero matrix.

\begin{remark}[Remark for FieldTNN]
In the above discussion, we only consider the TNN for the problem in non-tensor domains. 
However, similar technique can be employed to FieldTNN. 
Taking $d = 3$ as an example, one can easily construct 
\begin{align*}
\tilde{\psi}_{i,j,k}\big(g_{[a_{j},b_{j}]}\left(x_{j}\right); \theta_{i,j}\big) =
\begin{cases} 
\tilde{\psi}_{i,j,k}(a_{j}; \theta_{i,j}) = \hat{\psi}_{i,j}(a_{j};\theta_{i,j})(a_{j}-a_{j})(b_{j}-a_{j}) = 0, & x_{j} \leq a_{j}, \\
\tilde{\psi}_{i,j,k}\big(x_{j};\theta_{i,j}\big), & a_{j} < x_{j} < b_{j},\\
\tilde{\psi}_{i,j,k}(b_{j}; \theta_{i,j}) = \hat{\psi}_{i,j}(b_{j};\theta_{i,j})(b_{j}-a_{j})(b_{j}-b_{j}) = 0, & x_{j} \geq b_{j},
\end{cases}
\end{align*}
and
\begin{align*}
\Phi_{i,k}(\mathbf{x}; \Theta_{i}\big) & = \prod_{j = 1}^{3}\tilde{\psi}_{i,j,k}\left(g_{[a_{j},b_{j}]}\left(x_{j}\right); \theta_{i,j}\right),
\end{align*}
where $i,j = 1,2,3$ and $k = 1,\cdots,p$. 
Then, we have
\begin{equation*}
\begin{aligned}
\mathbf{\Phi}_{k}\left(\mathbf{x};\mathbf{\Theta} \right) =  
\left(
\begin{array}{l}
\Phi_{1,k}(\mathbf{x}; \Theta_{1})\\
\Phi_{2,k}(\mathbf{x}; \Theta_{2})\\
\Phi_{3,k}(\mathbf{x}; \Theta_{3})
\end{array}
\right).
\end{aligned}
\end{equation*}
Finally, we obtain the following FieldTNN $\mathbf{\Phi}\left(\mathbf{x}; \mathbf{\Theta}\right)$ which is compactly supported on $\prod\limits_{i=1}^{3}[a_i,b_i]$:
\begin{align}
& \mathbf{\Phi}(\mathbf{x} ; \mathbf{\Theta}) = \sum_{k = 1}^{p} u_{k} \mathbf{\Phi}_{k}(\mathbf{x} ; \mathbf{\Theta}).\label{eq:nonTensorFieldTNN} 
\end{align}
\end{remark}
\section{Solving Maxwell eigenvalue problems by FieldTNN}\label{sec:Maxwell}
In this section, we present FieldTNN-based machine learning approach to solve the Maxwell eigenvalue problems.
The focus is on applying these neural networks to compute the eigenvalues and eigenfunctions associated with Maxwell equations efficiently.
The methodology we propose is particularly effective in addressing the challenges of electromagnetic wave propagation in different domains, making it a powerful tool for practical computations.
For simplicity, we only discuss the case where $d=3$ in this section, the results for $d=2$ can be derived similarly.
\subsection{Quadrature scheme for variational form}
Recalling the variational form presented in \eqref{eq:vform2}:
Find $\mathbf{0} \neq \mathbf{E} \in \mathbf{H}_{0}(\operatorname{curl};\Omega)$
and $\lambda \in \mathbb{R}$, such that
\begin{align}\label{eq:VEP}
&\mathcal{A}(\mathbf{E}, \mathbf{V}) = \lambda \mathcal{B}(\mathbf{E}, \mathbf{V}),~\forall~\mathbf{V} \in \mathbf{H}_{0}(\operatorname{curl}; \Omega),
\end{align}
where
\begin{align*}
\mathcal{A}(\mathbf{E}, \mathbf{V})  = \left(\nabla \times \mathbf{E},\nabla \times \mathbf{V}\right), \ \mathcal{B}(\mathbf{E}, \mathbf{V}) = (\mathbf{E}, \mathbf{V}),
\end{align*}
Next, denote by $\mathbb{U}_{p}$ a subspace of $\mathbf{H}_{0}(\operatorname{curl};\Omega)$ such that
\begin{align*}
&\mathbb{U}_{p}	:= \mathbf{H}_{0}(\operatorname{curl};\Omega) \cap \mathbb{V}_{p},
\end{align*}
where $\mathbb{V}_{p}$ is a subspace spanned by FieldTNN \eqref{eq:FTNN_space} and $\mathbb{U}_{p}$ satisfies the boundary conditions or constraints required by the problem.
In this paper, in order to solve the Maxwell eigenvalue problems, it is necessary to construct FiledTNN satisfying the boundary conditions and divergence-free conditions to obtain $\mathbb{U}_p$, we will explain this in detail in Section \ref{sec:details}.

Using the subspace  $\mathbb{U}_p$ to discrete variational eigenvalue problem \eqref{eq:VEP}, 
one can derive the FieldTNN-Galerkin approximation: Find $0 \neq \mathbf{E}_{\scriptscriptstyle \text{NN}} \in \mathbb{U}_{p}$ and $ \lambda \in \mathbb{R}$, such that
\begin{align}
&\mathcal{A}(\mathbf{E}_{\scriptscriptstyle \text{NN}}, \mathbf{V}_{\scriptscriptstyle \text{NN}}) = \lambda \mathcal{B}(\mathbf{E}_{\scriptscriptstyle \text{NN}}, \mathbf{V}_{\scriptscriptstyle \text{NN}}),~\forall~\mathbf{V}_{\scriptscriptstyle \text{NN}} \in \mathbb{U}_{p}.\label{eq:FTNNG}
\end{align}
We can write the neural network approximation eigenfunction $\mathbf{E}_{\scriptscriptstyle \text{NN}}$ as the combination of FieldTNN basis functions of $\mathbb{U}_{p}$,
\begin{align}
&\mathbf{E}_{\scriptscriptstyle \text{NN}} = \sum_{k=1}^{p} u_{k} \mathbf{\widehat{\Phi}}_{k}\left(\mathbf{x} ; \mathbf{\Theta}\right),\label{eq:FTNNtrial}
\end{align}
where $\left\{\mathbf{\widehat{\Phi}}_{k}\left(\mathbf{x} ; \mathbf{\Theta}\right)\right\}_{k=1}^{p}$ is a set of basis of $\mathbb{U}_{p}$.
Inserting Eq. \eqref{eq:FTNNtrial} into problem \eqref{eq:FTNNG} and taking the test function $\mathbf{V}_{\scriptscriptstyle \text{NN}}$ through a set of basis in the approximation space $\mathbb{U}_{p}$, we can derive the a generalized matrix eigenvalue problem as follows,
\begin{align}
&\mathbf{S}\mathrm{U} = \lambda_{\scriptscriptstyle \text{NN}} \mathbf{M} \mathrm{U}, \label{eq:GMEP}
\end{align}
where $\mathbf{S}, \mathbf{M} \in \mathbb{R}^{p \times p}$ and $\mathrm{U} \in \mathbb{R}^{p \times 1}$. Here, the entries of matrices $\mathbf{S}$ and $\mathbf{M}$ and column vector $\mathrm{U}$ have the following form,
\begin{equation}\label{eq:Stiff}
\begin{aligned}
s_{\scriptscriptstyle ji}= & \mathcal{A}\left(\mathbf{\widehat{\Phi}}_i\left(\mathbf{x}; \mathbf{\Theta}\right), \mathbf{\widehat{\Phi}}_j\left(\mathbf{x} ; \mathbf{\Theta}\right)\right) \\
= & \sum_{n = 1}^{3} \sum_{s \neq n }^{3} \prod_{k \neq s}^{3} \int_{\Omega_{k}} \hat{\phi}_{n, k, i}\left(x_{\scriptscriptstyle k}; \theta_{n, k}\right) \hat{\phi}_{n, k, j}\left(x_{\scriptscriptstyle k}; \theta_{n, k}\right) \mathrm{d}x_{\scriptscriptstyle k} 
\int_{\Omega_{s}} \frac{\partial \hat{\phi}_{n, s, i}}{\partial x_{s}}\left(x_{s}; \theta_{n, s}\right) \frac{\partial \hat{\phi}_{n, s, j}}{\partial x_{s}}\left(x_{s}; \theta_{n, s}\right)  \mathrm{d}x_{s} \\
&-\sum_{n = 1}^{3} \sum_{s \neq n }^{3} \prod_{k \neq n,s}^{3} \int_{\Omega_k} \hat{\phi}_{n, k, i}\left(x_{\scriptscriptstyle k}; \theta_{n, k}\right) \hat{\phi}_{s, k, j}\left(x_{\scriptscriptstyle k} ; \theta_{n, k}\right) \mathrm{d} x_{\scriptscriptstyle k} \\
& \cdot \int_{\Omega_s}  \frac{\partial \hat{\phi}_{n, s, i}}{\partial x_{s}}\left(x_s ; \theta_{n, s}\right) \hat{\phi}_{s, s, j} \left(x_s ; \theta_{s, s}\right)  \mathrm{d}x_s \int_{\Omega_n} \hat{\phi}_{n, n, i} \left(x_n ; \theta_{n, n}\right) \frac{\partial \hat{\phi}_{s, n, j}}{\partial x_{n}}\left(x_n ; \theta_{s, n}\right)  \mathrm{d}x_{n},
\end{aligned}
\end{equation}
and 
\begin{align}
m_{\scriptscriptstyle ji} & = \mathcal{B}\Big(\mathbf{\widehat{\Phi}}_i\left(\mathbf{x} ; \mathbf{\Theta}\right), \mathbf{\widehat{\Phi}}_j\left(\mathbf{x} ;  \mathbf{\Theta}\right)\Big) \nonumber\\
& =\sum_{n=1}^{3} \prod_{k=1}^{3} \int_{\Omega_k} \phi_{n, k, i}\left(x_k ; \theta_{n, k}\right) \phi_{n, k, j}\left(x_k ; \theta_{n, k}\right) \mathrm{d} x_{k},\label{eq:Mass}
\end{align}
and
\begin{align}
\mathrm{U} &= \left(u_{1},u_{2},\cdots,u_{p}\right)^{\top}.
\end{align}
Given that only one-dimensional integrations are required in \eqref{eq:Stiff} and \eqref{eq:Mass},
to ensure the integrations within the loss functions are highly accurate, instead of Monte-Carlo method,
it is recommended to use a high-order one-dimensional quadrature scheme, such as Gauss-type rules, i.e., Gauss, Gauss-Radau and Gauss-Lobatto quadrature \cite{SpectralBook2011}. These can be expressed in a general form as follows:
\begin{align}
&\int_{a}^{b} f(x)dx \approx \sum_{\ell=1}^{N} w_{\ell} f(x_{\ell}),\label{eq:quadrature}
\end{align}
where $\{x_{\ell}\}_{\ell = 1}^{N}$ are the quadrature points and $\{w_{\ell}\}_{\ell = 1}^{N}$ are the corresponding weights.
For generality, we decompose each $\Omega_i$ into $M_i$ equal subintervals, where the length of each subinterval is $h_{i} = \frac{\left|\Omega_{i}\right|}{M_{i}}$, and select $N_{i}$ Legendre-Gauss points within each subinterval. 
The total quadrature points and corresponding weights are defined as follows,
\begin{align}
&\left\{x_i^{\left(\ell_{i}\right)}\right\}_{\ell_{i}=1}^{M_i N_i}, \quad\left\{w_i^{\left(\ell_{i}\right)}\right\}_{\ell_{i}=1}^{M_i N_i},\ i=1,2,3.\label{eq:LGPW}
\end{align}
Then using quadrature scheme \eqref{eq:quadrature} equipped with the Legendre-Gauss points \eqref{eq:LGPW},
entries of ``stiffness'' matrix $\mathrm{S}$ \eqref{eq:Stiff} and mass matrix $\mathrm{M}$ \eqref{eq:Mass} can be written in the following form:
\begin{align}
s_{\scriptscriptstyle ji} = &  \mathcal{A}\Big(\mathbf{\widehat{\Phi}}_i\left(\mathbf{x}; \mathbf{\Theta}\right), \mathbf{\widehat{\Phi}}_j\left(\mathbf{x} ; \mathbf{\Theta}\right)\Big) \nonumber\\
\approx &  \sum_{n = 1}^{3} \sum_{s \neq n }^{3} \prod_{k \neq s}^{3} \sum_{\ell_{k} = 1}^{M_{k} N_{k}} w_{\scriptscriptstyle k}^{(\ell_{k})}\hat{\phi}_{n, k, i}\left(x_{\scriptscriptstyle k}^{(\ell_{k})}; \theta_{n, k}\right) \hat{\phi}_{n, k, j}\left(x_{\scriptscriptstyle k}^{(\ell_{k})}; \theta_{n, k}\right)\nonumber\\
&\cdot \sum_{\ell_{s} = 1}^{M_{s} N_{s}} w_{s}^{(\ell_{s})}\frac{\partial \hat{\phi}_{n, s, i}}{\partial x_{s}}\left(x_{s}^{(\ell_{s})}; \theta_{n, s}\right) \frac{\partial \hat{\phi}_{n, s, j}}{\partial x_{s}}\left(x_{s}^{(\ell_{s})}; \theta_{n, s}\right) \nonumber\\
&-\sum_{n = 1}^{3} \sum_{s \neq n }^{3} \prod_{k \neq n,s}^{3} \sum_{\ell_{k} = 1}^{M_{k} N_{k}} w_{\scriptscriptstyle k}^{(\ell_{k})} \hat{\phi}_{n, k, i}\left(x_{\scriptscriptstyle k}^{(\ell_{k})} ; \theta_{n, k}\right) \hat{\phi}_{s, k, j}\left(x_{\scriptscriptstyle k}^{(\ell_{k})} ; \theta_{n, k}\right) \nonumber\\
& \cdot \sum_{\ell_{s} = 1}^{M_{s} N_{s}} w_{s}^{(\ell_{s})} \frac{\partial \hat{\phi}_{n, s, i}}{\partial x_{s}}\left(x_s^{(\ell_{s})} ; \theta_{n, s}\right) \hat{\phi}_{s, s, j} \left(x_s^{(\ell_{s})} ; \theta_{s, s}\right)\nonumber\\
&\cdot \sum_{\ell_{n} = 1}^{M_{n} N_{n}} w_{n}^{(\ell_{n})} \hat{\phi}_{n, n, i} \left(x_{n}^{(\ell_{n})} ; \theta_{n, n}\right) \frac{\partial \hat{\phi}_{s, n, j}}{\partial x_{n}}\left(x_{n}^{(\ell_{n})} ; \theta_{s, n}\right),\label{eq:StiffnessGQ}
\end{align}
and 
\begin{align}
m_{\scriptscriptstyle ji} & = \mathcal{B}\Big(\mathbf{\widehat{\Phi}}_i\left(\mathbf{x} ; \mathbf{\Theta}\right), \mathbf{\widehat{\Phi}}_j\left(\mathbf{x} ;  \mathbf{\Theta}\right)\Big) \nonumber\\
& = \sum_{n=1}^{3}\left(\widehat{\Phi}_{n,i}\left(\mathbf{x};\Theta_{n}\right), \widehat{\Phi}_{n,j}\left(\mathbf{x};\Theta_{n}\right)\right) \nonumber\\
& =\sum_{n=1}^{3} \prod_{k=1}^{3} \sum_{\ell_{k} = 1}^{M_{k} N_{k}} w_{\scriptscriptstyle k}^{(\ell_{k})} \phi_{n, k, i}\left(x_{\scriptscriptstyle k}^{(\ell_{k})} ; \theta_{n, k}\right) \phi_{n, k, j}\left(x_{\scriptscriptstyle k}^{(\ell_{k})} ; \theta_{n, k}\right),\label{eq:MassGQ}
\end{align}
respectively. 
\begin{remark}{\bf(Complexity analysis)}
It is noteworthy that the quadrature schemes \eqref{eq:StiffnessGQ} and \eqref{eq:MassGQ}, which utilize $M_i N_i$ one-dimensional quadrature points per dimension, are in fact equivalent to employing a full three-dimensional tensor product quadrature scheme, i.e.,
\begin{align*}
\mathrm{s}_{\scriptscriptstyle ji} & \approx  \sum_{\ell \in \mathcal{N}} w^{(\ell)} \nabla \times \mathbf{\widehat{\Psi}}_{i}(\mathbf{x}^{(\ell)}; \mathbf{\Theta}) \cdot \nabla \times \mathbf{\widehat{\Psi}}_{j}(\mathbf{x}^{(\ell)}; \mathbf{\Theta}), \\
\mathrm{m}_{\scriptscriptstyle ji} & \approx  \sum_{\ell \in \mathcal{N}} w^{(\ell)} \mathbf{\widehat{\Psi}}_{j}(\mathbf{x}^{(\ell)}; \mathbf{\Theta}) \mathbf{\widehat{\Psi}}_{i}(\mathbf{x}^{(\ell)}; \mathbf{\Theta}),
\end{align*}
where $3$-dimensional quadrature points and weights are defined as follows
\begin{equation}\label{eq:d_dimPW}
\begin{aligned}
&\left\{x^{(\ell)}\right\}_{\ell \in \mathcal{N}} = \left\{\left\{x_1^{\left(\ell_1\right)}\right\}_{\ell_1=1}^{M_1 N_1} \times \left\{x_2^{\left(\ell_2\right)}\right\}_{\ell_2=1}^{M_2 N_2} \times \left\{x_3^{\left(\ell_3\right)}\right\}_{\ell_3=1}^{M_3 N_3}\right\}, \\
&\left\{w^{(\ell)}\right\}_{\ell \in \mathcal{N}} = \left\{\left\{w_1^{\left(\ell_1\right)}\right\}_{\ell_1=1}^{M_1 N_1} \times\left\{w_2^{\left(\ell_2\right)}\right\}_{\ell_2=1}^{M_2 N_2} \times \left\{w_3^{\left(\ell_3\right)}\right\}_{\ell_3=1}^{M_3 N_3}\right\} .
\end{aligned}
\end{equation}
The quadrature points \eqref{eq:d_dimPW} used for $3$-dimensional integrations achieve an accuracy of 
$\mathcal{O}\left(\frac{h^2 \underline{N}}{(2 \underline{N})!}\right)$,
while the number of elements in $\mathcal{N}$ scales as $\mathcal{O}\left(M^3 N^3\right)$.
Here, $N = \max\{N_{1},\cdots,N_{d}\}$, $\underline{N} = \min\{N_{1},\cdots,N_{d}\}$, $M = \max\{M_{1},\cdots,M_{d}\}$, $\underline{M} = \min\{M_{1},\cdots,M_{d}\}$.
However, for large dimension $d$, this becomes prohibitively large, rendering the integration approximation unmanageable.
Due to the low-rank tensor structure of FieldTNNs, applying the full tensor quadrature scheme to compute the inner product of two FieldTNNs leads to the splitting schemes \eqref{eq:StiffnessGQ} and \eqref{eq:MassGQ}.
In these schemes, assembling the matrices $\mathrm{S}$ and $\mathrm{M}$ requires $\mathcal{O}\left(d^{3}MN\right)$ and $\mathcal{O}\left(d^{2}M N\right)$ operations, respectively.
Furthermore, as the full tensor quadrature points \eqref{eq:d_dimPW} is equivalent to the quadrature schemes \eqref{eq:StiffnessGQ} and \eqref{eq:MassGQ}, the quadrature schemes \eqref{eq:StiffnessGQ} and \eqref{eq:MassGQ} can also achieve an accuracy of $\mathcal{O}\left(\frac{h^2 \underline{N}}{(2 \underline{N})!}\right)$.
This demonstrates that the FieldTNN-based method proposed in this paper successfully overcomes the ``curse of dimensionality`` in the context of numerical integration.
\end{remark}

\subsection{Technical details}\label{sec:details}
In this subsection, we delve into the specific technical aspects of our approach. 
We begin by addressing the divergence-free condition, which is crucial for ensuring the physical accuracy of the solutions in electromagnetic wave simulations. 
Subsequently, we discuss the implementation of the boundary condition, which plays a significant role in defining the behavior of electromagnetic fields at the domain boundaries. 
These fundamental components are essential for the successful application of FieldTNN in solving Maxwell eigenvalue problem.

\subsubsection{Divergence-free condition}\label{subsec:div_free}
We notice that the solution $(\lambda, \mathbf{E})$ of the variational form \eqref{eq:vform2} can not guarantee that it is an eigenpair of the Maxwell eigenvalue problem \eqref{eq:Maxwelleig}. 
The eigenpairs that satisfy the divergence-free condition are named as ``real'' eigenpairs \cite{CMAME2008}.
The reason is that solutions to the variational form \eqref{eq:vform2} will include non-divergence-free ones, which are termed as spurious eigenpairs \cite{CMAME2008}.
Interestingly, provided the spurious eigenvalue is simple, the associated eigenvectors are automatically curl-free \cite{Costabel2003}.
Thus, in this subsection, we will present an approach to deal with divergence-free condition and then remove the spurious eigenpairs in the numerical results.
Actually, Ciarlet and Hechme proposed three computational strategies to remove spurious eigenmodes, interested readers can refer to \cite{CMAME2008} for more details. 
Based on the work of Ciarlet and Hechme, we consider the following fileter ratio $\rho$:
\begin{align*}
&\rho = \frac{\left| \mathbf{E}_{\scriptscriptstyle \text{NN}}\right|_{\mathbf{H}(\operatorname{div};\Omega)}^{2}}{\left| \mathbf{E}_{\scriptscriptstyle \text{NN}}\right|_{\mathbf{H}(\operatorname{curl};\Omega)}^{2}}.
\end{align*}
Actually, for the ``real'' eigenpairs, the eigenfunctions $\mathbf{E}$ is divergence-free, which implies that the divergence of the numerical eigenfunction $\mathbf{E}_{\scriptscriptstyle \text{NN}}$ should be relatively small. 
Hence, the numerator of fileter ratio $\rho$ will approach to $0$, resulting in a small ratio value (may be close to 0). 
In reverse, for the spurious eigenpair, since the eigenfunction is curl-free, the denominator of fileter ratio $\rho$ will approach to $0$, resulting in a large ratio value. 
Thus, we can specify the definition of single loss function \eqref{eq:singleLoss} used in this paper as follows:
\begin{align*}
&\mathcal{L}_{\scriptscriptstyle \text{Maxwell}}\left(\lambda_{k, \scriptscriptstyle \text{NN}}, \mathbf{E}_{k, \scriptscriptstyle \text{NN}}\right) = \lambda_{k, \scriptscriptstyle \text{NN}} + \beta \frac{\left| \mathbf{E}_{k, \scriptscriptstyle \text{NN}}\right|_{\mathbf{H}(\operatorname{div};\Omega)}^{2}}{\left| \mathbf{E}_{k, \scriptscriptstyle \text{NN}}\right|_{\mathbf{H}(\operatorname{curl};\Omega)}^{2}},~k=1,\cdots,p.
\end{align*}
For computing the leading $\mathcal{M}$ eigenpairs with $\mathcal{M}\leq p$, we can further define the following loss function (cf. \eqref{eq:totalLoss}),
\begin{align}
&\widehat{\mathcal{L}}\left(\mathbf{\Theta}\right)=\sum_{m=1}^{\mathcal{M}} \mathcal{L}_{\scriptscriptstyle \text{Maxwell}}\left(\lambda_{m, \scriptscriptstyle \text{NN}}, \mathbf{E}_{m, \scriptscriptstyle \text{NN}}\right).\label{eq:Loss}
\end{align}

\subsubsection{Boundary condition}\label{subsec:BC}
In the computation of Maxwell equation and its associated eigenvalue problem, one of common used boundary conditions is the perfect electric conductor (PEC) boundary condition:
\begin{align*}
&\mathbf{E} \times \mathbf{n}=\mathbf{0},\ \text{on} \ \partial \Omega.
\end{align*}
Taking the 3D unit cube cavity $\Omega=\left[a_{1}, b_{1}\right] \times \left[a_{2}, b_{2}\right] \times \left[a_{3}, b_{3}\right]$ as an example, the boundary of the domain $\Omega$ can be decomposed into the following $6$ faces:
\begin{align*}
\begin{array}{ll}
\Gamma_{1, a}:=  \left\{x_{1}=a_{1}\right\} \times\left(a_{2}, b_{2}\right) \times\left(a_{3}, b_{3}\right),&\Gamma_{1, b}:=  \left\{x_{1}=b_{1}\right\} \times\left(a_{2}, b_{2}\right) \times\left(a_{3}, b_{3}\right), \\
\Gamma_{2, a}:= \left(a_{1}, b_{1}\right) \times\left\{x_{2}=a_{2}\right\} \times\left(a_{3}, b_{3}\right),&\Gamma_{2, b}:=  \left(a_{1}, b_{1}\right) \times\left\{x_{2}=b_{2}\right\} \times\left(a_{3}, b_{3}\right), \\
\Gamma_{3, a}:= \left(a_{1}, b_{1}\right) \times\left(a_{2}, b_{2}\right) \times\left\{x_{3}=a_{3}\right\},&\Gamma_{3, b}:=  \left(a_{1}, b_{1}\right) \times\left(a_{2}, b_{2}\right) \times\left\{x_{3}=b_{3}\right\}.
\end{array}
\end{align*}
Then the corresponding $6$ outward normal vectors to boundaries $\Gamma_{i, a}$ and $\Gamma_{i, b}~(i=1,2,3)$ are
\begin{equation*}
\begin{array}{ll}
\mathbf{n}_{1}=(1,0,0), & -\mathbf{n}_{1}=(-1,0,0), \\
\mathbf{n}_{2}=(0,1,0), & -\mathbf{n}_{2}=(0,-1,0), \\
\mathbf{n}_{3}=(0,0,1), & -\mathbf{n}_{3}=(0,0,-1).
\end{array}
\end{equation*}
Obviously, we need to satisfy the boundary condition $\mathbf{E} \times \mathbf{n} = 0$ on the six faces mentioned above ($\Gamma_{i, a}$ and $\Gamma_{i, b}$,~$i=1,2,3$). 
We will specifically discuss the implementation of the boundary condition on $\Gamma_{1, a}$, and the remaining ones can be extended accordingly. On the face $\Gamma_{1, a}$, we have
\begin{equation*}
\mathbf{E} \times \mathbf{n}_{1}=\left(\begin{array}{l}
\mathrm{E}_{1} \\
\mathrm{E}_{2} \\
\mathrm{E}_{3}
\end{array}\right) \times\left(\begin{array}{l}
1 \\
0 \\
0
\end{array}\right)=\left|\begin{array}{ccc}
\mathbf{e}_{1} & \mathbf{e}_{2} & \mathbf{e}_{3} \\
\mathrm{E}_{1} & \mathrm{E}_{2} & \mathrm{E}_{3} \\
1 & 0 & 0
\end{array}\right|=\left(\begin{array}{c}
~0\\
~~\mathrm{E}_{3}\\
-\mathrm{E}_{2}
\end{array}\right),
\end{equation*}
Then the boundary condition $\mathbf{E} \times \mathbf{n}_{1} = 0$ implies that the following constraints to the face $\Gamma_{1, a}$:
\begin{align*}
&\mathrm{E}_{3}\left(x_{1}=a_{1}, x_{2}, x_{3}\right)=0,~\mathrm{E}_{2}\left(x_{1}=a_{1}, x_{2}, x_{3}\right)=0.
\end{align*}
The discussion on the $6$ faces yields $12$ constraints.
Thus, we construct a function that satisfies the boundary conditions by multiplying analytic functions $\gamma_{i}(x_{i})$, which satisfy $\gamma_{i}(a_{i}) = 0$, $\gamma_{i}(b_{i}) = 0$ and $\gamma_{i}(x_{i}) \neq 0$ when $x_{i} \in \left(a_{i},b_{i}\right)$, i.e.,
\begin{equation}\label{eq:FTNNtimesBC}
\boldsymbol{\widehat{\Phi}}\left(\mathbf{x}; \mathbf{\Theta}\right) = \sum_{k=1}^{p}\alpha_{k}\left(\begin{array}{c}
\widehat{\Phi}_{1,k}\left(\mathbf{x}; \Theta_{1}\right) \\
\widehat{\Phi}_{2,k}\left(\mathbf{x}; \Theta_{2}\right) \\
\widehat{\Phi}_{3,k}\left(\mathbf{x}; \Theta_{3}\right)
\end{array}\right)=\left(\begin{array}{l}
\widehat{\Psi}_{1,k}\left(\mathbf{x}; \Theta_{1}\right) \gamma_{2}\left(x_{2}\right) \gamma_{3}\left(x_{3}\right) \\
\widehat{\Psi}_{2,k}\left(\mathbf{x}; \Theta_{2}\right) \gamma_{1}\left(x_{1}\right) \gamma_{3}\left(x_{3}\right) \\
\widehat{\Psi}_{3,k}\left(\mathbf{x}; \Theta_{3}\right) \gamma_{1}\left(x_{1}\right) \gamma_{2}\left(x_{2}\right)
\end{array}\right).
\end{equation}
Obviously, $\boldsymbol{\widehat{\Phi}}\left(\mathbf{x}; \mathbf{\Theta}\right)$ satisfies the boundary condition in \eqref{eq:Maxwelleig}, it follows that
\begin{align*}
&\boldsymbol{\widehat{\Phi}}\left(\mathbf{x}; \mathbf{\Theta}\right) \in  \mathbf{H}_{0}(\operatorname{curl};\Omega) \cap \mathbb{V}_{p}.
\end{align*}

\subsection{Algorithm}\label{sec:algorithm}
In this section, we will present complete FieldTNN-based machine learning algorithm for computing leading $\mathcal{M}$ $(\leq p)$ Maxwell eigenpairs. 
The algorithm for the tensor computational domains is summarized in Algorithm \ref{algorithm:a1}, while the one for non-tensor domains is detailed in Algorithm \ref{algorithm:a2}.
\begin{algorithm}[!ht]
\caption{Solving Maxwell eigenvalue problem by FieldTNN in tensor domain.}
\label{algorithm:a1}
\begin{algorithmic} 
\STATE \textbf{Step} $\mathbf{1:}$ Initialization: Construct the initial FieldTNN $\boldsymbol{\widehat{\Phi}}\left(\mathbf{x}; \mathbf{\Theta}\right)$ as defined in \eqref{eq:FTNNtimesBC}, set maximum training steps $\mathcal{H}$, learning rate $\eta$, quadrature points $\left\{x_i^{\left(\ell_{i}\right)}\right\}_{\ell_{i}=1}^{M_i N_i}$ and quadrature weights $\left\{w_i^{\left(\ell_{i}\right)}\right\}_{\ell_{i}=1}^{M_i N_i}$ $(i=1,2,3)$ defined in \eqref{eq:LGPW}.
\STATE \textbf{Step} $\mathbf{2:}$ Assemble matrices $\mathrm{S}$ and $\mathrm{M}$ according to quadrature schemes \eqref{eq:StiffnessGQ} and \eqref{eq:MassGQ}. 
\STATE \textbf{Step} $\mathbf{3:}$ Compute the value of the loss function \eqref{eq:Loss}.
\STATE \textbf{Step} $\mathbf{4:}$ Compute the gradient of the loss function with respect to parameters $\mathbf{\Theta}$ by automatic differentiation and update parameters $\mathbf{\Theta}$ by \eqref{eq:update}.
\STATE \textbf{Step} $\mathbf{5:}$ Set $\hbar=\hbar+1$ and go to Step 2 for the next step until $\hbar = \mathcal{H}$.
\STATE \textbf{Step} $\mathbf{6:}$ Post-processing step: Solve the generalized matrix eigenvalue problem \eqref{eq:GMEP} to compute the leading $\mathcal{M}$ eigenpairs.
\end{algorithmic}
\end{algorithm}

\begin{algorithm}[!ht]
\caption{Solving Maxwell eigenvalue problem by FieldTNN in non-tensor domain.}
\label{algorithm:a2}
\begin{algorithmic} 
\STATE \textbf{Step} $\mathbf{1:}$ Domain decomposition: Decompose the computational domain $\Omega$ into several tensor subdomains $\Omega_{\ell}$.
\STATE \textbf{Step} $\mathbf{1:}$ Initialization: Construct the initial FieldTNN $\boldsymbol{\Phi}^{\ell}\left(\mathbf{x}; \mathbf{\Theta}\right)$, which defined in \eqref{eq:nonTensorFieldTNN}, for each subdomain $\Omega_{\ell}$, set maximum training steps $\mathcal{H}$, learning rate $\eta$, quadrature points $\left\{x_i^{\left(\ell_{i}\right)}\right\}_{\ell_{i}=1}^{M_i N_i}$ and quadrature weights $\left\{w_i^{\left(\ell_{i}\right)}\right\}_{\ell_{i}=1}^{M_i N_i}$ $(i=1,2,3)$ defined in \eqref{eq:LGPW}.
\STATE \textbf{Step} $\mathbf{2:}$ Assemble submatrices $\mathrm{S}_{\ell}$ and $\mathrm{M}_{\ell}$ according to quadrature schemes \eqref{eq:StiffnessGQ} and \eqref{eq:MassGQ} to form block matrices like \eqref{eq:SPMatrix}. 
\STATE \textbf{Step} $\mathbf{3:}$ Compute the value of the loss function \eqref{eq:Loss}.
\STATE \textbf{Step} $\mathbf{4:}$ Compute the gradient of the loss function with respect to parameters $\mathbf{\Theta}_{\ell}$ by automatic differentiation and update parameters $\mathbf{\Theta}_{\ell}$ by \eqref{eq:update}.
\STATE \textbf{Step} $\mathbf{5:}$ Set $\hbar = \hbar + 1$ and go to Step $3$ for the next step until $\hbar = \mathcal{H}$.
\STATE \textbf{Step} $\mathbf{6:}$ Post-processing step: Solve the generalized matrix eigenvalue problem to compute the first $\mathcal{M}$ eigenpairs.
\end{algorithmic}
\end{algorithm}
\section{Numerical examples}\label{sec:Numerical}
In this section,  we provide several examples to investigate the performance of the  FieldTNN-based machine learning method for Maxwell eigenvalue problems proposed in this paper.
These examples demonstrate that our approach successfully eliminates spurious eigenpairs, 
achieves the remarkable accuracy, 
and ensures that the resulting eigenfunctions satisfy the divergence-free condition. 
In order to illustrate the accuracy and convergence of our method, the relative error of the eigenvalues is defined as follows:
\begin{align*}
&\text{err}_{\lambda, k} := \frac{|\lambda_{k, \scriptscriptstyle \text{NN}} - \lambda_{k}|}{\lambda_{k}}.
\end{align*}
where $k = 1,2,\cdots,\mathcal{M}$ and $\lambda_{k}$ denotes the $k$-th reference eigenvalue.

\subsection{Square cavity}
In the first example, we consider the computational domain to be a square cavity $\Omega = [0,1]^{2}$ and assume $\varepsilon = \mu = 1$. The exact eigenvalues are given by
\begin{align*}
&\lambda_{ij} = \left(i^2 + j^2\right)\pi^{2},
\end{align*} 
where $i, j=0,1,2, \cdots$ and $i + j>0$. These eigenvalues correspond to the fundamental frequencies of the Maxwell equations in a square cavity \cite{BoffiMixed}. The associated eigenfunctions are in the form
\begin{align*}
&\mathbf{E}_{ij} = \Big(-j \cos (i \pi x_{1}) \sin (j \pi x_{2}),~i \sin (i \pi x_{1}) \cos (j \pi x_{2})\Big),~i, j=0,1,2, \cdots.
\end{align*}
Each subnetwork of the FiledTNN is an FNN with  $3$ hidden layers and each hidden layer has $100$ neurons.
The activation function is selected as sine function and the rank $p$ is set to be $50$.
The Adam optimizer is used with an initial learning rate of $0.0001$ for the first $100000$ epochs, followed by L-BFGS optimization for $5000$ steps with a learning rate of $0.1$ to refine the final result. 
For numerical integration, Legendre-Gauss quadrature scheme with $1600$ points per dimension is applied.

\begin{figure}[htbp]
	\centering
	\includegraphics[width=.3\textwidth]{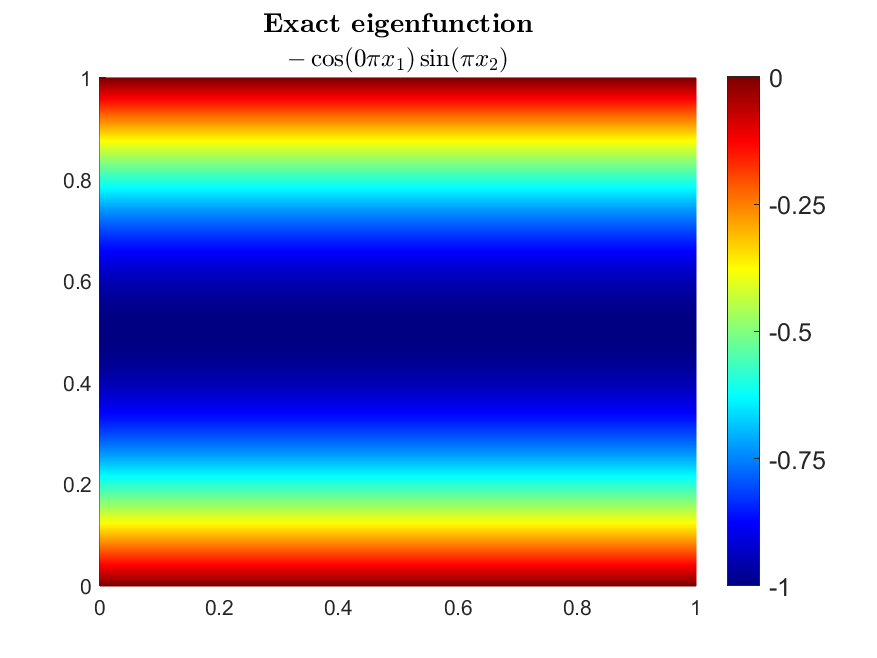}~
	\includegraphics[width=.3\textwidth]{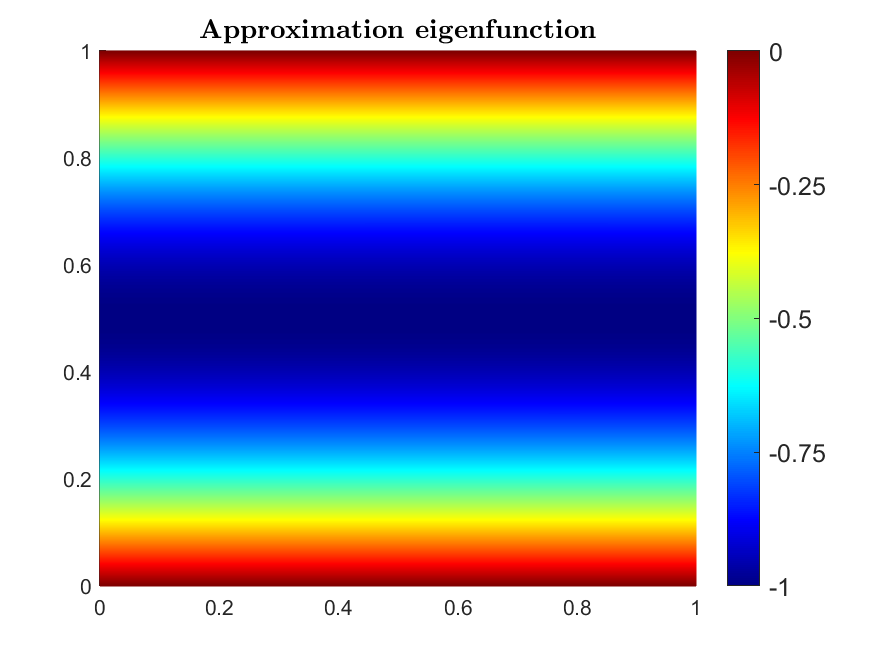}~
	\includegraphics[width=.3\textwidth]{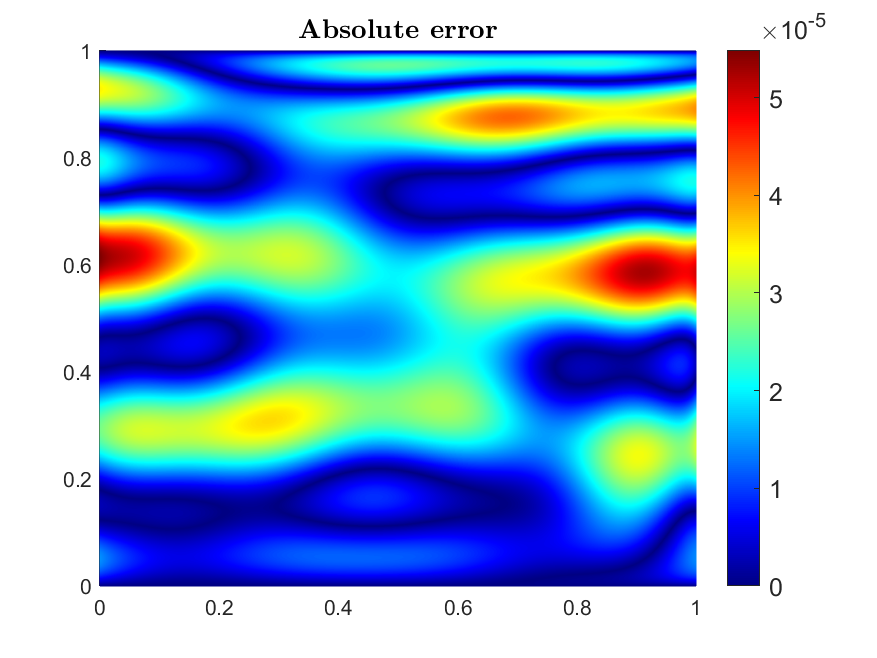}\\
	\includegraphics[width=.3\textwidth]{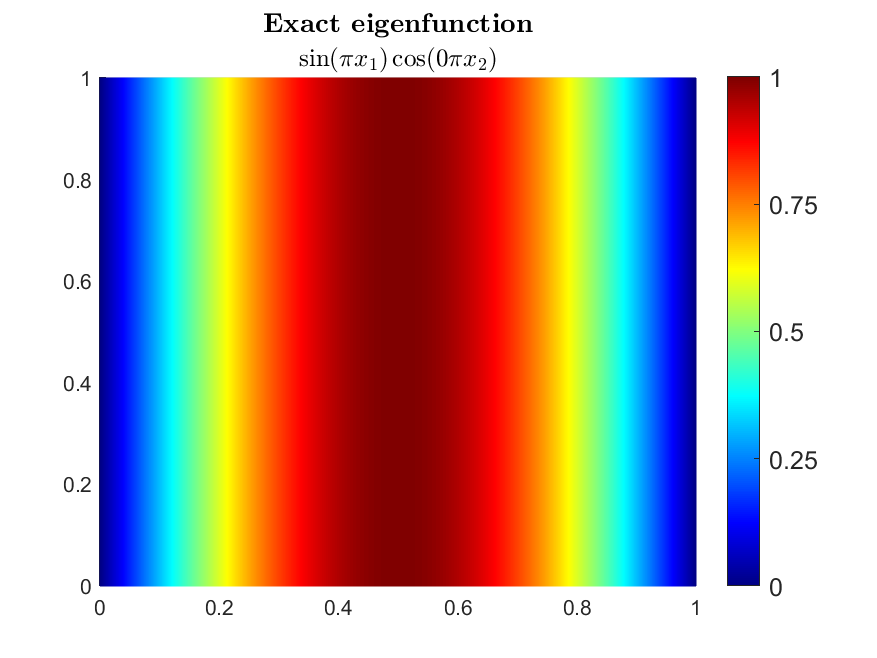}~
	\includegraphics[width=.3\textwidth]{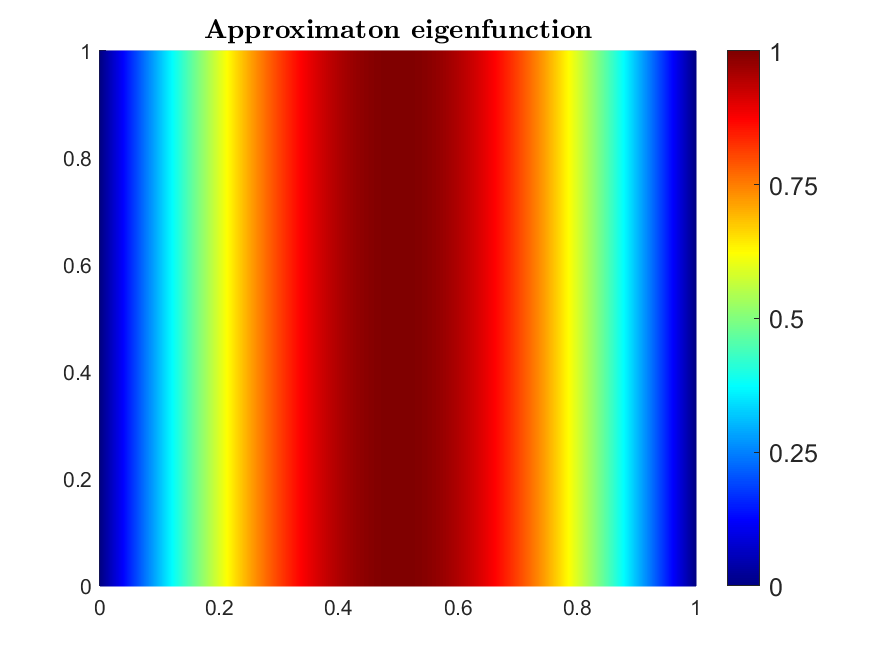}~
	\includegraphics[width=.3\textwidth]{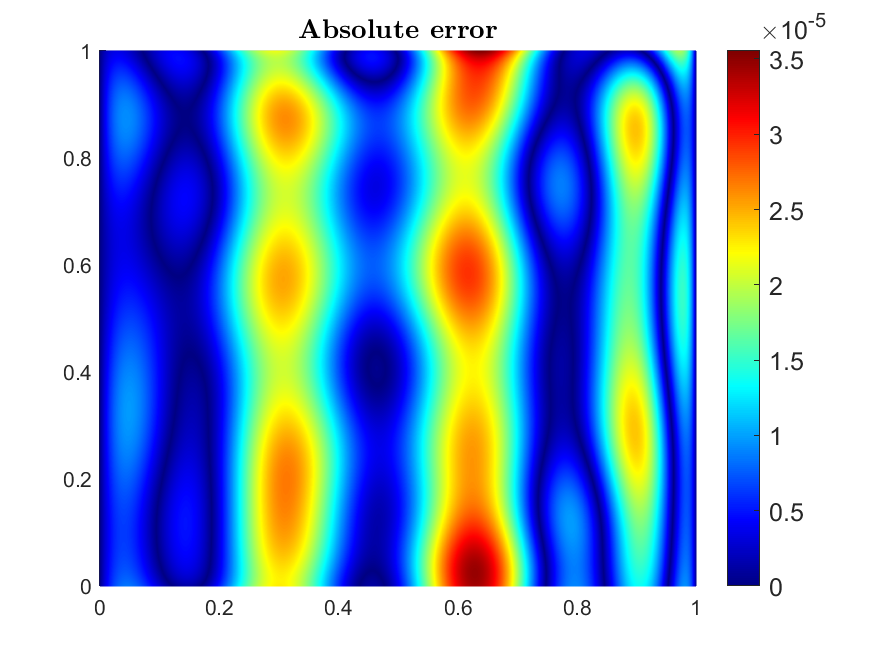}\\	
	\includegraphics[width=.3\textwidth]{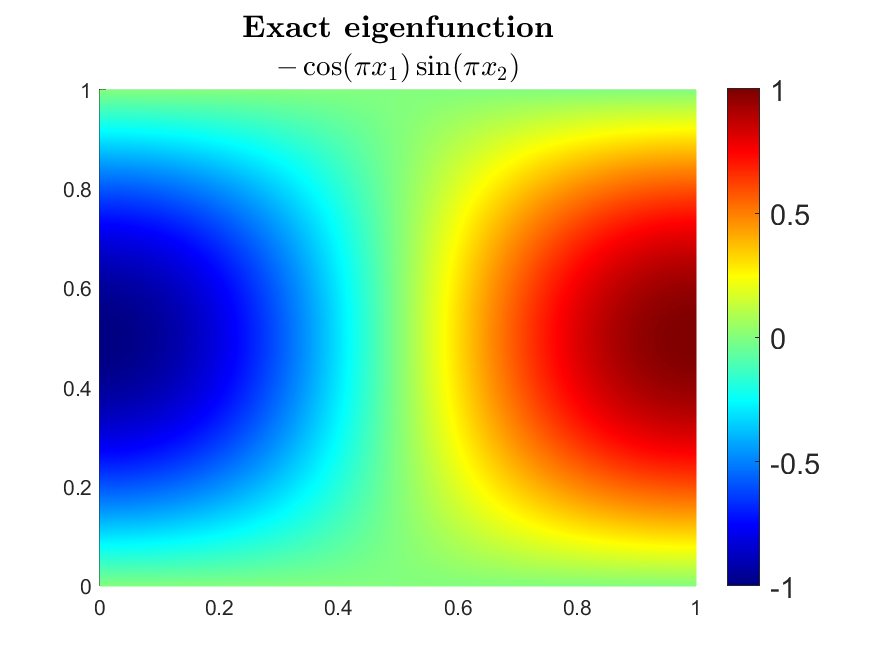}~
	\includegraphics[width=.3\textwidth]{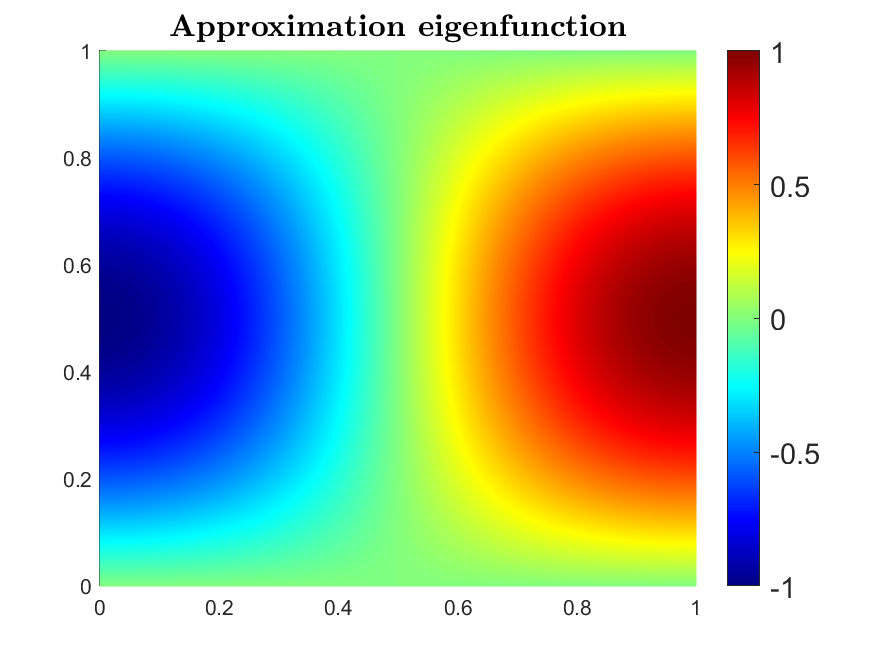}~
	\includegraphics[width=.3\textwidth]{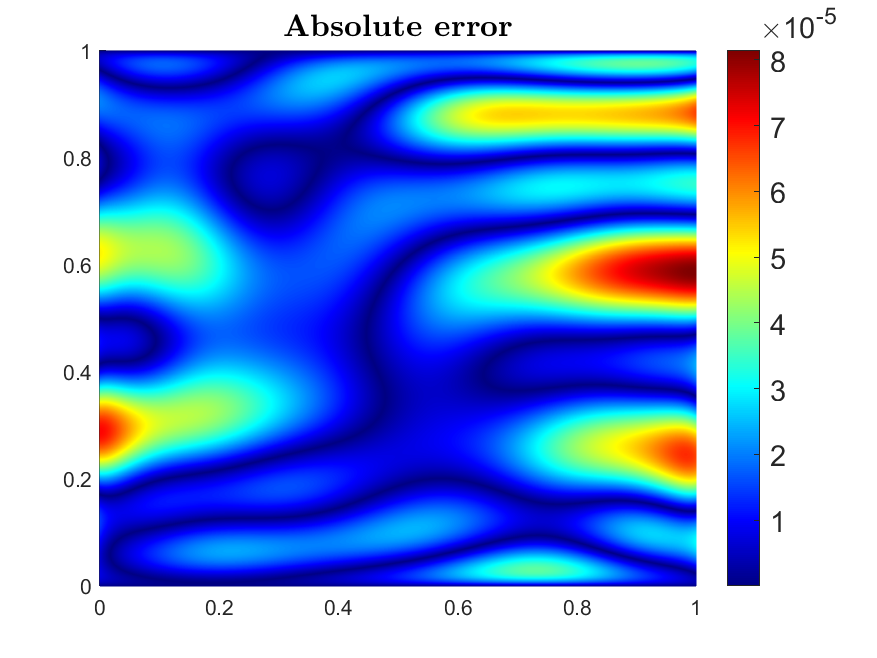}\\
	\includegraphics[width=.3\textwidth]{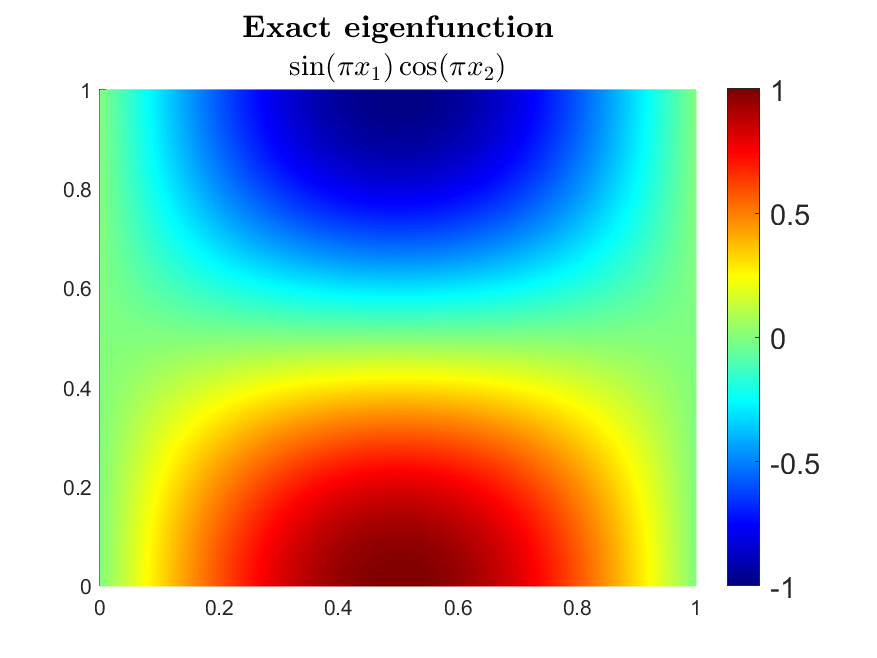}~
	\includegraphics[width=.3\textwidth]{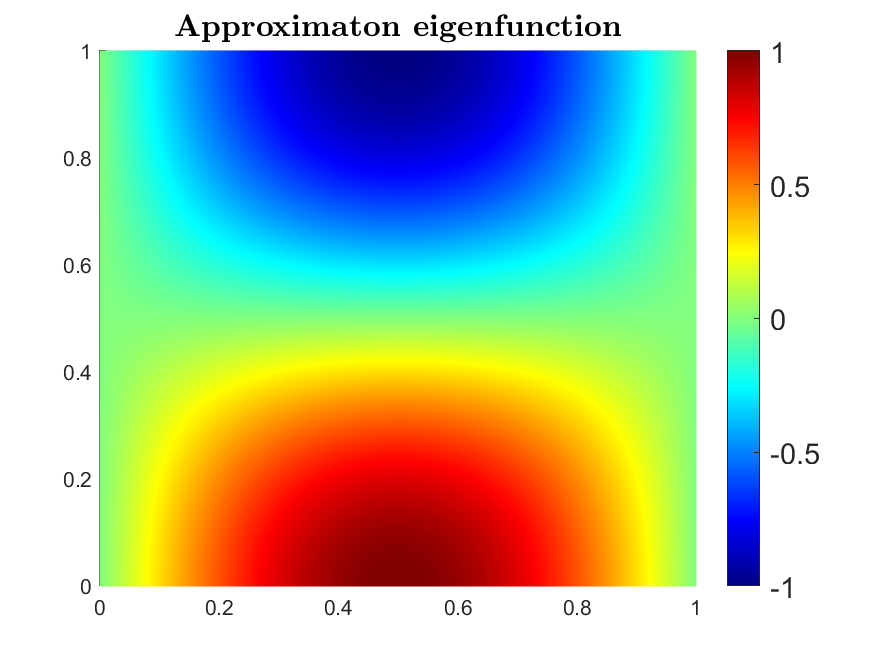}~
	\includegraphics[width=.3\textwidth]{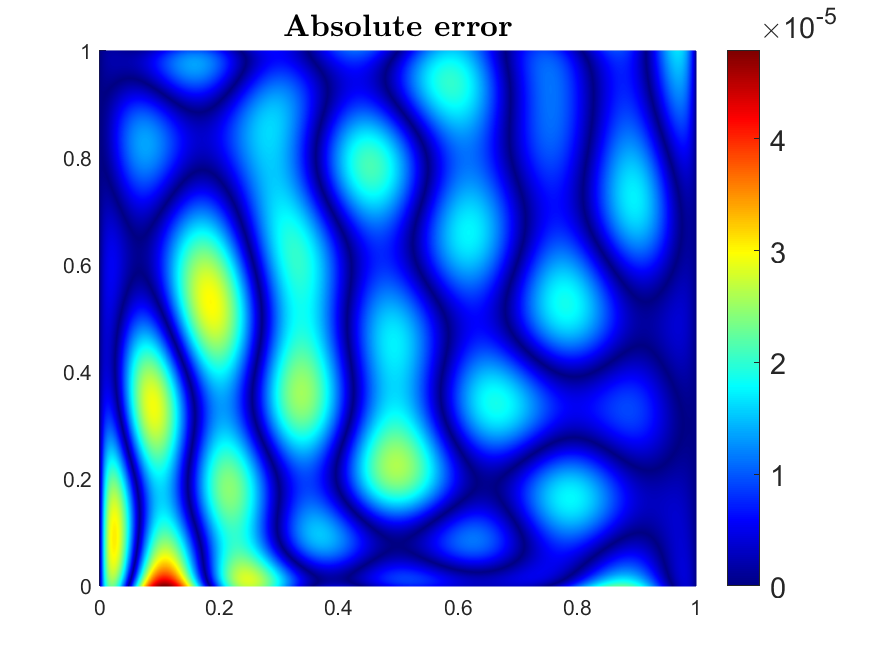}\\
	\includegraphics[width=.3\textwidth]{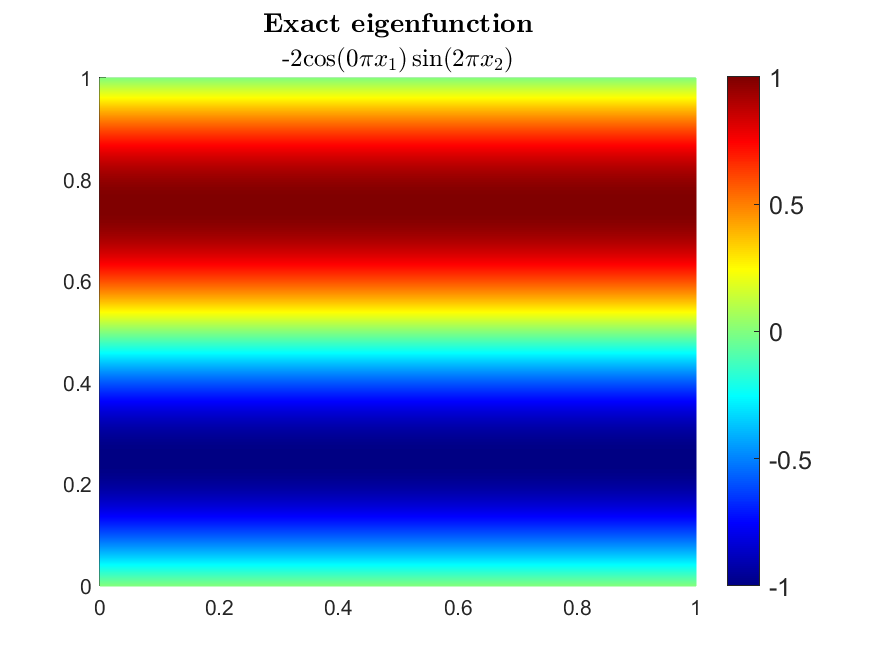}~
	\includegraphics[width=.3\textwidth]{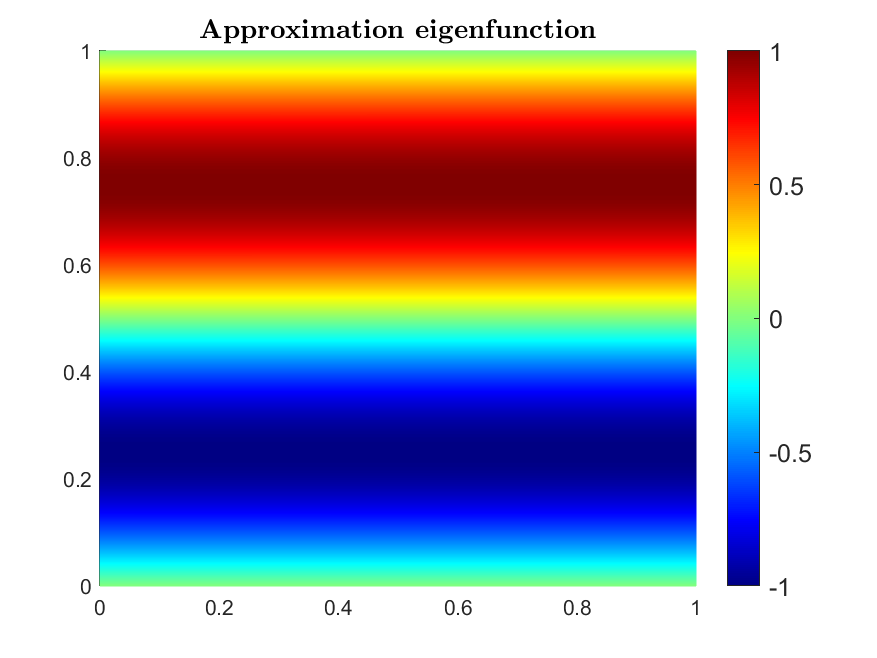}~
	\includegraphics[width=.3\textwidth]{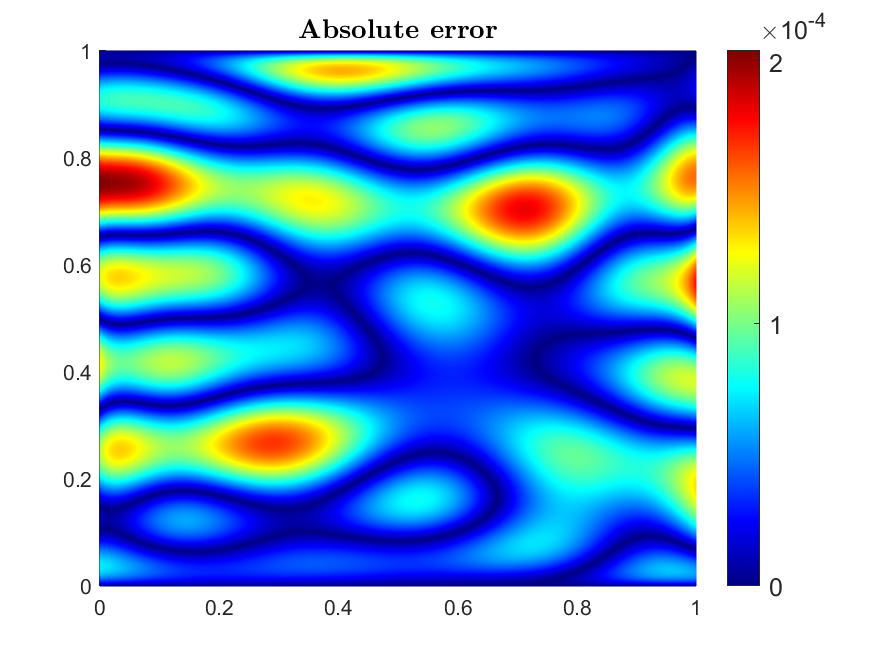}\\
	\includegraphics[width=.3\textwidth]{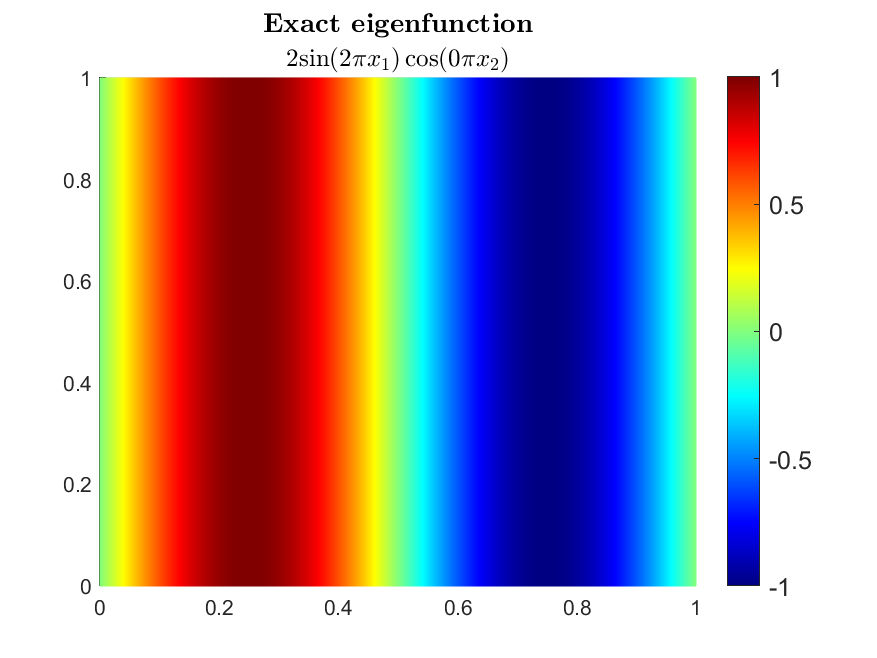}~
	\includegraphics[width=.3\textwidth]{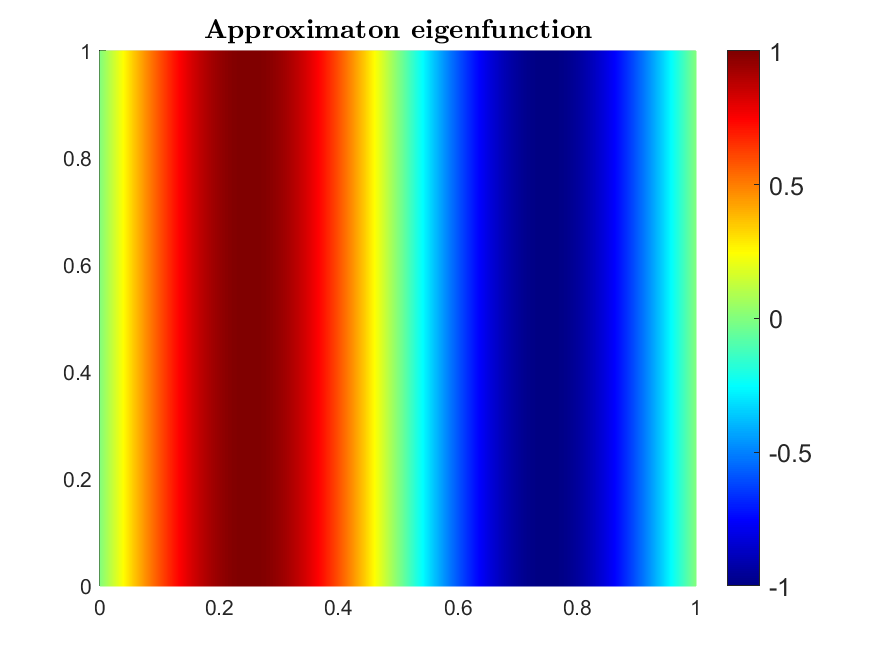}~
	\includegraphics[width=.3\textwidth]{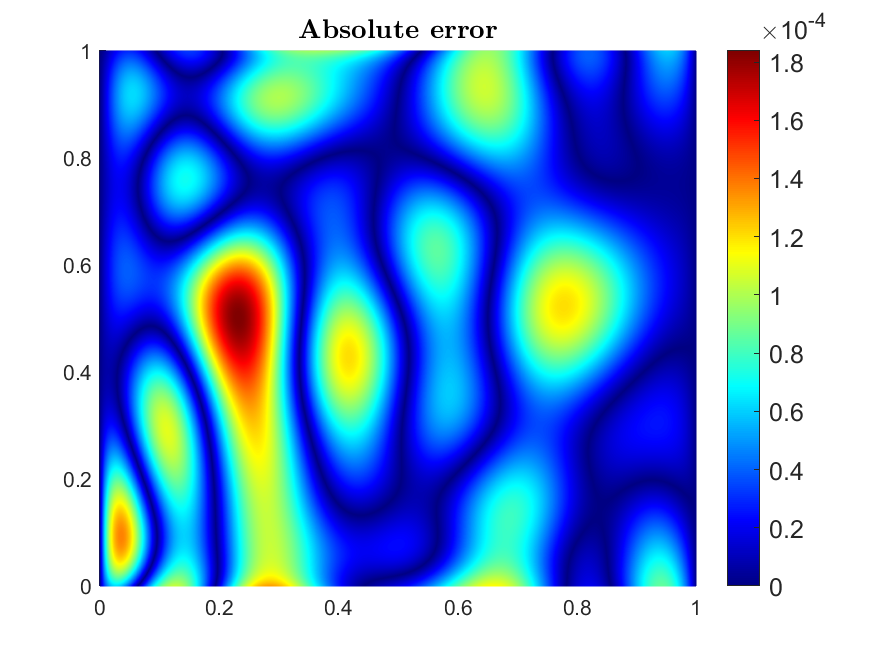}
\caption{Maxwell eigenfunctions in square cavity $\Omega = [0,1]^{2}$, exaction solutions (left column), FieldTNN approximation solutions (middle column) and associated absolute errors (right column).}
\label{fig:eigenfunctions_2D.esp1}
\end{figure}

\begin{figure}[htbp]
	\centering
	\includegraphics[width=.3\textwidth]{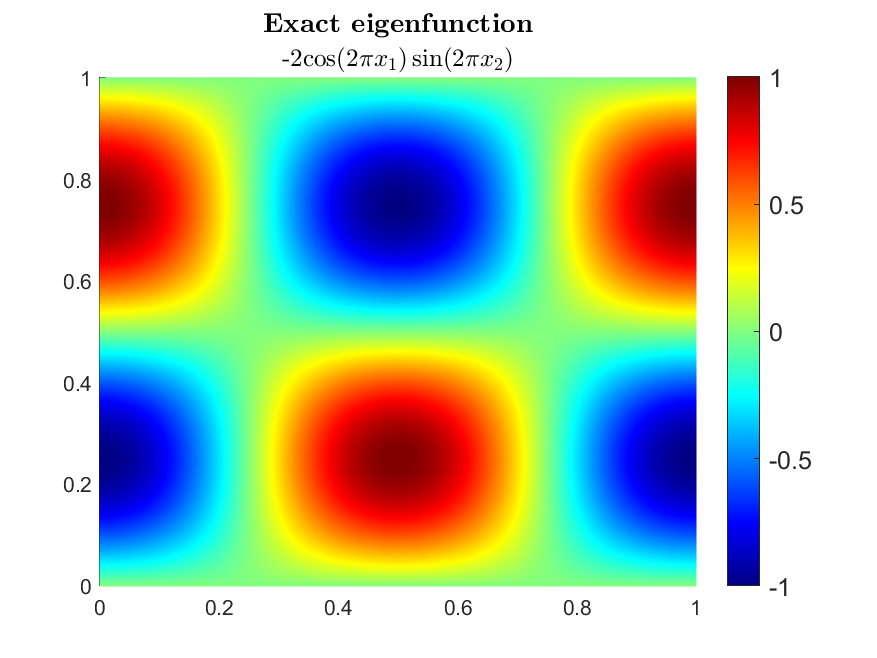}~
	\includegraphics[width=.3\textwidth]{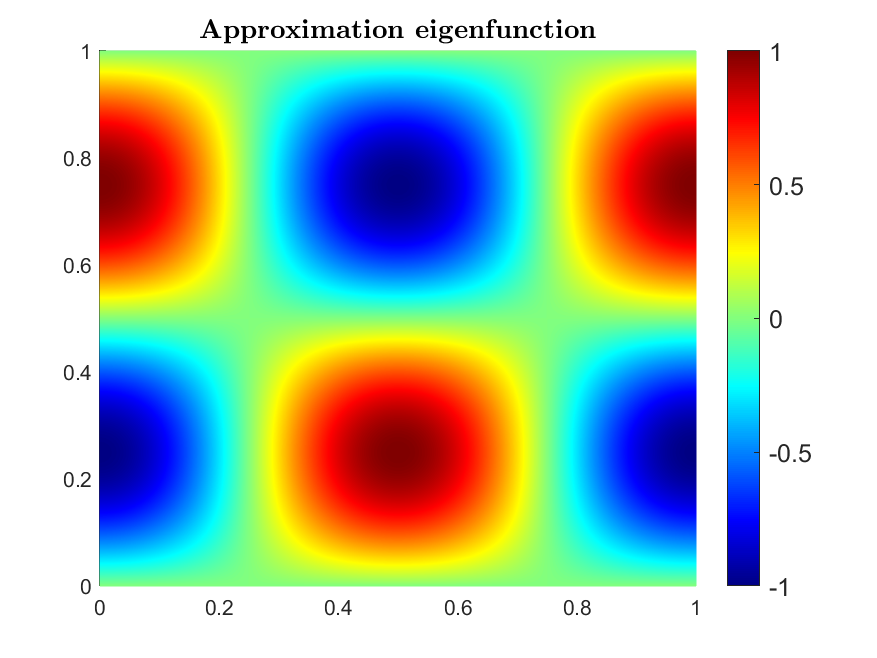}~
	\includegraphics[width=.3\textwidth]{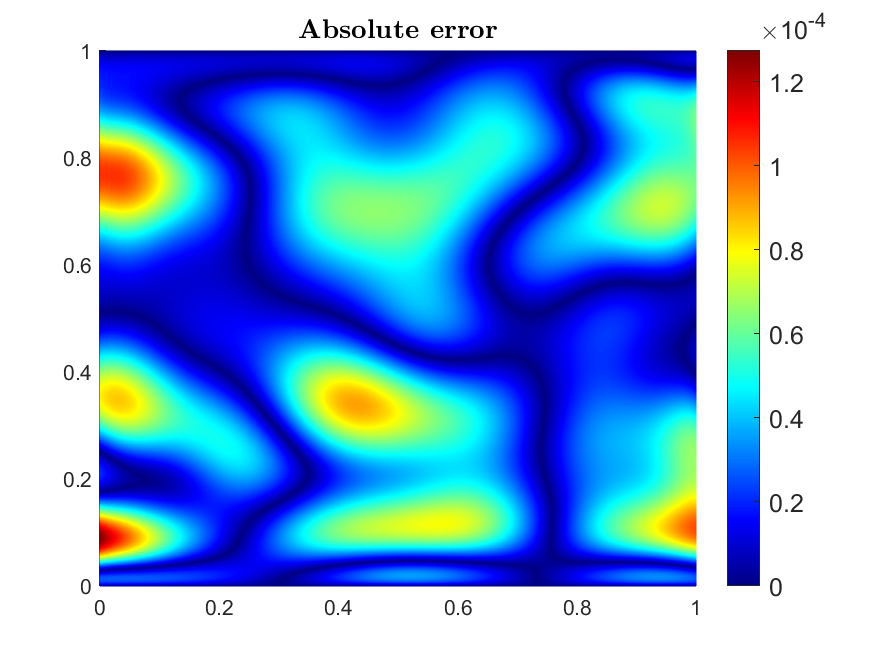}\\
	\includegraphics[width=.3\textwidth]{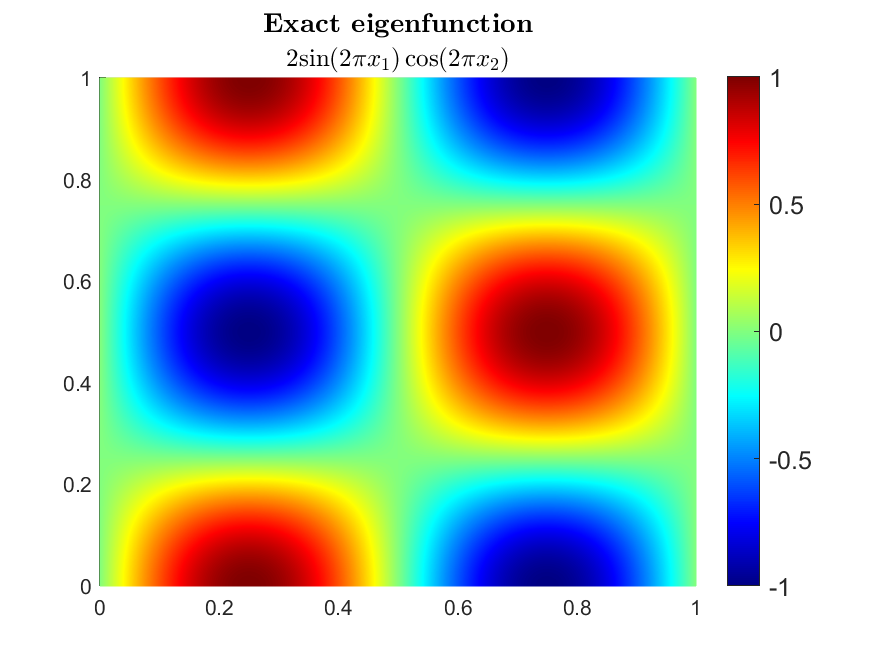}~
	\includegraphics[width=.3\textwidth]{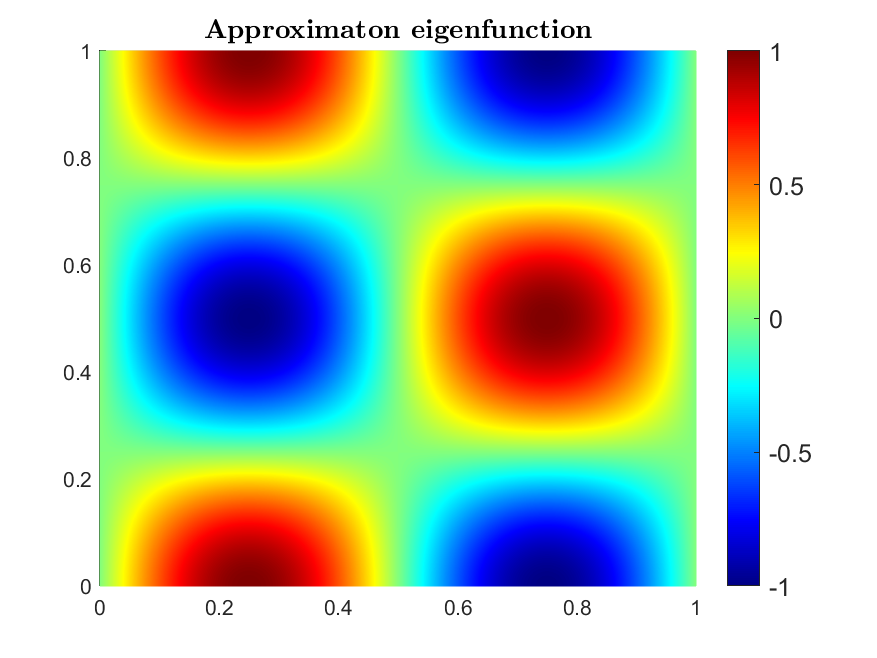}~
	\includegraphics[width=.3\textwidth]{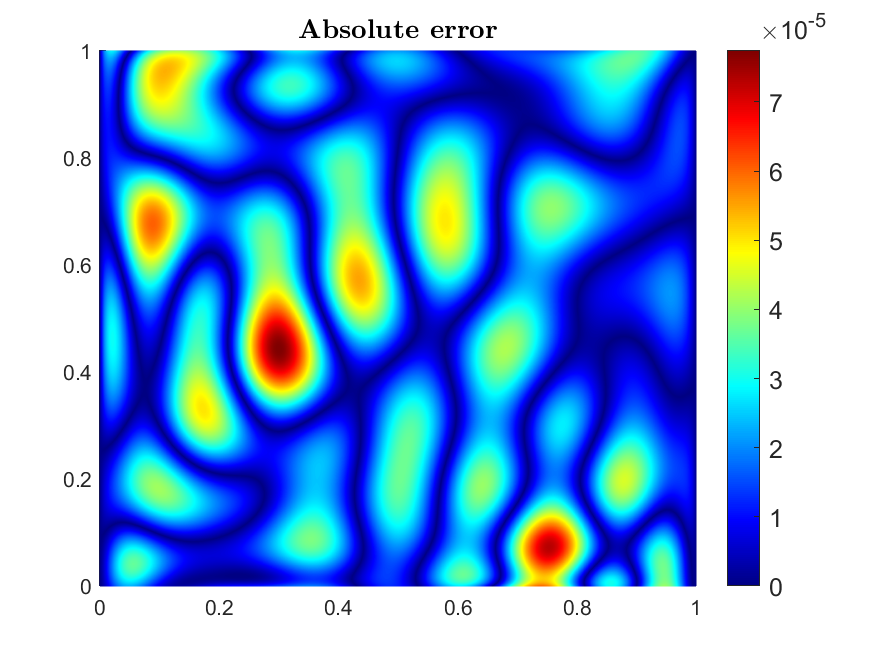}\\	
	\includegraphics[width=.3\textwidth]{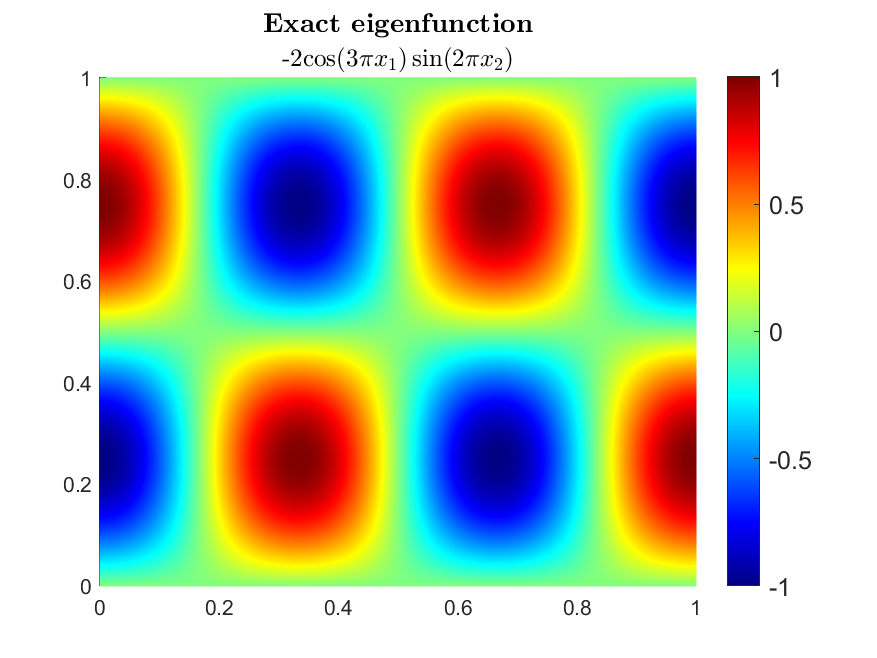}~
	\includegraphics[width=.3\textwidth]{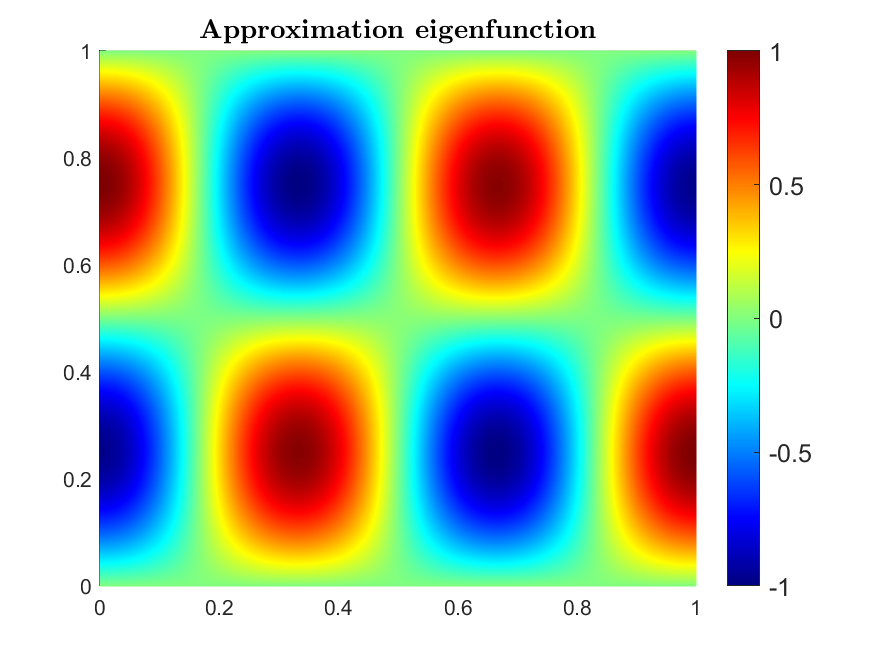}~
	\includegraphics[width=.3\textwidth]{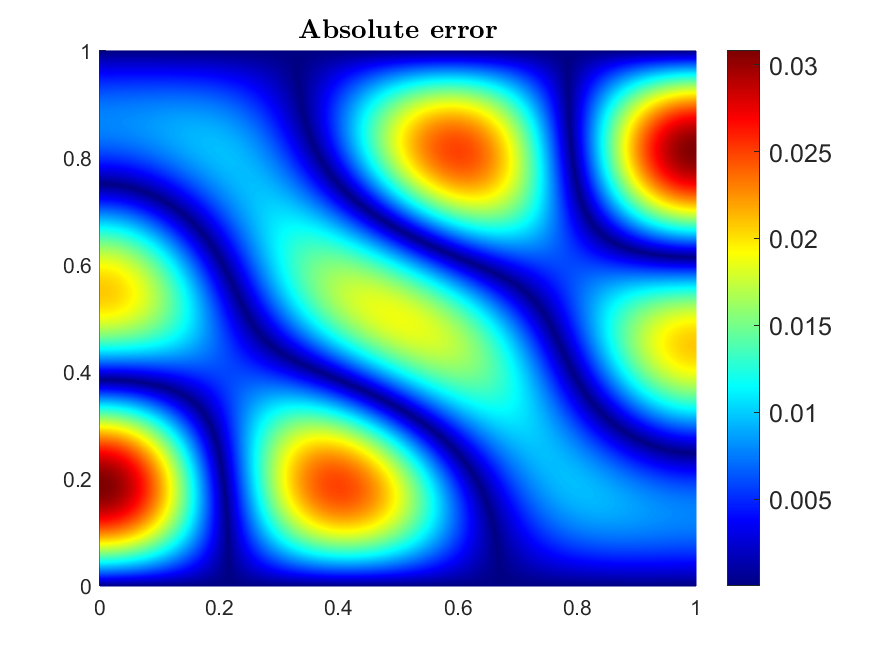}\\
	\includegraphics[width=.3\textwidth]{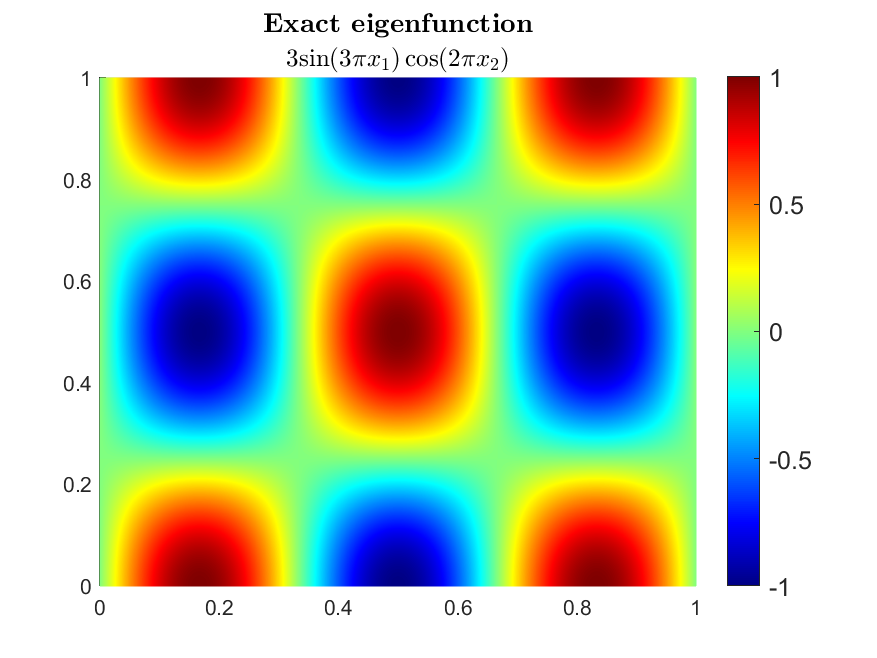}~
	\includegraphics[width=.3\textwidth]{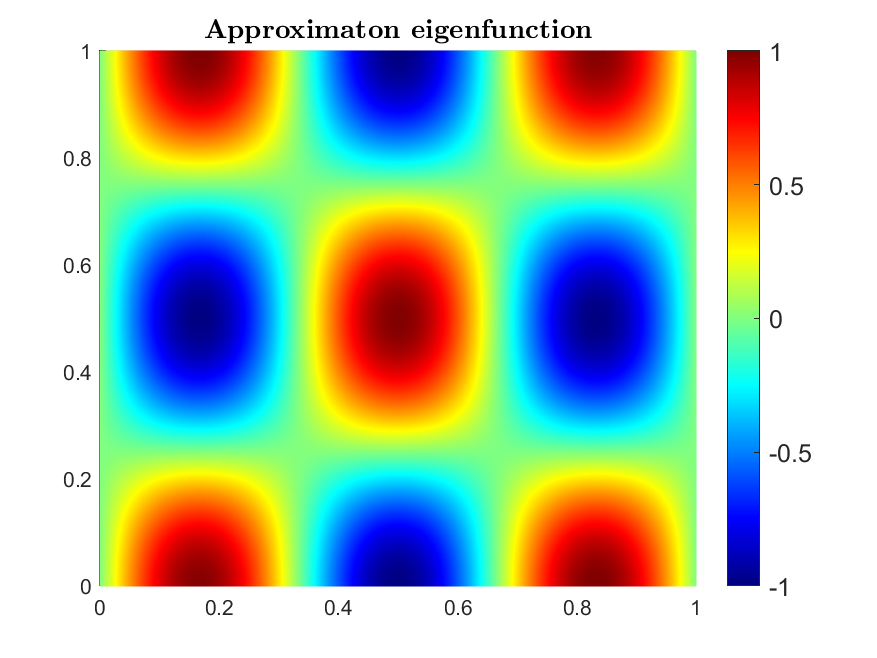}~
	\includegraphics[width=.3\textwidth]{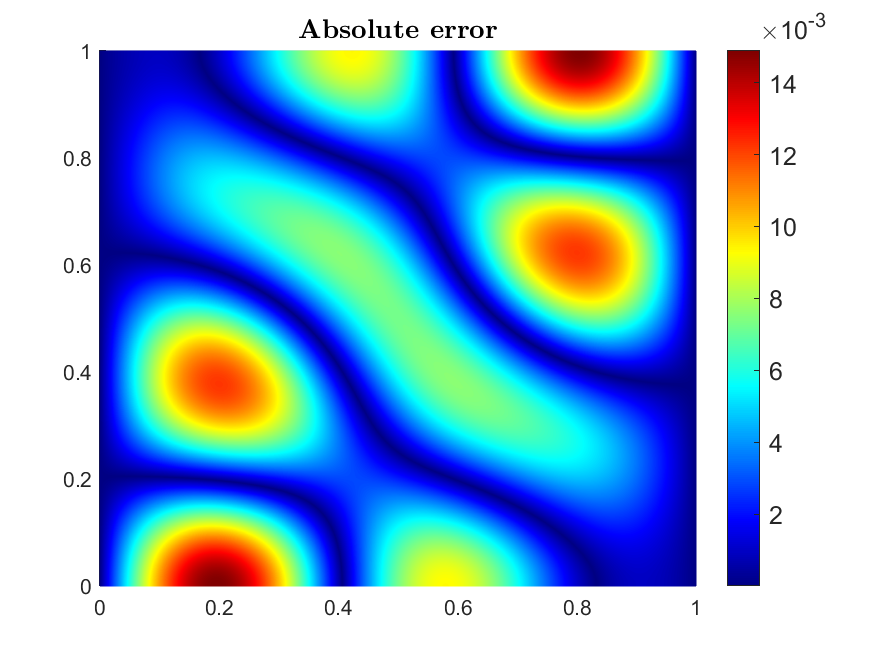}\\
	\includegraphics[width=.3\textwidth]{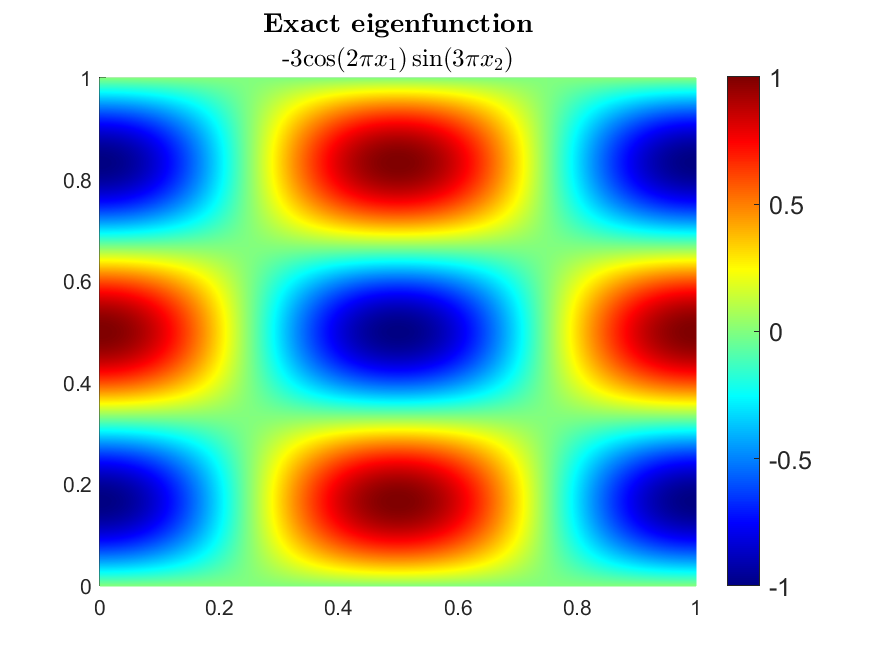}~
	\includegraphics[width=.3\textwidth]{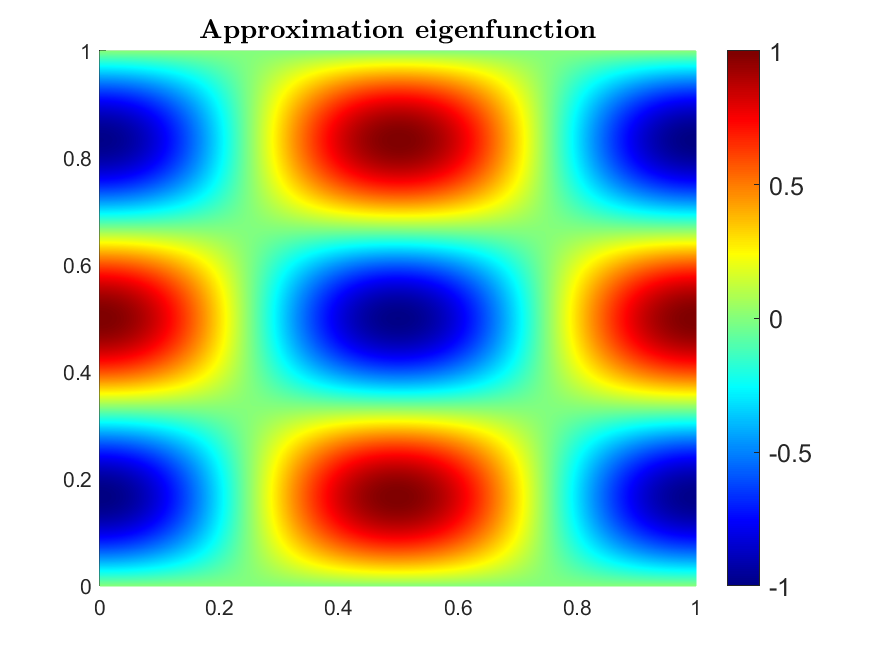}~
	\includegraphics[width=.3\textwidth]{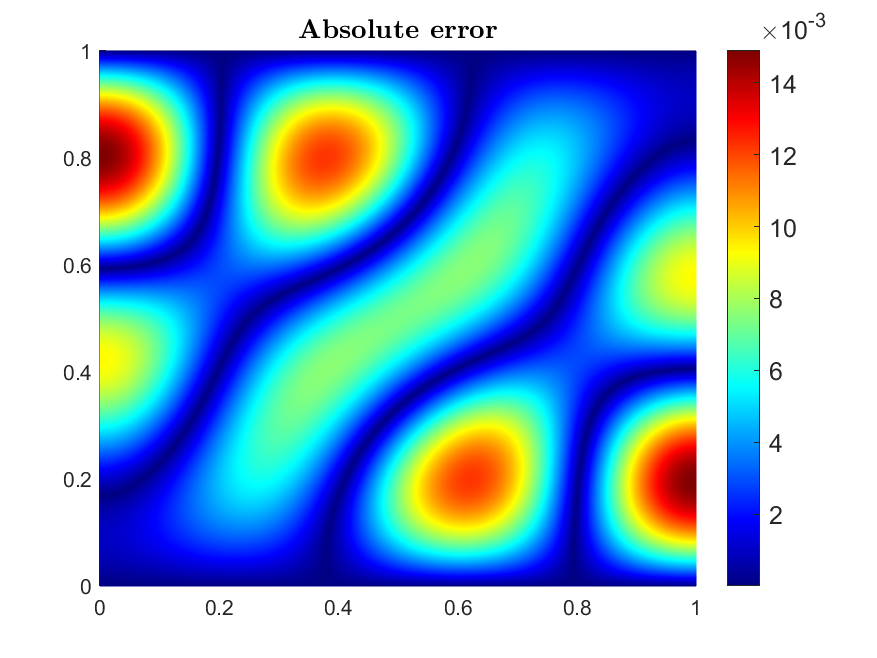}\\
	\includegraphics[width=.3\textwidth]{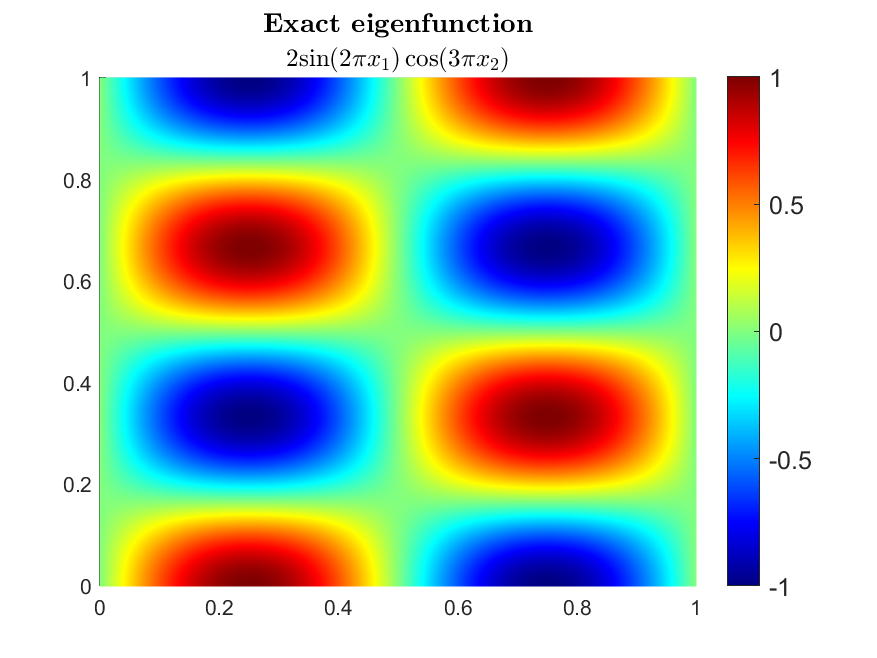}~
	\includegraphics[width=.3\textwidth]{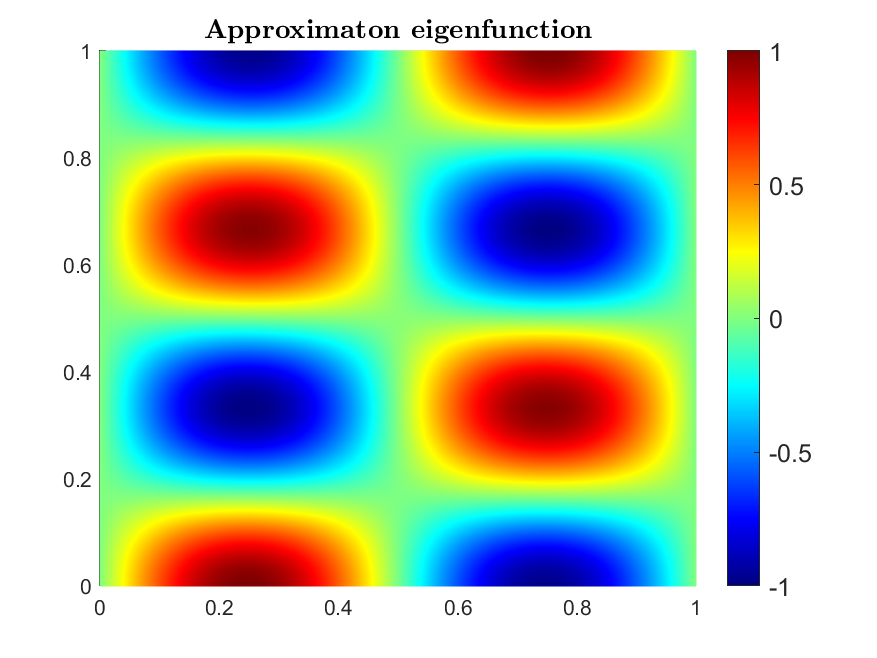}~
	\includegraphics[width=.3\textwidth]{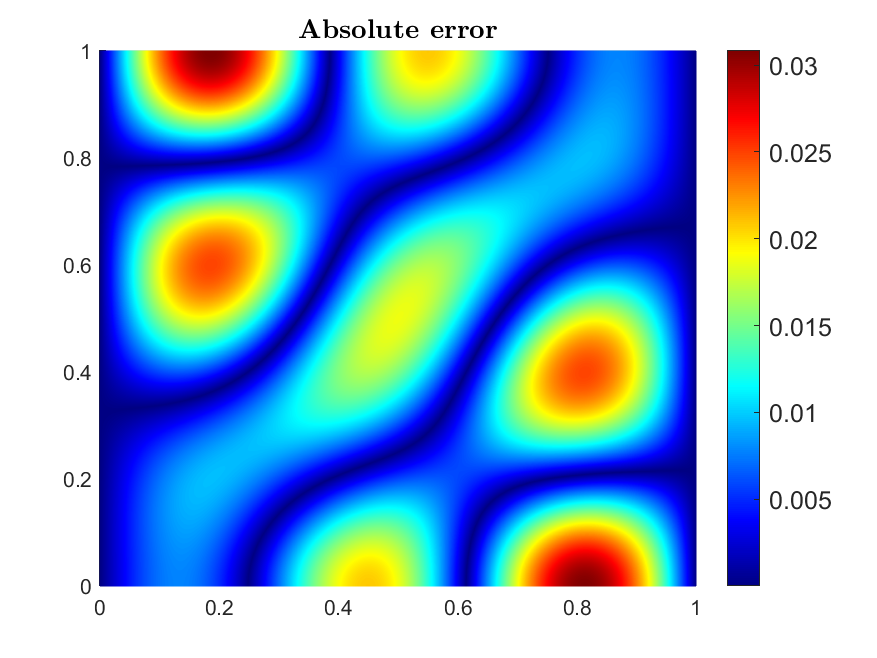}\\		
	\caption{Maxwell eigenfunctions in square cavity $\Omega = [0,1]^{2}$, exaction solutions (left column), FieldTNN approximation solutions (middle column) and associated absolute errors (right column).}
	\label{fig:eigenfunctions_2D.esp2}
\end{figure}

Table \ref{NE:E1_square_eig} presents the numerical results and the associated errors. As seen from the table, the proposed FiedTNN-based method yields highly accurate results. In addition, the values of 
$\left|\mathbf{E}_{\scriptscriptstyle \text{NN}}\right|_{\mathbf{H}(\div ; \Omega)}$ 
listed in Table \ref{NE:E1_square_eig} confirm that the eigenfunctions are divergence-free. 
For comparison, Table \ref{NE:E1_comparison} provides the numerical results obtained by using the two-grid edge element method \cite{SINUM_2gird2014}, 
demonstrating that our method outperforms it in terms of numerical precision. 
The corresponding eigenfunctions are illustrated in Figures \ref{fig:eigenfunctions_2D.esp1} - \ref{fig:eigenfunctions_2D.esp2}, further showing that the numerical solutions obtained by the FieldTNN-based machine learning method align well with the reference solutions and exhibit satisfactory accuracy.
\begin{table}[!ht]
\centering
\caption{The $18$ smallest Maxwell eigenvalues in the square cavity $\Omega = [0,1]^{2}$.}\label{NE:E1_square_eig}%
\begin{tabular}{cccc}
\hline
Exact eigenvalues &Approximation eigenvalues &$\text{err}_{\lambda}$ &$\left|\mathbf{E}_{\scriptscriptstyle \text{NN}}\right|_{\mathbf{H}(\div ; \Omega)}$\\
\hline     
9.869604401089~($\pi^{2}$)  &9.86960453843 &1.39e-08 &5.54e-10 \\
9.869604401089~($\pi^{2}$)  &9.86960465907 &2.61e-08 &1.15e-09 \\
19.73920880218~($2\pi^{2}$) &19.7392094187 &3.12e-08 &9.15e-10 \\
39.47841760436~($4\pi^{2}$) &39.4784180287 &1.07e-08 &4.39e-10 \\
39.47841760436~($4\pi^{2}$) &39.4784182384 &1.61e-08 &5.23e-10 \\
49.34802200545~($5\pi^{2}$) &49.3480223816 &7.62e-09 &1.19e-09  \\
49.34802200545~($5\pi^{2}$) &49.3480225844 &1.17e-08 &1.03e-09  \\
78.95683520871~($8\pi^{2}$) &78.9568361547 &1.20e-08 &9.77e-10 \\
88.82643960980~($9\pi^{2}$) &88.8264397765 &1.88e-09 &3.47e-10 \\
88.82643960980~($9\pi^{2}$) &88.8264400261 &4.69e-09 &6.77e-10 \\
98.69604401089~($10\pi^{2}$)&98.6960449485 &9.50e-09 &7.14e-10 \\
98.69604401089~($10\pi^{2}$)&98.6960445521 &5.48e-09 &1.12e-09 \\
128.3048572142~($13\pi^{2}$)&128.304858492 &9.96e-09 &1.10e-09 \\
128.3048572142~($13\pi^{2}$)&128.304858232 &7.93e-09 &1.48e-09 \\
157.9136704174~($16\pi^{2}$)&157.913670765 &2.20e-09 &5.70e-10 \\
157.9136704174~($16\pi^{2}$)&157.913670865 &2.84e-09 &6.58e-10 \\
167.7832748185~($17\pi^{2}$)&167.783275736 &5.47e-09 &7.87e-10  \\
167.7832748185~($17\pi^{2}$)&167.783276446 &9.70e-09 &1.17e-09 \\
\hline
\end{tabular}
\end{table}	

\begin{table}[!ht]
\centering
\caption{Comparison result computed by two-grid methods, square cavity (cf. \cite{SINUM_2gird2014}, Tables $1$-$2$), \\ $H$ and $h$ denote the mesh size of the coarse grid and the fine grid, respectively.}\label{NE:E1_comparison}
\begin{tabular}{ccccc}
\hline
$k$ & $H$& $h$& $\lambda_{k}^{h}$ &$|\lambda_{k}-\lambda_{k}^{h}|/\lambda_{k}$ \\
\hline     
1 &1/16 &1/256 &9.869597  &7.50e-07  \\
2 &1/16 &1/256 &9.869597  &7.50e-07  \\
3 &1/16 &1/256 &19.739291 &4.16e-06  \\
\hline
1 &1/64 &1/512 &9.869585  &1.97e-06  \\
2 &1/64 &1/512 &9.869602  &2.43e-07  \\
3 &1/64 &1/512 &19.739229 &1.02e-06  \\
\hline
\end{tabular}
\end{table}	

\subsection{2D L-shaped cavity}
In the second example, we consider the computational domain to be a 2D L-shaped cavity $\Omega = [-1,1]^{2}/\left\{[0,1]\times[-1,0]\right\}$, as shown in Figure \ref{fig:sec4} $(1)$, and suppose that $\varepsilon = \mu = 1$ in $\Omega$.
This problem is part of Dauge's benchmark computations for the Maxwell eigenvalue problem \cite{Dauge_Bechmark}.
As mentioned in \cite{Dauge_Bechmark,JSC2020_Maxwell}, the regularity of the Maxwell eigenfunctions associated with the leading $5$ eigenvalues is described as follows. 
The $1$-st and $5$-th eigenfunctions have a strong unbounded singularity, which belongs to $\left(H^{2 / 3-\epsilon}(\Omega)\right)^2$ for any $\epsilon>0$. 
Specifically, the first eigenfunction around the $(0,0)$ is proportional to $r^{-\frac{1}{3}}$, where $r$ is the radius in polar coordinate.
While the $2$-nd one belongs to $\left(H^{4 / 3-\epsilon}(\Omega)\right)^2$. 
The $3$-rd and $4$-th ones are analytic, with the exact eigenvalue being $\pi^2 \approx 9.869604401$.
Each subnetwork of the FiledTNN is an FNN with  $3$ hidden layers and each hidden layer has $50$ neurons.
The activation function is selected as sine function and the rank $p$ is set to be $100$.
The Adam optimizer is employed with a learning rate $1\text{e-}5$ and $100000$ epochs to produce the final result. 
The Legendre-Gauss quadrature scheme with $1600$ points is adopted in each dimension.

We present in Table \ref{NE:E2_lsd_eig} the numerical results and the corresponding errors, 
and list in Table \ref{NE:E2_comparison} the numerical results computed by finite element methods \cite{CAMA_2023} and mixed spectral-element methods \cite{CICP_2015} as well the relative errors. 
From Table \ref{NE:E2_lsd_eig}, we find that the proposed FiedTNN-based machine learning method not only computes eigenvalues accurately but also ensures the divergence-free condition of the eigenfunctions.

\begin{figure}[H]
\centering
\includegraphics[width=.3\textwidth]{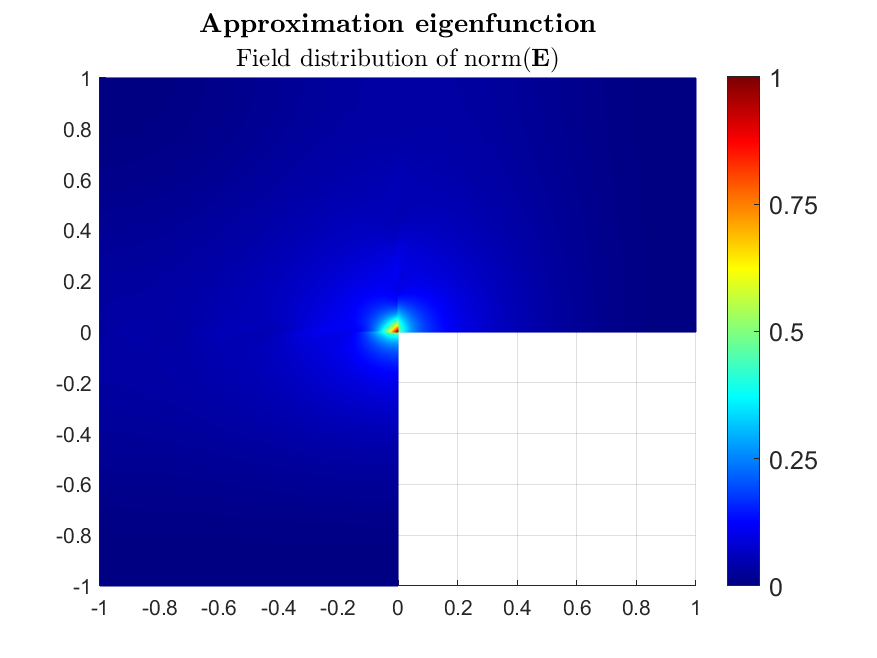}~
\includegraphics[width=.3\textwidth]{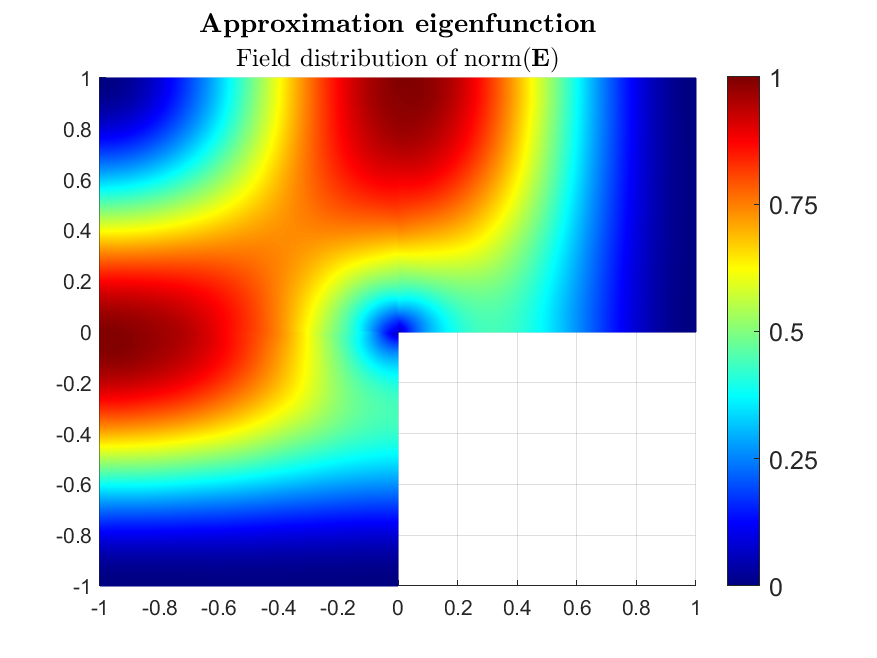}~
\includegraphics[width=.3\textwidth]{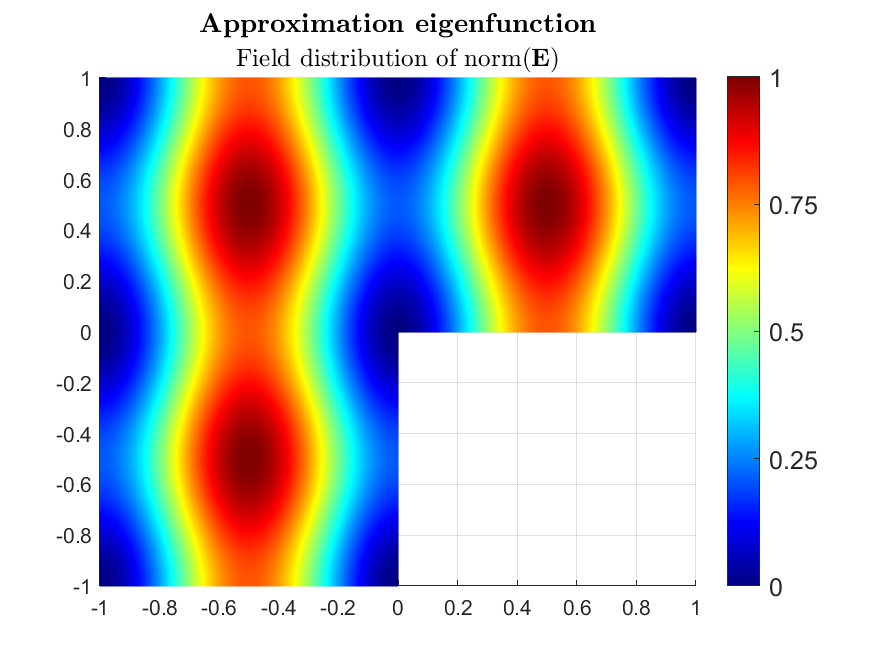}\\
\includegraphics[width=.3\textwidth]{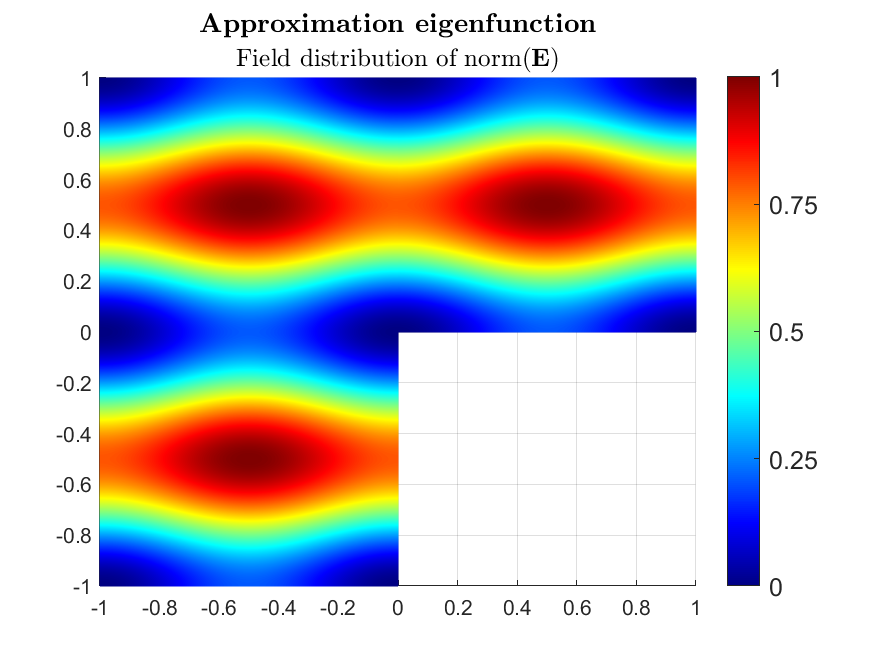}~
\includegraphics[width=.3\textwidth]{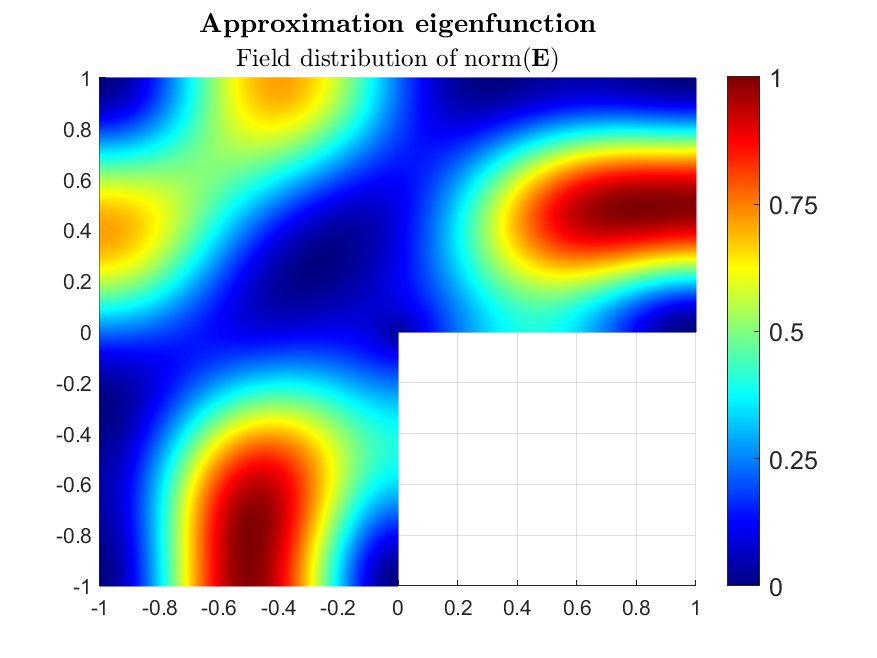}\\			
\caption{Norms of Maxwell eigenfunctions for the first $5$ eigenvalues in a 2D L-shaped cavity, arranged from top left to bottom right.}
\label{fig:eigf_LSD.esp}
\end{figure}

In addition, we also provide in Figure \ref{fig:eigf_LSD.esp} the field distribution of the leading $5$ Maxwell eigenfunctions $\big($norm$(\mathbf{E}_{\scriptscriptstyle \text{NN}})\big)$ in the L-shaped domain. 
It is observed from Figure \ref{fig:eigf_LSD.esp} that the norm of $1$st eigenfunctions is considerably large around the $(0,0)$ and the $3$-rd and the $4$-th ones are analytic, these observations match the conclusion we mention above.

\begin{table}[!ht]
\centering
\caption{The $5$ smallest Maxwell eigenvalues in 2D L-shaped cavity.}\label{NE:E2_lsd_eig} 
\begin{tabular}{cccc}
\hline
Benchmark &Approximation eigenvalues &$\text{err}_{\lambda}$ &$\left|\mathbf{E}_{\scriptscriptstyle \text{NN}}\right|_{\mathbf{H}(\div ; \Omega)}$\\
\hline     
1.47562182408 &1.47367812406  &1.32e-03 &2.46e-03\\
3.53403136678 &3.53401751530  &3.92e-06 &7.64e-04\\
9.86960440109 &9.86960441721  &1.63e-09 &1.33e-06\\
9.86960440109 &9.86960439972  &1.39e-10 &2.05e-06\\
11.3894793979 &11.38945178126 &2.42e-06 &1.16e-04\\
\hline
\end{tabular}
\end{table}	



\begin{table}[!ht]
\centering
\caption{Comparison results computed by FEM and mixed  SEM, $2$D L-shaped cavity.}\label{NE:E2_comparison}
\begin{tabular}{cccc}
\hline
~&Uniform mesh, $1/64$  &Nonuniform mesh, $1/256$  &Mixed SEM \\
No.&(cf. \cite{CAMA_2023}, Table $4$) &(cf. \cite{CAMA_2023}, Table $5$) &(cf. \cite{CICP_2015}, Table 3)\\
\hline     
1&1.49509\ (1.32e-02)  &1.48561\ (6.77e-03) &1.4752\ (3.20e-04)  \\
2&3.53609\ (5.83e-04)  &3.53414\ (3.07e-05) &3.5340\ (4.86e-07)  \\
3&9.88218\ (1.27e-03)  &9.87008\ (4.82e-05) &9.8697\ (1.41e-05)  \\
4&9.88218\ (1.27e-03)  &9.87008\ (4.82e-05) &9.8698\ (1.86e-05)  \\
5&11.40879\ (1.70e-03) &11.39012\ (5.62e-05) &11.3899\ (3.72e-05) \\
\hline
\end{tabular}
\end{table}	

\subsection{Inhomogeneous cavity} \label{eq:num_inhomo}
\begin{figure}[H]
\centering
\subfloat[{Inhomogeneous computational domain.}]{\includegraphics[width=0.4\linewidth]{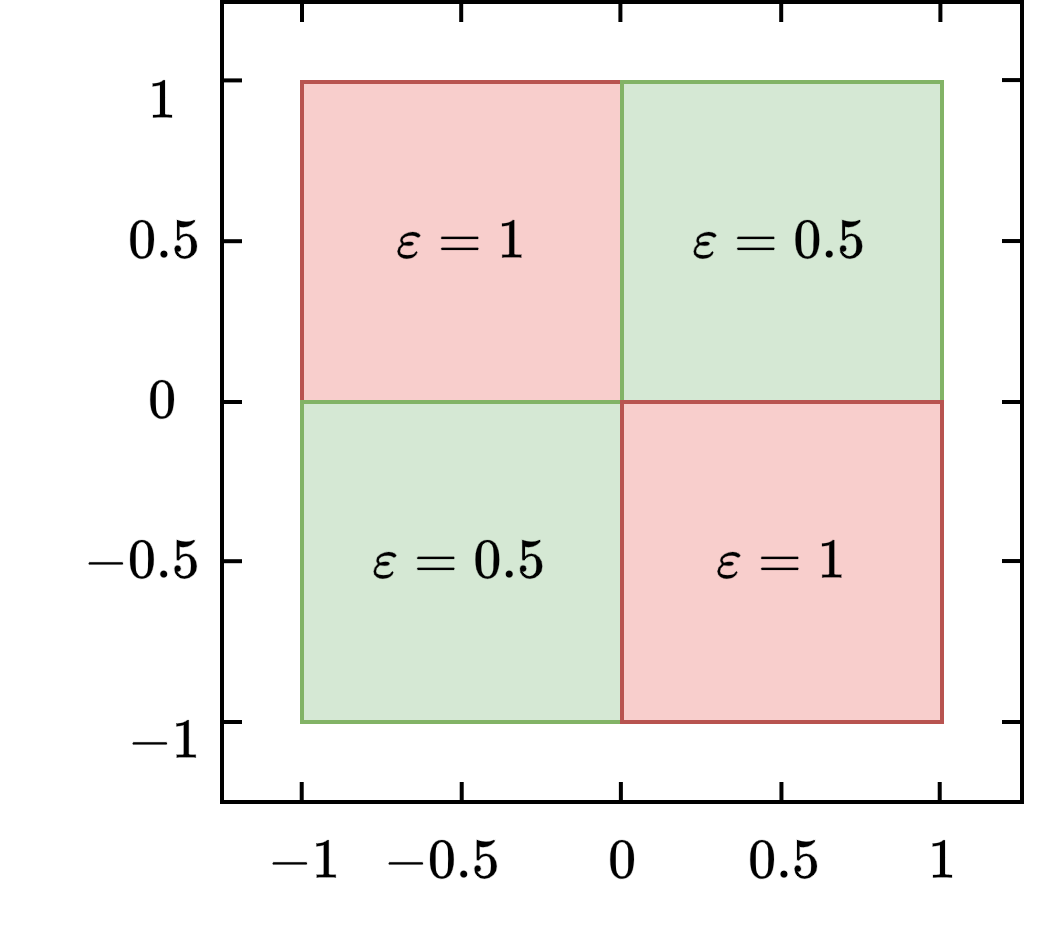}}\hfill 
\subfloat[{Decomposition of inhomogeneous cavity.}]{\includegraphics[width=0.4\linewidth]{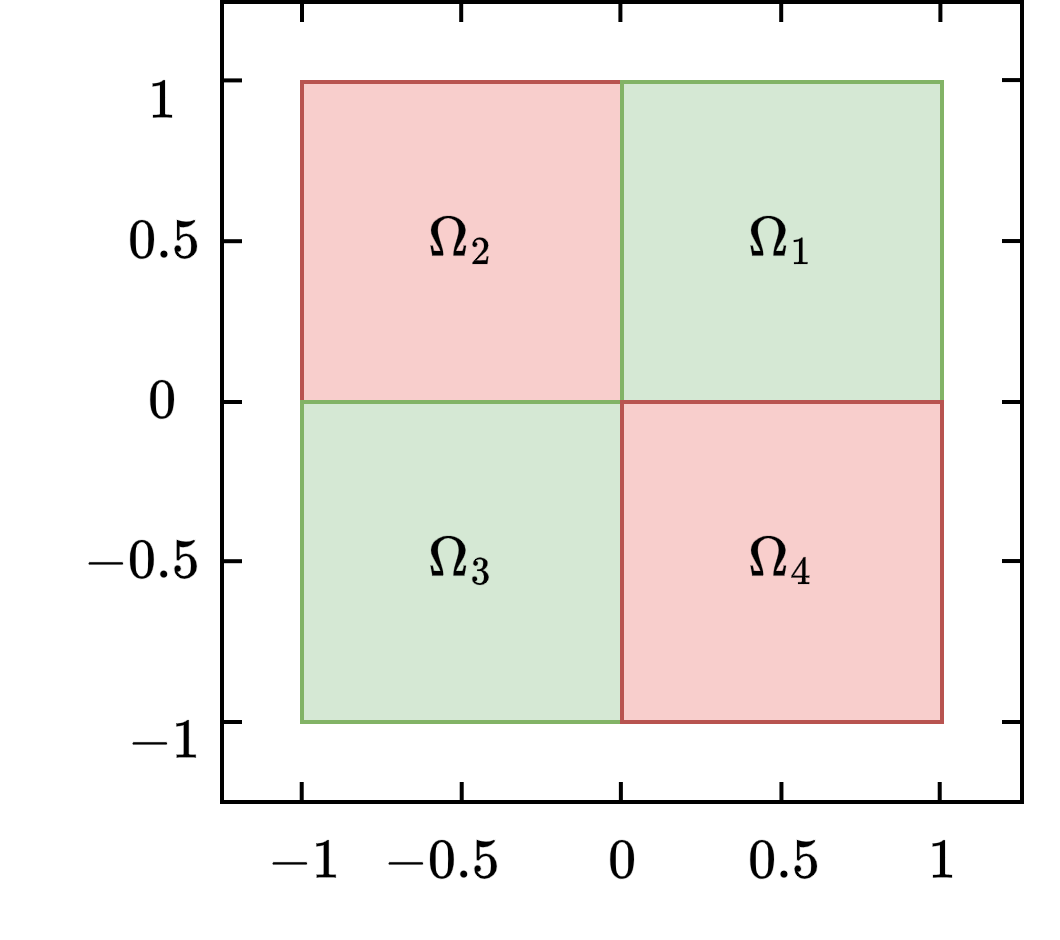}}\\
\caption{Computational domain for the numerical example in Section $\ref{eq:num_inhomo}$.}
\label{fig:Omega}
\end{figure}

In the third example, we examine an inhomogeneous medium within the computational domain $\Omega = [-1,1]^{2}$. 
As illustrated in Figure \ref{fig:Omega}, the cavity is divided into four squares. 
The squares $\Omega_{1}$ and $\Omega_{3}$ are composed of material $1$ with $\varepsilon_{1} = 0.5$ and $\mu_{1} = 1$, while the remaining squares are constituted of material $2$ with $\varepsilon_{2} = 1$ and $\mu_{2} = 1$. 
Here, $\Omega_{1} = [0,1]^2$, $\Omega_{2} = [-1,0] \times [0,1]$, $\Omega_{3} = [-1,0]^2$ and $\Omega_{4} = [0,1] \times [-1,0]$.

Each subnetwork of the FiledTNN is an FNN with  $3$ hidden layers and each hidden layer has $50$ neurons.
The activation function is selected as sine function and the rank $p$ is set to be $100$.
The Adam optimizer is employed with a learning rate $3\text{e-}4$ and $50000$ epochs to produce the final result. 
The Legendre-Gauss quadrature scheme with $1600$ points is adopted in each dimension.

We present in Table \ref{NE:E3_inhomo_eig} the numerical results and the corresponding errors, and list in Table \ref{NE:E3_comparison} the comparison results computed by the mixed spectral-element methods \cite{CICP_2015}. 
We observe from Table \ref{NE:E3_inhomo_eig} that the proposed machine learning method demonstrates good accuracy across the $10$ smallest approximation Maxwell eigenvalues in the inhomogeneous cavity. 
The errors in approximation eigenvalues are predominantly very small, suggesting the robustness and precision of the method.
The divergence of most eigenfunctions are also minimal, indicating stable numerical performance.
The comparison with mixed spectral element methods shows that the our method is competitive and achieves similar levels of accuracy. 
Further refinement and optimization of the method could potentially reduce the errors even more, making it a robust choice for solving Maxwell eigenvalue problems in inhomogeneous media.

\begin{figure}[htbp]
\centering
\includegraphics[width=.3\textwidth]{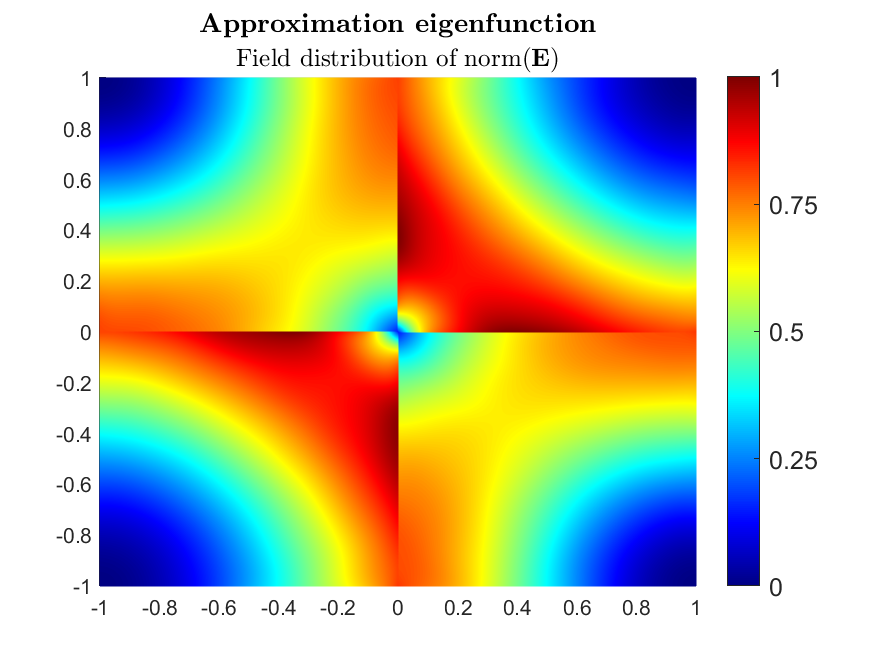}~
\includegraphics[width=.3\textwidth]{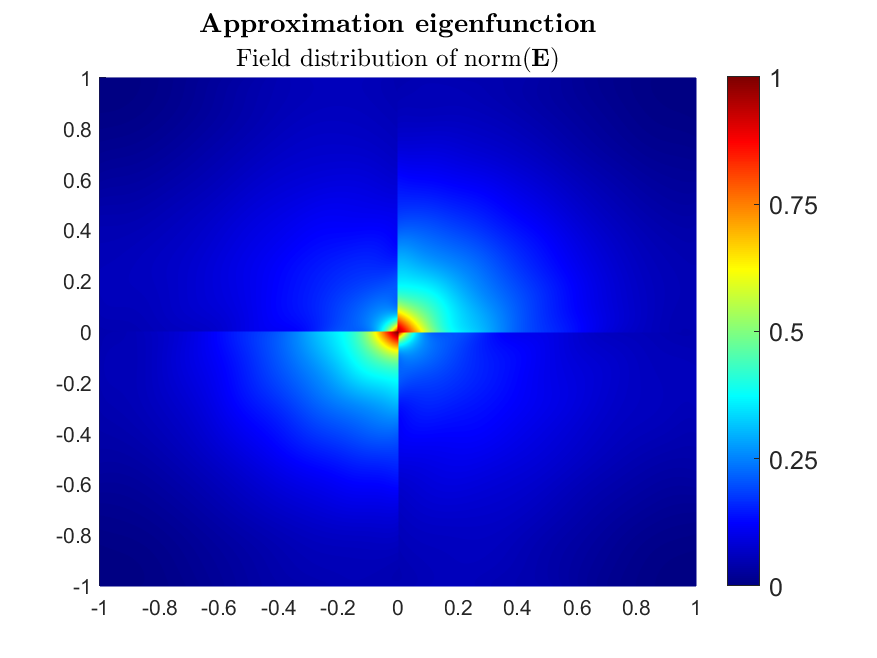}~
\includegraphics[width=.3\textwidth]{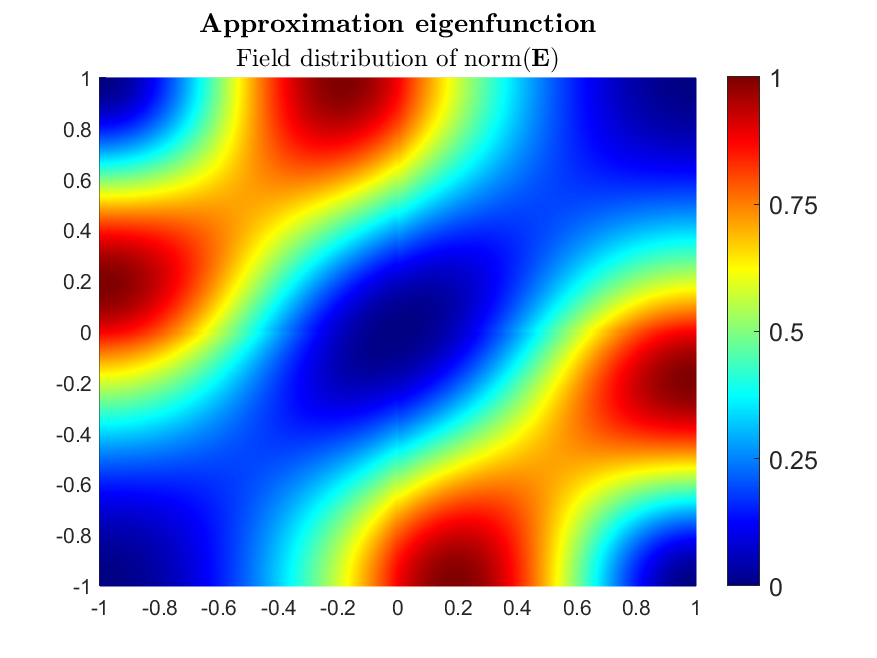}\\
\includegraphics[width=.3\textwidth]{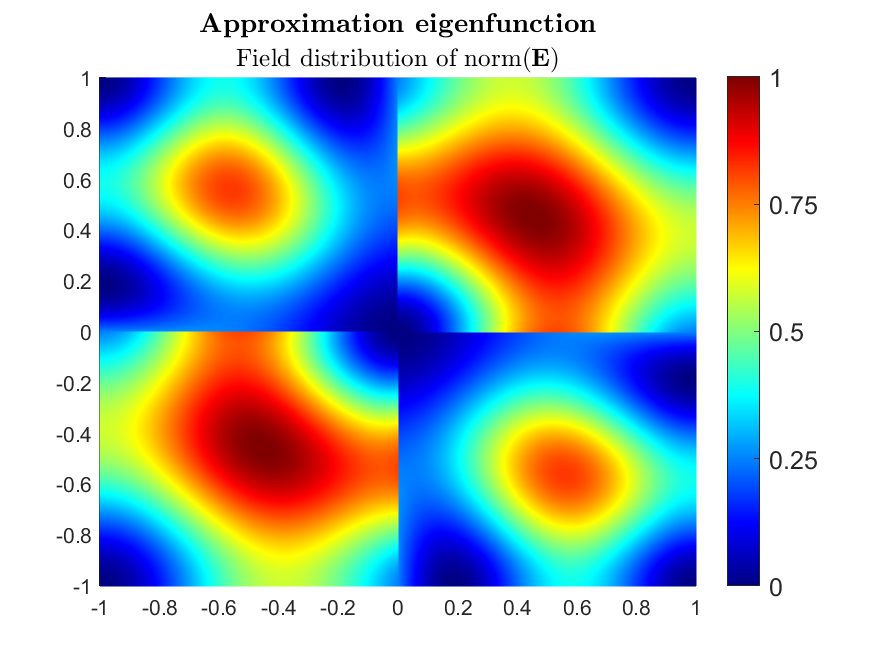}~
\includegraphics[width=.3\textwidth]{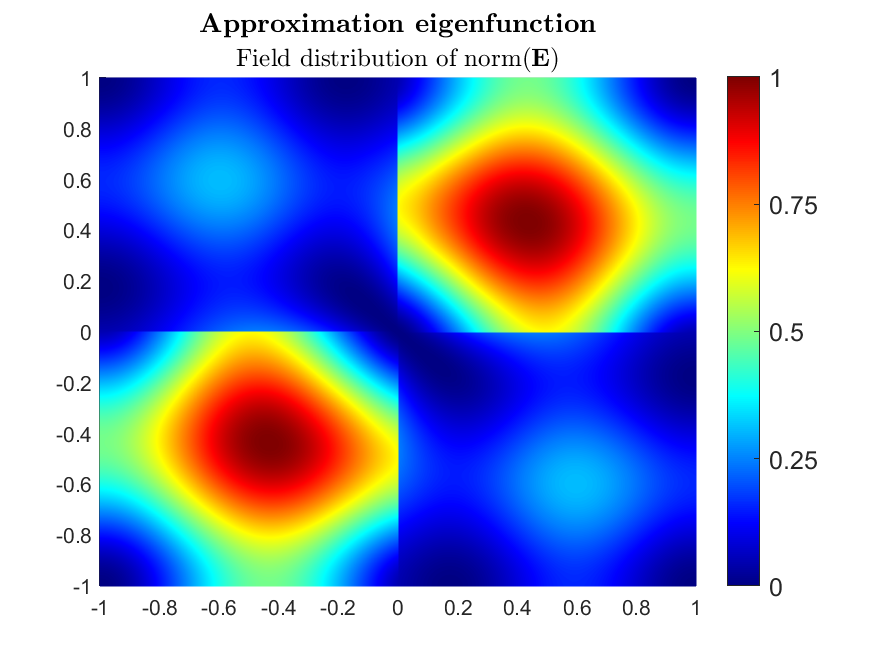}~
\includegraphics[width=.3\textwidth]{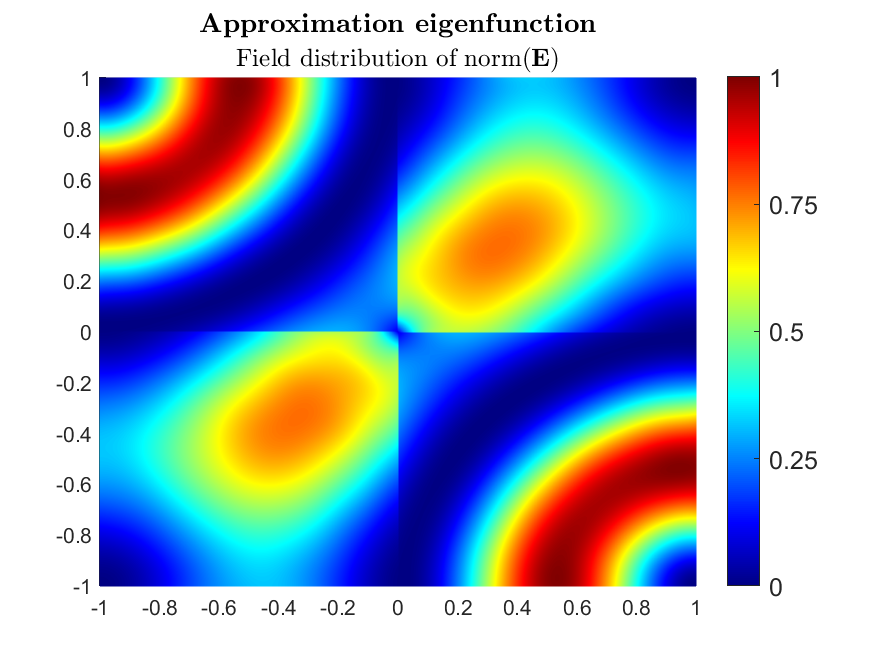}\\			
\caption{Norms of Maxwell eigenfunctions corresponding to the first $6$ eigenvalues in an inhomogeneous cavity, arranged from top left to bottom right.}
\label{fig:eigf_trans.esp}
\end{figure}

In addition, we also provide the field distribution of the leading $6$ Maxwell eigenfunctions $\big($norm$(\mathbf{E}_{\scriptscriptstyle \text{NN}})\big)$ in inhomogeneous cavity in Figure \ref{fig:eigf_trans.esp}, which illustrates that the norms of eigenfunctions are symmetric with respect to subdomains $\Omega_{1}$ and $\Omega_{3}$, as well as  subdomains $\Omega_{2}$ and $\Omega_{4}$, respectively. 

\begin{table}[!ht]
\centering
\caption{The $10$ smallest Maxwell eigenvalues in inhomogeneous cavity.}\label{NE:E3_inhomo_eig}
\begin{tabular}{cccc}
\hline
Benchmark &Approximation eigenvalues &$\text{err}_{\lambda}$ &$\left|\mathbf{E}_{\scriptscriptstyle \text{NN}}\right|_{\mathbf{H}(\div ; \Omega)}$\\
\hline     
3.317548763415 &3.31753973639  &2.72e-06 &4.51e-03\\
3.366324157260 &3.36583043172  &1.47e-04 &3.98e-02\\
6.186389562488 &6.1863895391   &3.78e-09 &7.72e-06\\
13.92632333103 &13.92632331809 &9.29e-10 &6.08e-06\\
15.08299096123 &15.08299088572 &5.01e-09 &1.29e-04\\
15.77886590819 &15.77882697523 &2.47e-06 &1.44e-03\\
18.64329693686 &18.64318258832 &6.13e-06 &1.87e-03\\
25.79753111031 &25.79752850409 &1.01e-07 &6.82e-05\\
29.85240067684 &29.85234879779 &1.73e-06 &6.15e-04\\
30.53785871253 &30.53522923416 &8.61e-05 &1.11e-02\\
\hline
\end{tabular}
\end{table}	

\begin{table}[!ht]
\centering
\caption{Comparison results computed by mixed SEM with $N=6$, inhomogeneous cavity (cf. \cite{CICP_2015}).}\label{NE:E3_comparison}%
\begin{tabular}{ccc}
\hline
$k$&$\lambda_{k}$  &relative error\\
\hline     
1&3.31754854&6.73e-08 \\
2&3.36628055&1.30e-05 \\
3&6.18638956&4.02e-10 \\
4&13.92632333&7.40e-11 \\
5&15.08299096&8.15e-11 \\
\hline
\end{tabular}
\end{table}	

\subsection{Cube cavity}
Next, we consider the Maxwell eigenvalue problem in a cubic cavity $\Omega = [0,1]^{3}$, and suppose that $\varepsilon = \mu = 1$ throughout the domain. 
The exact eigenvalues are known to take the form $\left\{k \pi^{2}\right\}$, where $k=k_{1}^{2}+k_{2}^{2}+k_{3}^{2}$ are non-negative integers satisfying $k_{1} k_{2}+k_{2} k_{3}+k_{3} k_{1}>0$. 
For fixed $k_{1},~k_{2}~k_{3}$, the corresponding eigenfunction is given by the following vector function:
\begin{align*}
\mathbf{E} = 
\left(\begin{array}{l}
\cos \left(k_1 \pi x_1\right) \sin \left(k_2 \pi x_2\right) \sin \left(k_3 \pi x_3\right) \\
\sin \left(k_1 \pi x_1\right) \cos \left(k_2 \pi x_2\right) \sin \left(k_3 \pi x_3\right) \\
\sin \left(k_1 \pi x_1\right) \sin \left(k_2 \pi x_2\right) \cos \left(k_3 \pi x_3\right)
\end{array}\right).	
\end{align*}

Each subnetwork of the FiledTNN is an FNN with  $3$ hidden layers and each hidden layer has $100$ neurons.
The activation function is selected as sine function and the rank $p$ is set to be $50$.
The Adam optimizer is used with an initial learning rate of $0.0003$ for the first $100000$ epochs, followed by L-BFGS optimization for the last $5000$ steps with a learning rate of $1$ to obtain the final result. 
A Legendre-Gauss quadrature scheme with $1600$ points for each dimension is utilized for numerical integration.

Table \ref{NE:E4_cube_eig} presents the $17$ smallest Maxwell eigenvalues and their corresponding relative errors, 
demonstrating the high accuracy of the FieldTNN-based machine learning method in solving the Maxwell eigenproblems within the cubic cavity. 
It is worth noting that our approach can accurately captures the multiplicity of eigenvalues, 
providing multiple approximations, all with errors of the same order of magnitude. 
Additionally, the eigenfunctions are confirmed to be divergence-free. 
For comparison, Table \ref{NE:E4_comparison} lists the results computed by using the two-grid method \cite{SINUM_2gird2014}, 
which indicates that our method achieves higher precision in the approximation eigenvalues.
\begin{table}[!ht]
\centering
\caption{The $17$ smallest Maxwell eigenvalues in the cube cavity $\Omega = [0,1]^{3}$.}\label{NE:E4_cube_eig}
\begin{tabular}{cccc}
\hline
Exact eigenvalues &Approximation eigenvalues &$\text{err}_{\lambda}$ &$\left|\mathbf{E}_{\scriptscriptstyle \text{NN}}\right|_{\mathbf{H}(\div ; \Omega)}$\\
\hline   
19.73920880218($2\pi^{2}$) & 19.7392088182 & 8.05e-10 & 1.67e-10 \\  
19.73920880218($2\pi^{2}$) & 19.7392088446 & 2.15e-09 & 2.07e-10 \\ 
19.73920880218($2\pi^{2}$) & 19.7392089280 & 6.37e-09 & 8.99e-11 \\  
29.60881320327($3\pi^{2}$) & 29.6088133301 & 4.28e-09 & 1.73e-09 \\ 
29.60881320327($3\pi^{2}$) & 29.6088134093 & 6.96e-09 & 1.50e-09 \\  
49.34802200545($5\pi^{2}$) & 49.3480220798 & 1.51e-09 & 2.53e-10 \\ 
49.34802200545($5\pi^{2}$) & 49.3480221351 & 2.63e-09 & 2.00e-10 \\  
49.34802200545($5\pi^{2}$) & 49.3480221031 & 1.98e-09 & 3.65e-10 \\  
49.34802200545($5\pi^{2}$) & 49.3480221510 & 2.95e-09 & 2.02e-10 \\  
49.34802200545($5\pi^{2}$) & 49.3480222062 & 4.07e-09 & 2.14e-10 \\  
49.34802200545($5\pi^{2}$) & 49.3480222598 & 5.15e-09 & 2.83e-10 \\  
59.21762640654($6\pi^{2}$) & 59.2176266206 & 3.61e-09 & 1.44e-09 \\  
59.21762640654($6\pi^{2}$) & 59.2176265619 & 2.62e-09 & 2.12e-09 \\ 
59.21762640654($6\pi^{2}$) & 59.2176266781 & 4.59e-09 & 2.34e-09 \\ 
59.21762640654($6\pi^{2}$) & 59.2176266090 & 3.42e-09 & 2.71e-09 \\ 
59.21762640654($6\pi^{2}$) & 59.2176267312 & 5.48e-09 & 3.20e-09 \\
59.21762640654($6\pi^{2}$) & 59.2176265419 & 2.29e-09 & 5.55e-09 \\                      
\hline
\end{tabular}
\end{table}	

\begin{table}[!ht]
\centering
\caption{Comparison result computed by two-grid methods, cube cavity (cf. \cite{SINUM_2gird2014}, Tables $5$), \\ $H$ and $h$ denote the mesh size of the coarse grid and the fine grid, respectively.}\label{NE:E4_comparison}
\begin{tabular}{ccccc}
\hline
$k$ & $H$& $h$& $\lambda_{k}^{h}$ &$|\lambda_{k}-\lambda_{k}^{h}|/\lambda_{k}^{h}$ \\
\hline     
1 &1/2 &1/8 &19.467320 &1.38e-02  \\
2 &1/2 &1/8 &19.693282 &2.33e-03  \\
3 &1/2 &1/8 &19.693283 &2.33e-03  \\
\hline
1 &1/4 &1/64 &19.734459 &2.41e-04  \\
2 &1/4 &1/64 &19.738345 &4.38e-05  \\
3 &1/4 &1/64 &19.738345 &4.38e-05  \\
\hline
\end{tabular}
\end{table}

\subsection{3D L-shaped cavity}
\begin{figure}[H]
\centering
\includegraphics[width=0.4\linewidth]{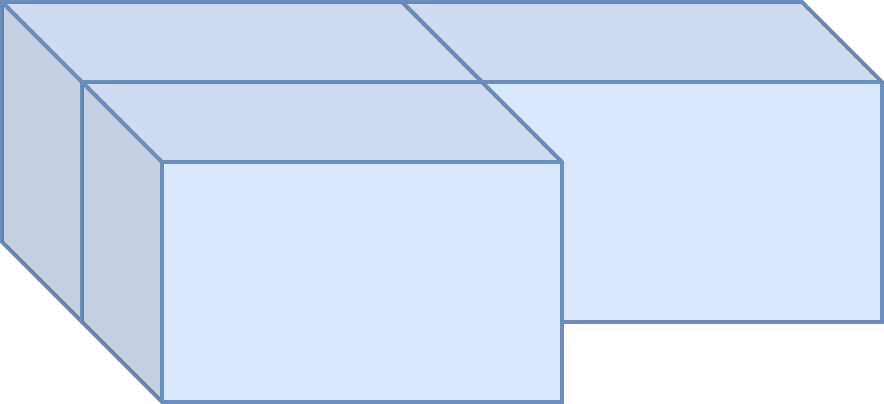}
\caption{Computational domain: 3D L-shaped cavity $\Omega = \{[-1,1]^{2}/[0,1]\times[-1,0]\}\times [0,1]$.}
\label{fig:Maxwell_3Dlsd}
\end{figure}

Finally, we calculate the Maxwell eigenvalues in a 3D L-shaped cavity $\Omega = \{[-1,1]^{2}/[0,1]\times[-1,0]\}\times [0,1]$, as shown in Figure \ref{fig:Maxwell_3Dlsd}. 
This problem is also part of Dauge's benchmark computations for the Maxwell eigenvalue problem \cite{Dauge_Bechmark}. 
The nonconvexity of the outer boundary causes some eigenfunctions to exhibit significant singularities, making it a challenging test as noted in \cite{IEEE_2013}.

Each subnetwork of the FiledTNN is an FNN with  $3$ hidden layers and each hidden layer has $50$ neurons.
The activation function is selected as sine function and the rank $p$ is set to be $100$.
The Adam optimizer is employed with a learning rate $3\text{e-}4$ and $50000$ epochs to produce the final result. 
The Legendre-Gauss quadrature scheme with $800$ points is adopted in each dimension.

The numerical results computed by FieldTNN-based machine learning method are provided in Table \ref{NE:E5_LSD_eig} and we observe from Table \ref{NE:E5_LSD_eig} that the numerical results are all match with the benchmark ones, which demonstrates that the proposed algorithm are efficient and robust. 
We also list the numerical eigenvalue in Table \ref{NE:E5_comparison} as a comparison, we find that FieldTNN-based machine learning method performs better than the mixed finite element method \cite{BoffiMixed}, and our method can achieve the same level of precision with the mixed spectral element method. 
Moreover, we should mention that, our method can achieve higher accurate numerical result for the case that eigenfunction with high regularity.
	
\begin{table}[!ht]
\centering
\caption{The $9$ smallest Maxwell eigenvalues in 3D L-shaped cavity.}\label{NE:E5_LSD_eig} 
\begin{tabular}{cccc}
\hline
Benchmark &Approximation eigenvalues &$\text{err}_{\lambda}$ &$\left|\mathbf{E}_{\scriptscriptstyle \text{NN}}\right|_{\mathbf{H}(\div ; \Omega)}$\\
\hline     
9.63972384472&9.66976383085   &3.12e-03 &2.56e-04\\
11.3452262252&11.33568261824  &8.41e-04 &1.35e-02\\
13.4036357679&13.40342600327  &1.57e-05 &1.99e-03\\
15.1972519265&15.19864892752  &9.19e-05 &1.20e-04\\
19.5093282458&19.54229680438  &1.69e-03 &4.29e-03\\
19.7392088022&19.73922077213  &3.99e-07 &2.75e-05\\
19.7392088022&19.73923061875  &6.06e-07 &3.03e-05\\
19.7392088022&19.73921668788  &1.11e-06 &3.94e-05\\
21.2590837990&21.25881721428  &1.25e-05 &4.50e-04\\
\hline
\end{tabular}
\end{table}	

\begin{table}[!ht]
\centering
\caption{Comparison result computed by FEM and mixed SEM, 3D L-shaped cavity.}\label{NE:E5_comparison}%
\begin{tabular}{cccc}
\hline
Mixed FEM (cf. \cite{CMAME2008})& err. &Mixed SEM (cf. \cite{IEEE_2013}) &err.\\
\hline     
9.6627&2.38e-3&9.6616370&2.27e-3  \\
10.5974&6.59e-2&11.3500164&4.22e-4 \\
13.3150&6.61e-3&13.4036732&2.79e-6 \\
13.6238&1.04e-1&15.1976786&2.81e-5 \\
15.9957&1.80e-1&19.5399025&1.57e-3 \\
16.7727&1.50e-1&19.7393602&7.67e-6  \\
18.6795&5.37e-2&19.7393602&7.67e-6  \\
19.7434&2.12e-4&19.7393602&7.67e-6  \\
19.7643&7.03e-2&-&- \\
\hline
\end{tabular}
\end{table}

\section{Conclusions}\label{sec:conclusion}
In this paper, we introduce FieldTNN, a novel network structure based on the TNN framework, designed to represent vector functions and solve electromagnetic field equations, with a particular focus on the Maxwell eigenvalue problem.
The contributions of this work can be summarized in the following three aspects.
Firstly, we have applied the TNN-based machine learning method to vector field equations for the first time. 
Secondly, we successfully extended the computational domain of TNN from tensor-type to non-tensor-type ones, significantly broadening the applicability of TNN to more complex regions. 
Lastly, to address the divergence-free condition in Maxwell eigenvalue problems, we designed a filtering algorithm that automatically eliminates eigenpairs that do not satisfy the divergence-free condition during the machine learning process.

Indeed, the proposed FieldTNN-based machine learning method can be developed to other type vector field equations, such as Maxwell equations, Stokes equations, magnetohydrodynamics equations and so on. Moreover, although our algorithm has extended the computational domain to non-tensor-type regions, it still faces certain limitations.
How to extend the TNN-based machine learning method proposed in this paper to the more complex computational domain is one of our future research.

\section*{Acknowledge}
\noindent This work is supported by the Strategic Priority Research Program of the Chinese Academy of Sciences (Grant No. XDB0640000), National Natural Science Foundation of China (Grant No. 12171026, 1233000214, U2230402 and 12031013), President Foundation of China Academy of Engineering Physics (YZJJ2Q2022017), the China Postdoctoral Science Foundation (GZB20240025), National Key Laboratory of Computational Physics (No.6142A05230501), National Center for Mathematics and Interdisciplinary Science, CAS. 

\section*{Conflict of interest statement} 
\noindent The authors declare that they have no competing interests.

\bibliographystyle{plain}
\bibliography{article.bib}

\end{document}